\documentclass{article}

\title{Mixed-Integer Program Sensitivity Analysis via Branch-and-Bound Trees}

\usepackage[affil-it]{authblk}

\author[1]{Anna Deza}
\author[2]{Fabian Badilla}
\author[2]{Santanu S. Dey}
\affil[1]{Mitsubishi Electric Research Laboratories (MERL), USA}
\affil[2]{H. Milton Stewart School of Industrial and Systems Engineering, Georgia Institute of Technology, USA} 
\date{}
\usepackage[utf8]{inputenc}
\usepackage{cancel}
\usepackage[T1]{fontenc}
\usepackage{algorithm}
\usepackage{algpseudocode}
\usepackage[margin=1in]{geometry}
\usepackage{amsmath, amssymb, amsthm}
\usepackage{graphicx}
\usepackage{booktabs}           
\usepackage[colorlinks,allcolors=blue]{hyperref}           
\usepackage[capitalize,nameinlink]{cleveref}

\usepackage{xcolor}             
\usepackage{enumitem}           
\usepackage{bm}
\usepackage{natbib}
\usepackage{tikz}
\usepackage{array}
\usepackage{bigstrut}
\usepackage[affil-it]{authblk}

\usepackage{pgfplots}
\pgfplotsset{compat=1.18}
\usepgfplotslibrary{groupplots}

\hypersetup{
  colorlinks=true
}

\newtheorem{theorem}{Theorem}
\newtheorem{lemma}{Lemma}

\newtheorem{proposition}{Proposition}
\theoremstyle{definition}
\newtheorem{definition}{Definition}

\newtheorem{remark}{Remark}

\newcommand{\R}{\mathbb{R}}
\newcommand{\Z}{\mathbb{Z}}

\newcommand{\TT}{\mathcal{T}}
\newcommand{\LL}{\mathbb{L}}

\newcommand{\val}[1][]{%
  \if\relax\detokenize{#1}\relax
    z
  \else
    z({#1})
  \fi
}

\newcommand\nominal{0}
\newcommand{\pert}[1][]{%
  \if\relax\detokenize{#1}\relax
    \Delta
  \else
    \Delta_{#1}
  \fi
}
\newcommand{\dualvar}{\lambda}
\newcommand{\dual}[1][]{%
  \if\relax\detokenize{#1}\relax
    f
  \else
    f({#1})
  \fi
}
\newcommand{\lagr}{\mathcal{L}}
\newcommand{\boundnorm}{M}

\def\keywordname{{\bf Keywords:}}
\providecommand{\keywords}[1]{\def\and{{\textperiodcentered} }%
\par\addvspace\baselineskip
\noindent\keywordname\enspace\ignorespaces#1}%

\usepackage[font=small,labelfont=bf]{caption}

\begin{document}

\maketitle

\begin{abstract}
We study how to efficiently compute valid dual bounds for mixed-integer linear programs (MILPs) with {finite bounds on all variables} under perturbations to the problem data, in an offline-online setting where a nominal instance is solved once and the resulting information is reused to evaluate subsequent perturbations. We begin by developing a framework for computing valid dual bounds for linear programs (LPs) under simultaneous perturbations to the objective coefficients, right-hand side, and constraint matrix, extending classical LP sensitivity analysis beyond right-hand-side perturbations. Our construction is based on a bounded Lagrangian dual evaluated at nominal optimal dual multipliers, and we prove that the gap between the resulting dual bound and the perturbed LP optimal value grows at most linearly with the perturbation magnitude under suitable conditions. We then extend this construction to MILPs by reusing the branch-and-bound tree and node-level dual information from the nominal solve, generalizing the classical sensitivity analysis framework of Schrage and Wolsey (1985) to perturbations of all problem data. To strengthen the resulting dual bounds, we employ the subtree-selection dynamic program of Schrage and Wolsey, adapted to our generalized dual framework. A key feature of the proposed approach is that its online computations are inherently well suited for GPU acceleration, enabling highly parallel evaluation of dual bounds. We also show how these bounds can be used to certify optimal variable fixings and accelerate subsequent reoptimization. Extensive computational experiments on generalized independent set, combinatorial auction, capacitated facility location, and unit commitment instances demonstrate that the proposed method produces remarkably tight dual bounds at extremely low computational cost. For single-parameter perturbations arising in ``what-if'' analysis, the computed dual bounds typically have gaps of a small fraction of a percent, with GPU evaluation times well below a millisecond per perturbation. Compared with the time Gurobi requires to obtain dual bounds at least as strong as those produced by our method, the GPU approach is approximately 2,000--240,000 times faster. The method also produces strong dual bounds under simultaneous perturbations to multiple coefficients, with corresponding speedups of approximately 120--50,000 times over Gurobi. Moreover, in experiments using exact perturbed optima as supplied primal cutoffs,  the inferred variable fixings can be used to accelerate exact reoptimization, although the computational benefits vary considerably across problem classes, ranging from speedups of up to a factor of 10 for generalized independent set instances to little or no improvement for unit commitment.
    \keywords{Sensitivity Analysis \and Mixed-Integer Linear Programming \and GPU-native algorithm}
\end{abstract}
\section{Introduction}

Mixed-integer programs (MIPs) arise in many operational decision-making problems where solutions must be obtained under tight time constraints. As new data arrive, these problems often need to be reoptimized within minutes, posing significant computational challenges even for state-of-the-art solvers. For example, the Unit Commitment problem, used daily by Independent System Operators to clear day-ahead electricity markets, must be solved within 15–30 minutes due to the limited window between receiving bids and clearing the market \cite{Chen2016Feb}. Fast reoptimization is also needed for what-if analysis, where planners ask how a solution would change if one input were different. For example, the OptiGuide framework, deployed in the Microsoft cloud supply chain, answers each such query by re-solving the underlying MILP~\cite{li2023large}.

These applications typically exhibit two characteristics that create opportunities to improve computational efficiency by leveraging historical data. First, these problems are solved repeatedly, often at regular intervals (e.g., hourly or daily), or large repositories of historical instances are available. This provides opportunities for systematic data collection and offline analysis that can be exploited to accelerate future solves. Second, instances of the same operational problem typically share a common structure, with the same variables and constraint matrix. Variations arise primarily through perturbations to the right-hand sides (RHS) and objective coefficients, reflecting changes in demand, costs, or resource availability.  In some applications, the coefficients of the constraint matrix may also change, although the overall problem structure remains largely unchanged. 

Traditional MIP solution methods, which largely focus on solving individual instances independently, do not systematically exploit these two characteristics: the availability of historical instances and the structural similarity between them. We are interested in developing sensitivity techniques for MIPs to efficiently reoptimize them efficiently after small perturbations to their data. While significant progress has been made in designing primal heuristics that generate high-quality feasible solutions under such perturbations \cite{Shetty1991Jan,Xavier2020Oct,Bertsimas2022Jul,Johnson2020Mar}, there remains a substantial gap regarding the efficient computation of dual bounds under general perturbations. Given the critical role of dual bounds in assessing solution quality and guiding branch-and-bound algorithms, improving their computation after small data perturbations offers a key opportunity.

In this work, we consider the mixed-integer linear program parametrized by $\pert = (\pert[A], \pert[b], \pert[c])$

\begin{equation}\label{eqn:zDelta}
   \val(\pert) = \min_{x\in\R^n} \left\{ \left(c+\pert[c]\right)^\top x \mid  (A+\pert[A])x \geq b + \pert[b], \ \ell \leq x \leq u, \ x_i \in \Z~ \forall i \in \mathcal I\right\} \tag{P($\pert$)},
\end{equation}
where we assume $\ell$ and $u$ are finite. We denote by $z_{LP}(\pert)$ the optimal value of its LP relaxation, obtained by dropping the integrality constraints.

We study how the optimal objective value changes under small perturbations to the problem data.  Let $P(\nominal)$ denote the nominal problem, with optimal value $\val(\nominal)=\val$. Motivated by applications requiring rapid reoptimization
of related instances, we consider a two-phase framework. In an \emph{offline} phase, substantial computational effort may be invested in solving $P(\nominal)$, solving auxiliary problems, and storing information derived from these computations to build the scaffolding for solving $P(\pert)$. In the \emph{online} phase, previously unknown perturbations $\pert$ arrive, and the corresponding instances must be processed efficiently.

Our goal is to exploit information obtained in the offline phase to infer the optimal value $\val(\pert)$ for new perturbations. While recovering $\val(\pert)$ exactly may not be possible in general, we seek a valid dual bound $\dual(\pert)$ satisfying $$\dual(\pert) \leq \val(\pert).$$
For $\dual(\pert)$ to be useful, we seek to balance two quantities:
\begin{itemize}
    \item the \textit{strength} of the dual bound, measuring how close the dual is to the value function,
    \item the \textit{computational effort} required to evaluate $\dual(\pert)$ given a perturbation $\pert$.
\end{itemize}

\paragraph{Related papers.} 

Sensitivity analysis is a classical topic in linear programming (LP), particularly for perturbations to the right-hand-side (RHS) vector. Perturbations of the constraint matrix also have a long history. Mills~\cite{mills1956marginal} and Freund~\cite{freund1985postoptimal} characterized the directional derivatives of the optimal value with respect to the matrix coefficients, and Robinson~\cite{robinson1977characterization} characterized when the optimal solution sets of an LP are stable under arbitrary small perturbations of all the data, including the constraint matrix; see also~\cite{robinson1975stability} for linear inequality systems. More recent work has also examined changes to the constraint matrix~\cite{istace2026analyzing,boveroux2026complexity,derval2025efficient}. Moreover, Derval et al.~\cite{derval2024sensitivity} bound the optimal value of LPs and MILPs whose constraint matrix depends linearly on a scalar parameter, using robust, coefficient-wise, and Lagrangian relaxations. For integer programming, Cook et al.~\cite{cook1986sensitivity} established bounds on the distance between optimal solutions under RHS perturbations, while K{\i}l{\i}n{\c c}-Karzan et al.~\cite{kilincc2009approximating} studied perturbations to objective coefficients that preserve optimality. Dawande and Hooker~\cite{Dawande2000Aug} developed an inference-based sensitivity analysis framework for mixed-integer linear programs (MILPs) under general data perturbations, which, in the binary case, can be interpreted through logical clauses associated with a search tree.

Other approaches exploit duality or information from the solution process to bound the optimal value of perturbed instances. Wolsey~\cite{Wolsey1981} developed the theory of price functions for integer programming duality and its use in sensitivity analysis. Cifuentes et al.~\cite{CifuentesDeyXu2026} showed that the copositive dual of the completely positive reformulation of a mixed-binary quadratic program is a strong dual, and used it to obtain LP-style sensitivity analysis with respect to right-hand-side perturbations. Most closely related to our work is the approach of Schrage and Wolsey~\cite{schrage1985sensitivity}, which reuses the branch-and-bound tree of the nominal problem to construct dual bounds for instances with perturbed RHS vectors. This idea was further developed by G{\"u}zelsoy and Ralphs~\cite{GuzelsoyRalphs2007}, who study dual functions for MILPs obtained from a branch-and-bound tree, and was implemented in the solver SYMPHONY~\cite{RalphsGuzelsoy2005}, which stores the tree from a nominal solve and uses it to warm-start subsequent solves. In a related direction, Gamrath et al.~\cite{GamrathHillerWitzig2015} warm-start a branch-and-bound solver from the final search frontier of a previous run when the objective changes or the feasible region shrinks. Our approach extends the framework of Schrage and Wolsey to perturbations of the constraint matrix and objective coefficients, in addition to the RHS, and recovers their method as a special case when only the RHS is perturbed.

\paragraph{Main Contributions.}
The main contributions of this paper are:
\begin{enumerate}
\item {Generalized LP sensitivity analysis:} We develop a dual function $f_{LP}(\pert)$ based on a  \emph{bounded Lagrangian relaxation}, evaluated at nominal optimal dual multipliers. Unlike classical LP sensitivity analysis that focuses mostly only on right-hand-side changes for dual bounds, this construction provides valid lower bounds under arbitrary simultaneous perturbations to the objective coefficients, right-hand side, and constraint matrix. We also establish that, under suitable mild conditions, the gap $z_{LP}(\pert)-f_{LP}(\pert)$ grows at most linearly with the magnitude of the perturbation (Theorem~\ref{thm:change.all}). 

\item {Generalized MILP dual bounds}: We extend this construction to MILPs by reusing the branch-and-bound tree and node-level LP dual information from the nominal solve. The resulting dual function provides valid bounds under general data perturbations (Proposition~\ref{valid.dual.ip}), extending the classical approach of Schrage and Wolsey~\cite{schrage1985sensitivity}. We develop procedures for selecting dual multipliers at nominally infeasible nodes and strategies for searching over subtrees to strengthen the bounds.

\item GPU-native online evaluation: We show that the proposed generalized MILP dual is inherently well suited for GPU computation. By exploiting precomputed nominal information, node-level dual bounds can be evaluated in parallel, while tree-based bound propagation can be implemented using simple vectorized operations. We also develop efficient implementations that exploit sparse and structured perturbations.

\item Optimality inference and variable fixing: We show how the generalized dual bounds can be used, together with a primal upper bound, to certify optimal binary variable fixings and tighten continuous variable bounds for perturbed instances. We further strengthen these inferences through iterative variable fixing, constraint propagation, and node-level implications. These procedures can also be implemented efficiently on GPUs and used to accelerate subsequent exact reoptimization.

\item Computational performance: We conduct extensive experiments on generalized independent set, combinatorial auction, capacitated facility location, and unit commitment instances. We consider both single-parameter perturbations to address ``what-if'' questions and simultaneous changes to multiple problem coefficients, tailored to each application. For simultaneous perturbations, the magnitude of individual changes ranges from 1\% to 10\%, while the {fraction of perturbed entries} ranges from 1\% to 10\% of the perturbable coefficients, depending on the problem class and type of perturbation. Our results show that the proposed approach produces remarkably tight dual bounds in very short times, typically requiring less than a millisecond per perturbation. The bounds are particularly strong for what-if analysis, with gaps often below a fraction of a percent, but remain strong even under simultaneous perturbations to many coefficients. Compared with Gurobi, the proposed GPU-based approach is approximately 2,000--240,000 times faster in the what-if setting and 120--50,000 times faster under multi-parameter perturbations, when measuring the time needed to obtain dual bounds of the same quality. Moreover, in experiments where the optimal value of the perturbed problem is used as the primal cutoff, the proposed variable fixing procedures can substantially accelerate exact reoptimization with Gurobi, although the improvements vary across problem classes: they are particularly significant for generalized independent set (up to a factor of about 10), while providing little or no benefit for unit commitment.

\end{enumerate}

The rest of the paper is organized as follows. Section~\ref{sec:generalized.lp.dual} develops the generalized dual for LPs, and Section~\ref{sec:generalized.mip.dual} extends it to MILPs using the nominal branch-and-bound tree. Section~\ref{sec:inferring.opt.cond} shows how these bounds can be used to fix variables and tighten variable bounds, and Section~\ref{sec:GPU.acceleration} describes the GPU implementation of our algorithm. Section~\ref{sec:numerical} describes the experimental setup, Sections~\ref{sec:whatif} and~\ref{sec:dense} report results for what-if and distributed perturbations, and Section~\ref{sec:conclusion} concludes. 

\section{Generalized dual for linear programs}\label{sec:generalized.lp.dual} To obtain dual bounds for perturbed MILPs, we first develop a methodology for obtaining dual bounds for linear programs under general perturbations, which we consider in this section. In this section, we always assume $\mathcal I = \emptyset$, even if not explicitly stated. 
\subsection{Bounded Lagrangian dual} 
Our proposed generalized dual function is based on a partial Lagrangian relaxation, which we now introduce. Consider the Lagrangian relaxation of $P(\pert)$ obtained by dualizing the constraints $(A+\pert[A])x \geq b+\pert[b]$, while retaining the bounds on $x$:
    \begin{align*}
        \lagr_{\pert}(\dualvar) := & \min_{\ell \leq x \leq u} \sum_{i=1}^n \left((c+\pert[c])_i - \dualvar^\top(A+\pert[A])^i\right)x_i + \dualvar^\top(b+\pert[b])\\
        = & \sum_{i=1}^n\left[\left( (c + \pert[c])_i - \dualvar^\top (A+\pert[A])^i \right)^+\ell_i + \left((c+\pert[c])_i  - \dualvar^\top (A+\pert[A])^i \right)^-u_i\right] + \dualvar^\top(b+\pert[b]),
    \end{align*}
where $(t)^+:=\max\{t,0\}$ and $(t)^-:=\min\{t,0\}$. Now, let 
$$\dualvar_{\pert}^* \in \arg\max_{\dualvar \geq 0} \lagr_{\pert}(\dualvar)$$
denote a Lagrangian dual maximizer for the $\pert$-perturbed problem. In particular, when $\Delta = 0$, $\lambda^*_0$ denotes the $\lambda$-component of a basic optimal solution of the LP dual of $P(0)$, which is available after solving the nominal LP.

We can now state the proposed dual function:

\begin{equation*}
    \dual_{LP}(\pert) := \lagr_{\pert}(\dualvar_{\nominal}^*).
\end{equation*}
The validity of $\dual_{LP}(\pert)$ follows directly from Lagrangian duality.

\begin{proposition}\label{valid.dual.lp}
    If $P(\nominal)$ is feasible, then $\dual_{LP}(\pert)$ is a valid generalized dual, that is
    $$\dual_{LP}(\pert) \leq \val_{LP}(\pert), ~ \forall \pert \in \R^{m\times n}\times\R^m \times \R^n.$$
\end{proposition}

\begin{proof}
    Since $P(\nominal)$ is feasible and bounded, $\dualvar^*_{\nominal}$ is well-defined and a non-negative vector. Therefore, $\lagr_{\pert}(\dualvar_{\nominal}^*)$ is a valid lower bound.
\end{proof}
The construction of the generalized LP dual bound closely parallels the classical LP dual bound. In fact, the generalized LP dual reduces to the classical one when perturbations are restricted to the right-hand side, $\pert = (0, \pert[b], 0)$. The underlying idea is the same: evaluate a dual function at a dual-feasible, though potentially suboptimal, point obtained from the nominal solve. The key difference lies in the use of the Lagrangian dual while retaining bounds on $x$, as the associated feasible region is simply a box and thus independent of problem data (except variable bounds), unlike the classical LP dual whose feasible region depends on $(A, c)$. Consequently, a dual solution obtained from the nominal problem remains admissible under generalized perturbations.

We conclude with an important remark on the construction of the bounded Lagrangian relaxation.
\begin{remark}
An alternative dual function, $f'_{LP}(\Delta)$, can be derived by also dualizing the bound constraints $\ell \le x \le u$, while still minimizing over the box $\ell \le x \le u$, and evaluating the resulting Lagrangian at the nominal optimal multipliers of all constraints. We remark that this option is dominated by the proposed dual, $\dual'_{LP}(\pert) \leq \dual_{LP}(\pert),$ as shown in Appendix \ref{appendix:alternative.dual}, and further, evaluating $\dual'_{LP}(\pert)$ requires storing an additional $2n$ scalars over $\dual_{LP}(\pert)$.
\end{remark}

\subsection{Worst-case analysis}

The generalized LP dual provides a valid lower bound under arbitrary perturbations, but validity alone does not guarantee that the bound remains useful when the problem data change. Since our approach relies on reusing  nominal dual multipliers, it is important to understand how much the quality of the bound can deteriorate under small perturbations.

We now study how much $\dual_{LP}(\pert)$ may underestimate $\val_{LP}(\pert)$ in the worst case. Under reasonable feasibility assumptions on the nominal and perturbed problems, we next show that the gap between the generalized dual and the value function grows at most linearly with the size of the perturbation. The bound boils down to evaluating the cost incurred by evaluating the bounded Lagrangian dual at $\dualvar^*_{\nominal}$ rather than the optimal $\dualvar^*_{\pert}$.

The constants governing these bounds depend on the conditioning of the constraint system defined by $A$. 
Specifically, this is captured through the Hoffman constant, whose definition we now recall.

\begin{definition}
    For $A \in \R^{m\times n}$, the \textit{Hoffman constant} is defined as
    $$\kappa_A = \max\left\{\|B^{{-1}}\|_\infty:B \text{ is an invertible submatrix of } \begin{bmatrix} A^\top, & I, & -I
\end{bmatrix}\right\},$$
    {where $\|C\|_\infty = \max_{i,j \in[n]}|C_{ij}|$ denotes the largest absolute entry in the matrix.}
\end{definition}
We also write $\|C\|_{1,1} := \sum_{i,j}|C_{ij}|$ for the sum of the absolute values of the entries of $C$. We now state the main results of this section, which bound the gap between the value and dual functions in the LP setting. Let $F(A,b) := \{x \in \R^n: Ax \geq b, \ell \leq x \leq u\}$.

\begin{theorem}\label{thm:change.all}
    Suppose $P(0)$ and $P(\pert)$ are feasible and $\mathcal I = \emptyset$. Let $x^0 \in F(A,b)$. Then for all $\pert$ satisfying:
    \begin{eqnarray}[\pert[b]]_i + M \sum_{j=1}^n|[\pert[A]]_{ij}| \leq (Ax^0-b)_i, \quad \forall i \in [m], \label{eq:int}
    \end{eqnarray}
    we have that:
    \begin{align*}
        \val_{LP}(\pert) - \dual_{LP}(\pert) \leq  \kappa_A\|c\|_1\|\pert[b]\|_1 + 2\boundnorm\|\pert[c]\|_1 + 2\boundnorm \kappa_A\|c\|_1\|\pert[A]\|_{1,1}.
        \end{align*}
    where $\boundnorm := \max\{\|\ell\|_\infty, \|u\|_\infty\}.$ Moreover, if $\pert[A] = \bm{0}$, then the result holds without  assumption (\ref{eq:int}).
\end{theorem}
\begin{remark}
Assumption~\eqref{eq:int} is a quantitative Slater-type condition: it requires a point $x^0\in F(A,b)$ whose slack in each row is large enough to absorb the worst-case change of that row over the box $[\ell,u]$. Some assumption of this type is unavoidable when the constraint matrix is perturbed. Indeed, Robinson~\cite{robinson1977characterization} showed that the optimal solution sets of an LP are stable under arbitrary small perturbations of all the data, including the constraint matrix, if and only if both the primal and the dual optimal solution sets are bounded, and boundedness of the dual optimal set is closely tied to the existence of a Slater point. 
We note that Theorem~\ref{thm:change.all} differs from~\cite{robinson1977characterization} in that it bounds the gap of a specific dual bound, evaluated at the nominal multiplier, rather than the distance between optimal solution sets, and that its constants are explicit.
\end{remark}
The remainder of this subsection is dedicated to proving Theorem~\ref{thm:change.all}. We begin by introducing a useful standard lemma, which bounds the size of an optimal solution to a feasible LP.
\begin{lemma}\label{lemma:bounded.soln}
 Assume $P(\nominal)$ is feasible. Any optimal basic feasible solution $\dualvar^*_{\nominal} \in \arg\max_{\dualvar \geq 0} \lagr_{\nominal}(\dualvar)$ is bounded: 
    $$\|\dualvar^*_{\nominal}\|_\infty \leq \kappa_A \|c\|_1.$$
\end{lemma}

\begin{proof}
    Since the primal is feasible and bounded, the dual is feasible and has an optimal basic feasible solution. Any optimal basic feasible solution of the dual is the solution to a system of linear equations of the form, 
    $$Tw = c,$$
    where $T$ is an invertible submatrix of $[A^\top, I, -I]$. In other words, $\|w\|_{\infty} = \|T^{-1}{c}\|_{\infty} \leq \|T^{-1}\|_{\infty}\|{c}\|_1 \leq \kappa_{A}\|c\|_1.$ 
\end{proof}

We now provide the proof for the main theorem. The proof separates the worst-case effect of the perturbation into three components. Objective perturbations are directly controlled by the boundedness of the variables. Right-hand side perturbations are controlled through the magnitude of the optimal dual multipliers. The main technical challenge is the perturbation to the constraint matrix which changes the geometry of the feasible region. We handle this using the assumption stated in the theorem, converting matrix changes into a conservative right-hand side perturbation, reducing the analysis to right-hand-side sensitivity.

\begin{proof}[Proof of Theorem \ref{thm:change.all}]
    For $\pert = (\pert[A], \pert[b], \pert[c])$ we start by decomposing $\val_{LP}(\pert)-\dual_{LP}(\pert):$
    \begin{align*}
        \val_{LP}(\pert[A], \pert[b], \pert[c]) - \dual_{LP}(\pert[A], \pert[b], \pert[c])  & = \val_{LP}(\pert[A], \pert[b], \pert[c]) - \val_{LP}(\pert[A], \pert[b], 0)\\
        & + \val_{LP}(\pert[A], \pert[b], 0) - \val_{LP}(0, 0, 0) \\
        & + \val_{LP}(0,0,0) - \dual_{LP}(\pert[A], \pert[b], \pert[c]).
    \end{align*}
    Note that $\val_{LP}(\pert[A], \pert[b], 0)$ is finite since the perturbed problem $P(\pert)$ was assumed to be feasible. We bound each term in order:
    \begin{itemize}
        \item Let $x^*$ be an optimal solution to the problem corresponding to $\val_{LP}(\pert[A], \pert[b], 0)$. We have that
        $$\val_{LP}(\pert[A],\pert[b], \pert[c]) - \val_{LP}(\pert[A], \pert[b], 0) \leq (c+\pert[c])^\top x^* - c^\top x^* = \pert[c]^\top x^*.$$. 
       \item Let $r_i = M\sum_{j=1}^n|[\pert[A]]_{ij}|$. For every $x\in[\ell,u]$ and every $i\in[m]$, $|(\pert[A]x)_i| \leq \sum_{j=1}^n|[\pert[A]]_{ij}||x_j| \leq M\sum_{j=1}^n|[\pert[A]]_{ij}| = r_i,$
        that is,
        \begin{equation}
        -r\leq \pert[A]x\leq r\qquad\text{for all } x\in[\ell,u]. \label{eq:r.bound}
        \end{equation} We will first show that
        \begin{equation} \emptyset \neq F(A, b+\pert[b]+r) \subseteq F(A+\pert[A], b+\pert[b]). \label{eq:subset.polytope}\end{equation}
        By the stated assumption, $(\pert[b])_i + M\sum_{j=1}^n |(\pert[A])_{ij}| = (\pert[b])_i + r_i \leq (Ax^0 - b)_i$ for $i \in [m]$, which implies that $x^0 \in F(A, b+\pert[b]+r)$, that is $F(A,b+\pert[b]+r)\neq \emptyset$. Now if $x \in F(A, b+\pert[b]+r),$ then 
    $Ax \geq b + \pert[b]+r,$
    and hence
    $$(A+\pert[A])x = Ax + \pert[A]x \geq b+\pert[b] + r -r = b+\pert[b],$$
    where the inequality follows from (\ref{eq:r.bound}).
    Thus, $x \in F(A+\pert[A], b + \pert[b]),$ proving \eqref{eq:subset.polytope}.
    
    Note now that \eqref{eq:subset.polytope} implies that
    $z_{LP}(0, \pert[b]+r, 0) \geq z_{LP}(\pert[A], \pert[b], 0),$ as the former term is the result of minimizing the same objective over a subset of the feasible region of the latter. Define $\Delta_+ = (0,\pert[b]+r,0)$ as the perturbation of the nominal problem.  The problem $P(\Delta_+)$ is feasible as it contains $x^0$, and is bounded, so by LP duality there is a basic optimal solution $\lambda_+^*\in\arg\max_{\lambda\geq 0}\lagr_{\Delta_+}(\lambda)$ with $\lagr_{\Delta_+}(\lambda_+^*)=\val_{LP}(0,\pert[b]+r,0)$.
 Moreover, for every $\lambda\geq 0$,
        $$\lagr_{\Delta_+}(\lambda) = \min_{\ell\leq x\leq u}\big\{(c-A^\top\lambda)^\top x\big\} + \lambda^\top(b+\pert[b]+r) = \lagr_0(\lambda) + \lambda^\top(\pert[b]+r).$$ Then we obtain,
    \begin{align*}
        z_{LP}(\pert[A], \pert[b], 0) - z_{LP}(0,0,0) &\leq z_{LP}(0,\pert[b]+r,0) - z_{LP}(0,0,0)\\
        &= \mathcal{L}_{\Delta_+}(\lambda_+^*) - \mathcal{L}_0(\lambda^*_0)\\
        &\leq \mathcal{L}_{\Delta_+}(\lambda^*_+) - \mathcal{L}_0(\lambda^*_+)\\
        &=(\lambda_+^*)^\top(\pert[b]+r).
    \end{align*} The equality in the second line is strong duality for $P(\Delta_+)$ and $P(0)$. The inequality in the third line holds because $\lambda_0^*$ maximizes $\lagr_0$ over $\lambda\geq 0$ and $\lambda_+^*\geq 0$, so $\lagr_0(\lambda_0^*)\geq\lagr_0(\lambda_+^*)$. The last equality is the identity displayed above.
    \item We now bound the final term, $z_{LP}(0,0,0)-f_{LP}(\pert).$ 
    Since $\lambda_0^*\geq 0$, \eqref{eq:r.bound} gives $-(\lambda_0^*)^\top\pert[A]x\geq-(\lambda_0^*)^\top r$ for all $x\in[\ell,u]$. We use this to lower bound the generalized dual:
\begin{align*}
f_{LP}(\pert) = \mathcal{L}_{\pert}(\lambda_0^*)
&=
\min_{\ell\leq x\leq u}
\Bigl\{
(c-A^\top\lambda_0^*)^\top x
+(\lambda_0^*)^\top b
+\pert[c]^\top x
-(\lambda_0^*)^\top\pert[A]x
+(\lambda_0^*)^\top\pert[b]
\Bigr\}\\
& \geq  \min_{\ell\leq x\leq u}
\Bigl\{
(c-A^\top\lambda_0^*)^\top x
+(\lambda_0^*)^\top b\Bigr\} 
+\min_{\ell \leq x \leq u}\Bigl\{\pert[c]^\top x
-(\lambda_0^*)^\top\pert[A]x
+(\lambda_0^*)^\top\pert[b]
\Bigr\} \\
& \geq \mathcal{L}_0(\lambda_0^*)
+\min_{\ell\leq x\leq u}\pert[c]^\top x
+(\lambda_0^*)^\top\pert[b]
-(\lambda_0^*)^\top r.
\end{align*}
Using $\val_{LP}(0,0,0)=\lagr_0(\lambda_0^*)$ by strong duality, we bound the final term:
\begin{align*}
\val_{LP}(0,0,0)-\dual_{LP}(\pert)
& = \mathcal{L}_0(\lambda_0^*) - \mathcal{L}_{\Delta}(\lambda_0^*) \\
& \leq 
-\min_{\ell\leq x\leq u}\pert[c]^\top x
-(\lambda_0^*)^\top\pert[b]
+(\lambda_0^*)^\top r.
\end{align*}
\end{itemize}
Collecting all the terms, we obtain:
\begin{align*}
    z_{LP}(\pert)-f_{LP}(\pert) &\leq \pert[c]^\top x^* + (\lambda_+^*)^\top(\pert[b]+r) -(\lambda^*_0)^\top\pert[b] + (\lambda_0^*)^\top r -\min_{\ell\leq x \leq u}\pert[c]^\top x\\
    & = \left[ \pert[c]^\top x^* -\min_{\ell\leq x \leq u}\pert[c]^\top x \right] + \left[ (\lambda_+^*-\lambda_0^*)^\top \pert[b]\right] + \left[(\lambda_+^*+\lambda_0^*)^\top r\right].
\end{align*}
We again bound terms in order. 
\begin{itemize}
    \item For the first term, recall that $x^* \in \arg\min_{x \in F(A+\pert[A], b+\pert[b])}c^\top x.$ Since $F(A+\pert[A], b+\pert[b]) \subseteq [\ell,u]$, we obtain \begin{align*}
{\pert[c]}^\top x^*
-\min_{\ell\le x\le u}{\pert[c]}^\top x
&\le
\max_{\ell\leq x\leq u}\pert[c]^\top x
-\min_{\ell\leq x\leq u}\pert[c]^\top x\\
&=\sum_{j=1}^n(u_j-\ell_j)|\pert[c]|_j \leq 2\boundnorm\|\pert[c]\|_1.
\end{align*}
\item For the final two terms, recall that by Lemma~\ref{lemma:bounded.soln} applied to the basic optimal multipliers of $P(0)$ and of $P(\Delta_+)$, we have, $ 0\le\lambda_0^*,\lambda_+^*\le\kappa_A\|c\|_1\mathbf{1}.$ 
Thus $\|\lambda^*_+-\lambda^*_0\|_\infty \leq \kappa_A\|c\|_1$, and $\|\lambda_+^* + \lambda_0^*\|_\infty \leq 2 \kappa_A \|c\|_1$. We obtain, $$(\lambda_+^*-\lambda_0^*)^\top \pert[b] + (\lambda_+^*+\lambda_0^*)^\top r\ \leq \kappa_A\|c\|_1\|\pert[b]\|_1 + 2\kappa_A\|c\|_1M\|\pert[A]\|_{1,1}.$$
\end{itemize}
The final bound is obtained by summing all the terms.

It remains to treat the case $\pert[A]=\bm{0}$ without assumption \eqref{eq:int}. Assumption \eqref{eq:int} is used only to show that $F(A,b+\pert[b]+r)\neq\emptyset$, which guarantees that $P(\Delta_+)$ has a basic optimal dual solution $\lambda_+^*$. If $\pert[A]=\bm{0}$, then $r=0$ and $\Delta_+=(0,\pert[b],0)$. In this case $F(A,b+\pert[b]+r)=F(A+\pert[A],b+\pert[b])$, which is the feasible region of $P(\pert)$ and is nonempty by hypothesis. Hence $P(\Delta_+)$ is feasible, and it is bounded because its feasible region lies in $[\ell,u]$. Thus $\lambda_+^*$ exists, and every remaining step of the proof holds verbatim with $r=0$. In particular, the bound reads $\val_{LP}(\pert)-\dual_{LP}(\pert)\leq\kappa_A\|c\|_1\|\pert[b]\|_1+2\boundnorm\|\pert[c]\|_1$.
\end{proof}

We next show that Theorem~\ref{thm:change.all} is tight: each of the three coefficients in the bound is attained in the limit, and assumption~\eqref{eq:int} cannot be dropped when the constraint matrix is perturbed. For a perturbation $\pert$, let
$$B(\pert) := \kappa_A\|c\|_1\|\pert[b]\|_1 + 2\boundnorm\|\pert[c]\|_1 + 2\boundnorm \kappa_A\|c\|_1\|\pert[A]\|_{1,1}$$
denote the right-hand side of the bound in Theorem~\ref{thm:change.all}.
\begin{proposition}[Tightness of Theorem~\ref{thm:change.all}]\label{prop:tightness}
The following hold.
\begin{enumerate}
    \item For each $\tau\in\{A,b,c\}$,
    $$\sup\left\{\frac{\val_{LP}(\pert)-\dual_{LP}(\pert)}{B(\pert)}\right\}=1,$$
    where the supremum is taken over all instances and all nonzero perturbations $\pert$ that are supported on the block $\pert[\tau]$ and satisfy the hypotheses of Theorem~\ref{thm:change.all}, including assumption~\eqref{eq:int}. The supremum is approached by one-dimensional instances whose nominal dual problem has a unique optimal solution, with $\kappa_A$ arbitrarily large. In particular, none of the three coefficients in Theorem~\ref{thm:change.all} can be replaced by a smaller value.
    \item Assumption~\eqref{eq:int} cannot be dropped when $\pert[A]\neq\bm{0}$. There exists an instance and perturbations $\pert=(\pert[A],0,0)$ with $\|\pert[A]\|_{1,1}\downarrow 0$, such that $P(\pert)$ is feasible and $\val_{LP}(\pert)-\dual_{LP}(\pert)$ is bounded away from zero. In particular, the conclusion of Theorem~\ref{thm:change.all} fails for these perturbations.
\end{enumerate}
\end{proposition}
\begin{proof}
    \begin{enumerate}
    \item By Theorem~\ref{thm:change.all}, each ratio is at most $1$, so it suffices to exhibit instances and perturbations for which the ratio converges to $1$. Fix any $K\geq 1$, let $s=1/K\in(0,1]$, and for $0<\varepsilon<1/2$ set $\delta=\varepsilon^2$. Consider the LP
$$\min\{x : sx\geq s(\delta-1),\ -1\leq x\leq 1\},$$
that is, $n=m=1$, $A=[s]$, $b=s(\delta-1)$, $c=1$, $\ell=-1$ and $u=1$. Then $\boundnorm=1$ and $\|c\|_1=1$. The invertible submatrices of $[A^\top, I, -I]^\top=(s,1,-1)^\top$ are $[s]$, $[1]$ and $[-1]$, and since $s\leq 1$ we obtain $\kappa_A=1/s=K$. The feasible region of $P(\nominal)$ is $[\delta-1,1]$, so $\val_{LP}(\nominal)=\delta-1$. For $\dualvar\geq 0$,
$$\lagr_{\nominal}(\dualvar)=\min_{-1\leq x\leq 1}(1-s\dualvar)x+s(\delta-1)\dualvar=-|1-s\dualvar|+s(\delta-1)\dualvar.$$ This function is concave and piecewise linear, with slope $s\delta>0$ for $\dualvar<1/s$ and slope $s(\delta-2)<0$ for $\dualvar>1/s$. Hence $\dualvar_{\nominal}^*=1/s$ is its unique maximizer. It is also basic, since it solves $Tw=c$ with $T=[s]$. Consequently, for every $\pert=(\pert[A],\pert[b],\pert[c])\in\R^3$,
$$\dual_{LP}(\pert)=\lagr_{\pert}(1/s)=\min_{-1\leq x\leq 1}\left(\pert[c]-\frac{\pert[A]}{s}\right)x+\delta-1+\frac{\pert[b]}{s}=\delta-1+\frac{\pert[b]}{s}-\left|\pert[c]-\frac{\pert[A]}{s}\right|.$$
Now consider the following three perturbations. For each of them, $\dual_{LP}(\pert)=\delta-1-\varepsilon=\varepsilon^2-1-\varepsilon$.
\begin{itemize}
    \item $\pert=(0,0,-\varepsilon)$. The feasible region is unchanged and $1-\varepsilon>0$, so $\val_{LP}(\pert)=(1-\varepsilon)(\delta-1)$.
    \item $\pert=(0,-s\varepsilon,0)$. The perturbed constraint is $x\geq\delta-1-\varepsilon$. Since $\varepsilon\geq\varepsilon^2=\delta$, this lower bound is at most $-1$, so the feasible region is $[-1,1]$ and $\val_{LP}(\pert)=-1$.
    \item $\pert=(s\varepsilon,0,0)$. The perturbed constraint is $s(1+\varepsilon)x\geq s(\delta-1)$, that is, $x\geq(\varepsilon^2-1)/(1+\varepsilon)=\varepsilon-1\in(-1,1)$, so $\val_{LP}(\pert)=\varepsilon-1$.
    \end{itemize}
Now it is easy to verify that:
    \begin{eqnarray*}
    \begin{array}{c|c|c|c}
\pert & \val_{LP}(\pert)-\dual_{LP}(\pert) & B(\pert) & \dfrac{\val_{LP}(\pert)-\dual_{LP}(\pert)}{B(\pert)} \\ \hline
(0,0,-\varepsilon) & \varepsilon(2-\varepsilon^2) & 2\varepsilon & 1-\varepsilon^2/2 \\
(0,-s\varepsilon,0) & \varepsilon-\varepsilon^2 & \varepsilon & 1-\varepsilon \\
(s\varepsilon,0,0) & 2\varepsilon-\varepsilon^2 & 2\varepsilon & 1-\varepsilon/2
\end{array}
    \end{eqnarray*}

where we used $B(0,0,-\varepsilon)=2\boundnorm\varepsilon$, $B(0,-s\varepsilon,0)=\kappa_A\|c\|_1 s\varepsilon$ and $B(s\varepsilon,0,0)=2\boundnorm\kappa_A\|c\|_1 s\varepsilon$.

All three perturbations satisfy the hypotheses of Theorem~\ref{thm:change.all}. Indeed, $P(\nominal)$ and $P(\pert)$ are feasible. Moreover, with $x^0=1\in F(A,b)$ we have $Ax^0-b=s(2-\delta)$, and assumption~\eqref{eq:int} reads $\pert[b]+\boundnorm|\pert[A]|\leq s(2-\delta)$. The left-hand side equals $0$, $-s\varepsilon$ and $s\varepsilon$ for the three perturbations, respectively, and $s\varepsilon<s<s(2-\delta)$ since $\varepsilon<1/2$.

As $\varepsilon\downarrow 0$, all three ratios converge to $1$. This holds for arbitrarily large $\kappa_A$. Finally, suppose that one of the three coefficients in Theorem~\ref{thm:change.all} were replaced by $\theta$ times its value for some $\theta<1$. Choosing $\varepsilon$ small enough that the corresponding ratio in the table exceeds $\theta$ yields a perturbation, supported on the corresponding block, that violates the modified bound.
    \item Consider the nominal LP
$$\min\{-x : 0\cdot x\geq 0,\ -1\leq x\leq 1\},$$
that is, $A=[0]$, $b=0$, $c=-1$. The invertible submatrices of $(0,1,-1)^\top$ are $[1]$ and $[-1]$, so $\kappa_A=1$; moreover $\|c\|_1=1$, $\boundnorm=1$ and $\val_{LP}(\nominal)=-1$. Since $\lagr_{\nominal}(\dualvar)=\min_{-1\leq x\leq 1}(-x)=-1$ for every $\dualvar\geq0$, every $\dualvar\geq 0$ is optimal for the nominal dual.
For any $\eta>0$, consider the matrix perturbation $\pert=(-\eta,0,0)$. The perturbed constraint becomes $-\eta x\geq 0$, that is, $x\leq 0$. Hence $P(\pert)$ is feasible and
$$\val_{LP}(\pert)=0,\qquad \dual_{LP}(\pert)=\min_{-1\leq x\leq 1}(-x)=-1.$$
Therefore,
$$\val_{LP}(\pert)-\dual_{LP}(\pert)=1$$
for every $\eta>0$, even though $\|\pert[A]\|_{1,1}=\eta\downarrow 0$. In particular, since $B(\pert)=2\boundnorm\kappa_A\|c\|_1\|\pert[A]\|_{1,1}=2\eta$, the conclusion of Theorem~\ref{thm:change.all} fails for every $\eta<1/2$, and more generally no bound that is linear in the perturbation size, with a constant depending only on the nominal instance, can hold. This does not contradict Theorem~\ref{thm:change.all}: since $Ax^0-b=0$ for every $x^0\in[-1,1]$, assumption~\eqref{eq:int} reads $\eta = 0$, which fails for every $\eta>0$.
    \end{enumerate}
\end{proof}

\section{Generalized MILP dual}\label{sec:generalized.mip.dual}

In this section, we construct generalized dual bounds for the mixed-integer setting using the LP results from Section~\ref{sec:generalized.lp.dual} as building blocks. The high-level idea is to reuse two sources of information from a nominal solve via branch-and-bound: the nominal tree, which provides a valid disjunction, and nominal LP node multipliers that provide cheap lower bounds for each disjunct using ideas from Section~\ref{sec:generalized.lp.dual}.

Suppose the nominal problem  $P(\nominal)$ has been solved by branch-and-bound, yielding a tree $\TT$ with node set $V(\TT)$.
We represent variable fixings (or lack thereof) at any node $v$ with $\ell_{j}^v \leq x_j \leq u_{j}^{v}, ~ j \in \mathcal{I}$.
Each node $v$ is associated with an LP, and thus admits its own LP value function $z_{LP}^v(\pert)$:
$$\val_{LP}^v(\pert) = \min_{x\in\mathbb{R}^n}\Big\{ (c+\pert[c])^\top x~|~(A+\pert[A])x \geq b+\pert[b], ~ \ell^v \leq x_i \leq u^v\Big\}.$$ 
Suppose that for each node $v$, we are given a valid node-level generalized dual bound $\dual_{LP}^v(\pert) \leq \val_{LP}^v(\pert)$. Let $\LL(\TT) \subseteq V(\TT)$ denote the set of leaf nodes of tree $\TT$. If $\TT$ is a pure branch-and-bound tree (that is, a tree with no cuts or preprocessing), we define the tree-based generalized dual as:
\begin{align}\label{eqn:min.leaf}\dual(\pert) = \min_{v \in \LL(\TT)}\dual_{LP}^v(\pert).\end{align}
\begin{proposition}\label{valid.dual.ip}
    If $P(\nominal)$ is feasible, then $\dual(\pert)$ is a valid generalized dual, that is
    $$\dual(\pert) \leq \val(\pert), ~ \forall \pert \in \R^{m\times n} \times \R^m \times \R^n.$$ 
\end{proposition}
\begin{proof}\label{prop:valid.dual.mip}
    The leaves of a pure branch-and-bound tree define a valid partition of the integer-feasible domain, independent of problem data $A$, $b$, and $c$. Thus, for any pure branch-and-bound tree $\TT$, the following holds:
    $$\val(\pert) \geq \min_{v \in \LL(\TT)}\val_{LP}^v(\pert) \geq \min_{v \in \LL(\TT)}\dual_{LP}^v(\pert) = f(\pert).$$
\end{proof}
The use of branch-and-bound trees for sensitivity analysis in mixed-integer programming has been explored previously. In particular, \cite{schrage1985sensitivity} develop a related approach for deriving tree-based dual bounds. However, their reliance on LP duality to construct node-level bounds restricts their analysis to right-hand-side perturbations, whereas our framework accommodates perturbations to all problem data.

Proposition~\ref{valid.dual.ip} leads to a natural
offline-online procedure. Offline, we solve the nominal MIP using pure branch-and-bound, storing dual information at every tree node. Online, given perturbation $\pert$, we obtain $\dual(\pert)$ by evaluating the generalized LP duals $\dual_{LP}^v(\pert)$ associated with each leaf node and return their minimum.

\subsection{Offline node-level dual choice}\label{sec:dualselection}
It remains to specify what dual information we store at each node $v$ such that we can obtain valid dual bounds $\dual_{LP}^v(\pert), ~\forall v \in V(\TT)$ online. As before, let $\dual_{LP}^v(\pert) = \lagr_{\pert}^v(\dualvar^v).$ Nodes $v \in V(\TT)$ are either nominally primal feasible or nominally primal infeasible, as we assume boundedness. Thus for nominally primal feasible nodes, we have the straightforward choice of setting $\dualvar^v = (\dualvar^{*}_{\nominal})^v,$  where $(\dualvar^{*}_{\nominal})^v$ is a basic optimal dual solution of the nominal node LP. We have not explored the selection of a specific optimal dual solution. 

However, when a node is primal infeasible, it is in turn dual unbounded, so the dual optimal solution is not well-defined. We therefore introduce a heuristic that seeks to balance the strength of the nominal dual with its sensitivity to perturbations. Specifically, we seek multipliers with small norm that nevertheless certify a sufficiently strong nominal lower bound, since sensitivity to data perturbations is determined in part by multiplier size. For some target threshold $\gamma \geq 0,$ we set $\dualvar^v = (\dualvar^{\infty}_{\nominal})^v$ obtained by solving (where the superscript $\infty$ indicates that the nominal node LP is infeasible, so that its dual is unbounded): 
    \begin{align*}
(\dualvar^{\infty}_{\nominal})^v \in \arg\min_{\lambda \geq 0}\Big\{ \|\dualvar\|_1 : \mathcal{L}^v_{\nominal}(\lambda) \geq z(\nominal) + \gamma|z(\nominal)|\label{constr:min.lagrangian.obj}\Big\}.
    \end{align*}
Recall that $z(\nominal)$ denotes the optimal value of the nominal MILP $P(\nominal)$. This value is available when selecting multipliers for nominally infeasible nodes, since these multipliers are computed only after $P(\nominal)$ has been solved to optimality by branch-and-bound in the offline phase. In our experiments, we set $\gamma=0.5$, which we found to work well empirically.

\subsection{Online bound enhancement via full-tree information}\label{subsec:bound.enhancement.techniques}

The tree-based generalized MIP dual \eqref{eqn:min.leaf} is valid for any pure branch-and-bound tree  (by Proposition~\ref{prop:valid.dual.mip}). We use the nominal branch-and-bound solve to determine both $\TT$ and the node-level dual multipliers $\dualvar^v$, with the expectation that this nominal information will provide strong bounds under small perturbations. While the stored multipliers can yield tight bounds for individual node LP relaxations (see Theorem~\ref{thm:change.all}), the resulting tree-based dual can be brittle: even a single weak leaf bound $\dual_{LP}^v(\pert)$ can significantly degrade the overall bound. This suggests that using the entire nominal tree may not always be the best choice.

To address this issue, we consider selecting a subtree $\TT' \subseteq \TT$ obtained by pruning the nominal tree. We restrict attention to such subtrees to limit online computational overhead, since node-level dual information has been stored offline only for nodes of $\TT$. Importantly, if the perturbed node LP relaxations were solved exactly, the full tree would yield the strongest bound among all its subtrees, since LP relaxation bounds are nondecreasing under branching. However, bounds obtained from stored dual multipliers need not satisfy this monotonicity, and pruning certain branches may therefore improve the resulting dual bound. A natural question, also asked by \cite{schrage1985sensitivity}, is: what is the strongest dual bound obtainable over all such subtrees?

\begin{equation}\label{eqn:opt.tree}
\max_{\TT' \subseteq \TT} \min_{v\in\LL(\TT')} \dual_{LP}^v(\pert).
\end{equation}

A natural dynamic program may be used to compute \eqref{eqn:opt.tree}. For each non-terminal node $v$, let $L(v)$ and $R(v)$ denote its left and right children respectively. Define
\begin{equation}\label{eq:full.propagation}
    BB^v(\pert) := \begin{cases}
        \dual_{LP}^v(\pert), & \text{if node $v$ is a terminal node,}\\
        \max\{\dual_{LP}^v(\pert), \min(BB^{L(v)}(\pert), BB^{R(v)}(\pert))\}, & \text{else.}
    \end{cases}
\end{equation}
$BB^v(\pert)$ is the best bound that can be achieved by a subtree rooted at $v$ under perturbation $\pert$: either we prune the subtree at $v$ itself, or we retain both children and take the minimum of the best bounds achieved at the subtrees rooted at the left and right children. The value of problem \eqref{eqn:opt.tree} is thus the best bound obtained at the subtree rooted at the root, $BB^{\mathrm{root}}(\pert)$, {see Theorem 1 in \cite{schrage1985sensitivity} for a proof}. In the future sections we will refer to this bound as \textbf{full-propagation}. Computing $BB^{\mathrm{root}}(\pert)$ online requires traversing every node in $\TT$, which can be computationally expensive. However, as we show in Section~\ref{sec:GPU.acceleration}, the computations at nodes within the same tree level can be performed in parallel, making the procedure well-suited for GPU acceleration.

\section{Inferring optimality conditions}\label{sec:inferring.opt.cond}
The generalized MIP dual provides an inexpensive lower bound on the \textit{optimal value} of a perturbed problem. We now show how the same nominal information can also be used to infer properties of \textit{optimal solutions}. In particular, we show how generalized MIP dual bounds can be used to certify \emph{variable fixings and tighten variable bounds without solving additional LPs}. Besides providing insight into how optimal solutions behave under perturbations, these reductions can be passed to subsequent solves to reduce their computational cost~\cite{savelsbergh1994preprocessing}.

In Sections~\ref{subsubsec:int.var.fix} and \ref{subsubsec:cont.var.tightening}, we show how generalized MIP dual bounds can be used to certify binary variable fixings and tighten continuous variable bounds through optimality-based arguments. Then, in Section~\ref{subsubsec:constr.bound.prop}, we discuss how constraint propagation can further tighten variable bounds.

Although the results of Sections~\ref{sec:generalized.lp.dual} and~\ref{sec:generalized.mip.dual} hold for general finite bounds, in this section we assume for simplicity that all integer variables are binary, and we denote their index set by $\mathcal{B} := \mathcal{I}$. This assumption is satisfied by all instances in our experiments, and the approach extends to general integer variables. Continuous variables $i \in \mathcal{C} := [n] \setminus \mathcal{I}$ retain general finite bounds, which are not modified by branching.

\subsection{Integer variable fixing}\label{subsubsec:int.var.fix}

Optimality-based variable fixing is a classic technique for eliminating variable assignments by comparing lower bounds on restricted problems against a known primal solution. We show how to efficiently obtain such bounds using the generalized MILP dual.

For variable $i\in\mathcal{B},$ node $v \in V(\TT)$, we define:
\begin{align}
    \dual^v_{LP}(\pert; x_i=0) := \min_{\substack{\ell^v\leq x \leq u^v\\x_i=0}} \big\{(c+\pert[c]-(A+\pert[A])^\top \lambda^v)^\top x+(b+\pert[b])^\top\lambda^v\big\}
\end{align}
By Lagrangian duality, $\dual^v_{LP}(\pert;x_i=0)$ is a valid lower bound on node $v$'s LP relaxation with the additional restriction that $x_i=0$. Let $\bar{c}^v(\Delta) = (c + \Delta_c) - (A + \Delta_A)^\top \lambda^v$. We can directly compute $\dual^v_{LP}(\pert; x_i=0)$ using the node's bounds on $x_i$:
\begin{align}
    f^v_{LP}(\pert;x_i=0) &= \begin{cases}
        f^v_{LP}(\pert) & \text{if } \ell_i^v = u_i^v = 0,\\
        f^v_{LP}(\pert) - [\bar{c}^v_i(\pert)]^-  & \text{if } \ell_i^v = 0, u_i^v=1,\\
        + \infty & \text{if } \ell_i^v = u_i^v = 1.
    \end{cases}
\end{align}
The expression follows directly from the variable bounds $\ell_i^v, u_i^v$. If $x_i=0$ is already enforced in node $v$, then the dual bound remains unchanged.
If $x_i$ is free, {its contribution to the Lagrangian minimum changes from $\min\{0, \bar{c}_i^v(\pert)\} = [\bar{c}_i^v(\pert)]^-$ to zero; the correction is therefore $-[\bar{c}_i^v(\pert)]^-.$}
Finally, if $x_i=0$ contradicts the node bounds, the resulting minimization is over the empty set, and thus its value is $+\infty$.

Following the same logic, the following is a valid lower bound on node $v$'s LP relaxation with the constraint $x_i=1$:
\begin{align}
    f^v_{LP}(\pert;x_i=1) &= \begin{cases}
        +\infty & \text{if } \ell_i^v = u_i^v = 0,\\
        f^v_{LP}(\pert) + [\bar{c}^v_i(\pert)]^+  & \text{if } \ell_i^v = 0, u_i^v=1,\\
        f^v_{LP}(\pert) & \text{if } \ell_i^v = u_i^v = 1.
    \end{cases}
\end{align}

To obtain a lower bound on $P(\pert)$ with the additional constraint, $x_i = s$, $s \in \{0,1\}$, we apply the same tree-based construction as in Section \ref{sec:generalized.mip.dual}, specifically the full-propagation procedure. Define:
\begin{align*}
    BB^v(\pert;x_i=s) := \begin{cases}
        \dual_{LP}^v(\pert;x_i=s), & \text{if } v \in \LL(\TT)\\
        \max\left\{\dual_{LP}^v(\pert; x_i=s), \min\left(BB^{L(v)}(\pert;x_i=s), BB^{R(v)}(\pert;x_i=s)\right)\right\}, & \text{else.}
    \end{cases}
\end{align*}
We let:
$$\dual(\pert; x_i=s) := BB^{\text{root}}(\pert;x_i=s).$$
By the same pure branch-and-bound tree reasoning as for the generalized MILP dual, $\dual(\pert; x_i=s)$ is a valid lower bound on the optimal value of $P(\pert)$ restricted to $x_i=s$. 

We next describe a standard bound-based variable-fixing criterion~\cite{savelsbergh1994preprocessing} using the dual bounds obtained from our framework. Given a primal upper bound on $P(\pert)$, these bounds can identify binary variables whose values can be fixed in all optimal solutions, potentially reducing the computational cost of solving the perturbed problem.

\begin{proposition}\label{prop:binary.variable.fixing}
    Let $\bar{z}_{\pert}$ be a primal upper bound on $P(\pert)$. Then, for any $i\in\mathcal{B}$, any optimal solution $x^*_{\pert}$ to $P(\pert)$ satisfies:
    \begin{align*}
        (x^*_{\pert})_i = \begin{cases}
            0, & \text{if } \dual(\pert; x_i=1) > \bar{z}_{\pert},\\
            1, & \text{if } \dual(\pert; x_i=0) > \bar{z}_{\pert}.
        \end{cases}
    \end{align*}
\end{proposition}
\begin{proof}
    Suppose $\dual(\pert;x_i=1)>\bar{z}_{\pert}$. Since $\dual(\pert;x_i=1)$ is a valid lower bound on the problem restricted to $x_i=1$, every feasible solution with $x_i = 1$ has objective value at least $\dual(\pert; x_i = 1) > \bar{z}_{\pert} \ge z(\pert)$, so it is not optimal. Thus, every optimal solution must satisfy $x_i^*=0$. The second case follows symmetrically.
\end{proof}

\paragraph{Iterative variable fixing.}
    {Incorporating already certified fixings into $f(\pert)$ and $f(\pert;x_i=s)$ can strengthen the bounds and certify further fixings. The procedure proceeds in rounds: in each round, every binary variable not yet fixed is tested, with all fixings certified in earlier rounds imposed when evaluating the bounds.}
    Suppose we have certified that every optimal solution satisfies $x_i^* = \bar{x}_i$ for all $i \in \bar{\mathcal{B}} \subseteq \mathcal{B}$. For each node $v \in V(\TT)$, let $\mathcal{B}_v^{\text{free}} := \{i \in \bar{\mathcal{B}} : \ell_i^v =0, u^v_i=1\}$ denote the set of indices of previously fixed variables that remain free at node $v$. The restricted node-level dual is given by:
    \begin{align*}
        f^v_{LP}(\pert; x_{\bar{\mathcal{B}}} = \bar{x}) = \begin{cases}
        + \infty & \text{if } \exists i \in \bar{\mathcal{B}} \text{ s.t. } \bar{x}_i \not \in [\ell_i^v, u_i^v]\\
        f_{LP}^v(\pert)+\sum_{i\in \mathcal{B}_v^{\text{free}}} [(2\bar{x}_i-1)\bar{c}_i^v(\pert)]^+ & \text{else.}
    \end{cases}
    \end{align*}
    The correction follows from the fact that fixing a binary variable changes its contribution to the Lagrangian minimum from $\min\{0,\bar{c}_i^v(\pert)\}$ to $\bar{x}_i\bar{c}_i^v(\pert)$.
    
    Incorporating certified fixings can strengthen the generalized MIP dual in two ways: eliminating nodes that are incompatible with the fixings and strengthening the dual bound of the remaining nodes. These improvements may certify additional fixings, giving rise to the following iterative procedure:
    \begin{enumerate}
        \item Initialize $\bar{\mathcal{B}}=\emptyset$.
        \item For each $i\in \mathcal{B}\setminus \bar{\mathcal{B}}$, evaluate $f(\pert;x_i=s, x_{\bar{\mathcal{B}}} = \bar{x})$, $s\in\{0,1\}$, and identify any fixings using Proposition~\ref{prop:binary.variable.fixing}, which applies since every optimal solution satisfies $x_{\bar{\mathcal{B}}} = \bar{x}$.
        \item Add all newly certified fixings to $\bar{\mathcal{B}}.$ If additional fixings were found, return to step 2; else, terminate.
    \end{enumerate}
    {The first round, with $\bar{\mathcal{B}} = \emptyset$, is the \textit{single-pass} procedure. The \textit{binary iteration} procedure continues until a complete round certifies no new binary fixings.}

\subsection{Continuous variable bound tightening} \label{subsubsec:cont.var.tightening}
We now extend optimality-based reasoning to tighten the bounds of continuous variables.
Let $\mathcal{C} = [n] \setminus \mathcal{I}$ denote the indices of continuous variables. For any node $v$ and $t \in [\ell^v_i, u^v_i],$ the restricted node-level duals yield:
\begin{align*}
    &f^v_{LP}(\pert;x_i\geq t) = f^v_{LP}(\pert) + [\bar{c}^v_i(\pert)]^+(t - \ell_i^v),\\
    &f^v_{LP}(\pert;x_i\leq t) = f^v_{LP}(\pert) + [-\bar{c}^v_i(\pert)]^+(u^v_i - t).
\end{align*}
For example, when $\bar c_i^v(\pert)>0$, the Lagrangian minimization sets $x_i$ to its lower bound $\ell_i^v$. Imposing $x_i\geq t$ therefore increases the dual bound by $\bar c_i^v(\pert)(t-\ell_i^v)$.

We can use these identities to tighten variable bounds using a primal upper bound $\bar z_{\pert}$. Consider any feasible solution $x$ of the perturbed MILP $P(\pert)$ that satisfies the branching constraints at node $v$ and has objective value at most $\bar{z}_{\pert}$. Since this solution also satisfies the restriction $x_i\geq x_i$, validity of the restricted dual bound gives 
\begin{align*}
f_{LP}^v(\pert) +[\bar c_i^v(\pert)]^+(x_i-\ell_i^v) =f_{LP}^v(\pert;x_i\geq x_i) \leq \bar z_{\pert}.
\end{align*}
Consequently, if $\bar c_i^v(\pert)>0$, rearranging yields
$$x_i\leq \ell_i^v+ \frac{\bar z_{\pert}-f_{LP}^v(\pert)} {\bar c_i^v(\pert)}. $$
Thus, the upper bound on $x_i$ can be tightened to
$$\bar{u}^v_i := \begin{cases}
    \min\left\{ u^v_i,\ \ell^v_i + \dfrac{\bar{z}_{\pert} - \dual^v_{LP}(\pert)}{\bar{c}^v_i(\pert)} \right\}, & \text{if } \bar{c}^v_i(\pert) > 0,\\
    u^v_i, & \text{otherwise.}
\end{cases}$$
To obtain a global upper bound on $x_i$, we propagate node-level upper bounds throughout the nominal tree. Define:
\begin{equation}\label{eq:cont.var.ub.tree.prop}
    \mathrm{UB}^v_i = \begin{cases}
        \bar{u}^v_i, & \text{if } v \in \mathbb{L}(\TT),\\
        \min\left\{ \bar{u}^v_i, \max\left\{ \mathrm{UB}^{L(v)}_i, \mathrm{UB}^{R(v)}_i \right\} \right\}, & \text{else.}
    \end{cases}
\end{equation}
An analogous procedure derives a lower bound $\mathrm{LB}^v_i$. The resulting bounds on every optimal solution $x^*_i$ are:
$$\mathrm{LB}^{\mathrm{root}}_i \le x^*_i \le \mathrm{UB}^{\mathrm{root}}_i.$$

\subsection{Constraint-based bound propagation}\label{subsubsec:constr.bound.prop}
We can further tighten variable bounds by applying constraint
propagation \cite{savelsbergh1994preprocessing} to the bounds obtained from the generalized MIP dual.
Consider a node $v$ with current variable bounds
$\ell^v \leq x \leq u^v$, which may already incorporate
previously inferred restrictions. For any constraint $r\in[m]$ of $P(\pert)$ and variable $i \in [n]$, we have:
\begin{align*}
    (A + \pert[A])_{ri} x_i &\geq (b+\pert[b])_r - \sum_{j\neq i}(A+\pert[A])_{rj}x_j\\
    & \geq (b+\pert[b])_r - \sum_{j\neq i}\max\{(A+\pert[A])_{rj}\ell_j^v, (A+\pert[A])_{rj}u_j^v\}.
\end{align*}
Thus, if $(A+\pert[A])_{ri} > 0$, constraint $r$ implies the following lower bound on $x_i$:
$$\hat{\ell}_{ir}^v = \frac{(b+\pert[b])_r-
\sum_{j\neq i}\max\{(A+\pert[A])_{rj}\ell_j^v,
       (A+\pert[A])_{rj}u_j^v\}}{(A+\pert[A])_{ri}}.$$
Considering all constraints, we obtain the tightened lower bound:
$$\hat{\ell}^v_i = \max_{r:(A+\pert[A])_{ri}>0}\hat\ell_{ir}^v.$$
An analogous argument can be made to derive a constraint-based upper bound $\hat{u}^v_i$ for $x_i$ in node $v$. This bound can be refined by propagating node-bounds throughout the tree like \eqref{eq:cont.var.ub.tree.prop}, as well as doing iterative rounds of variable fixing, continuous variable tightening, and further constraint propagation.

\paragraph{Propagation with node implications} The baseline constraint propagation procedure outlined above incorporates variable bounds implied by tree propagation, which hold \textit{globally} over the entire tree. However, an individual node's nominal dual solution may imply stronger \textit{local} bounds, allowing for stronger domain reductions within its subtree.

For example, suppose that for some binary variable $x_i$ and node $v$, it holds that $$f_{LP}^v(\pert;x_i=0)>\bar{z}_{\pert}.$$
This certifies that $x_i=1$ for any optimal solution of $P(\pert)$ lying within the subtree rooted at $v$. This local fixing may hold even when the propagated implication does not, that is, even if $$BB^{\text{root}}(\pert;x_i=0)\leq \bar{z}_{\pert}.$$ 
Enforcing $\ell_i^v=1$ while performing constraint propagation in node $v$, and inheriting this fixing throughout its descendants, can yield further domain reductions.
\section{GPU acceleration}\label{sec:GPU.acceleration}
In this section, we discuss a key aspect of the proposed generalized MILP dual: its tractable evaluation using GPUs.
We begin with a brief overview of GPU computing, aimed at readers unfamiliar with GPUs and the characteristics of algorithms that make them amenable to GPU acceleration. 

Broadly speaking, while CPUs are designed to efficiently execute a relatively small number of complex tasks, GPUs are designed to perform many simpler operations simultaneously.
GPUs excel when the same computation, say $g(\cdot)$, can be applied independently to many inputs $g(x_i)$, $i\in[n]$, as is often the case with matrix and elementwise operations.

Another important consideration is memory.
GPUs have their own memory, separate from the CPU's main memory, which is typically more limited in capacity. 
Transferring data between CPU and GPU memory can introduce
significant overhead.
Ideally, such transfers are limited, with frequently reused information kept in GPU memory whenever possible.
Further, GPU performance depends on how data is stored and accessed.
For example, accessing consecutive blocks of memory is more efficient than accessing scattered locations.

Parallelism and memory usage introduce a natural trade-off: organizing computations into large shared matrix operations allows one to leverage substantial parallelism and data reuse, but may require prohibitively large intermediate arrays.
One way to combat this is to divide computations into smaller batches.
While this can keep memory usage under control, each batch introduces additional overhead (such as launching GPU kernels).
An efficient GPU implementation therefore seeks to balance parallelism, memory usage, and data movement.

We highlight three questions that guided our implementation, which we hope may also be useful in identifying opportunities for GPU acceleration in other optimization algorithms:
\begin{itemize}
    \item \textbf{Parallelism:} Which computations are independent and can be performed simultaneously? Can computations be grouped into larger operations, rather than executed individually?
    \item \textbf{Data reuse:} What information can be shared, precomputed, and/or cached across current and future computations?
    \item \textbf{Memory management:} How can data be stored, accessed, and transferred efficiently? How should computations be batched to balance parallelism and memory usage?
\end{itemize}
In the following subsections, we show how these principles are used to implement ideas from Sections~\ref{sec:generalized.lp.dual}-\ref{sec:inferring.opt.cond} on GPUs.

\subsection{Parallel node evaluation}\label{subsec:gpu.node.eval}

Recall that for each node $v \in V(\TT)$, we let $\bar{c}^v(\pert) := c + \pert[c] - (A+\pert[A])^\top\dualvar^v$, and $\bar{b}^v(\pert) := (\dualvar^v)^{\top}(b + \pert[b])$. Using this notation, we can write the node-level LP dual as:
\begin{equation}
    \dual^v_{LP}(\pert) = (u^v)^\top[\bar{c}^v(\pert)]^- + (\ell^v)^\top[\bar{c}^v(\pert)]^+ + \bar{b}^v(\pert) \label{eq:LP.dual.dense.gpu}
\end{equation}
These LP dual evaluations are independent across all nodes and share the same structure, allowing us to evaluate all $N$ nodes in parallel using matrix operations. To see this explicitly, collect the following matrices:
$L = \begin{bmatrix}\ell^1 & \ell^2 & \cdots& \ell^N\end{bmatrix}, ~ U = \begin{bmatrix}u^1&  u^2 & \cdots & u^N\end{bmatrix}, ~ \bar{C}(\pert) = \begin{bmatrix} \bar{c}^1(\pert) & \bar{c}^2(\pert) & \cdots & \bar{c}^N(\pert)\end{bmatrix}$ all in $\mathbb{R}^{n\times N}$, and $\bar{b}(\pert) = \begin{bmatrix} \bar{b}^1(\pert) & \bar{b}^2(\pert) & \cdots & \bar{b}^N(\pert)\end{bmatrix}^\top \in \mathbb{R}^N.$ The vector holding all node-level dual bounds is:
$$\mathbf{\dual}_{LP}(\pert) = \left(U\odot[\bar C(\pert)]^- + L\odot[\bar C(\pert)]^+\right)^\top\mathbf 1_n + \bar b(\pert),$$
where $\odot$ denotes the Hadamard product.

It remains to show how to compute $\bar{C}(\pert)$ and $\bar{b}(\pert)$ efficiently for each perturbation. Collect the stored nominal dual multipliers in $$\Lambda = \begin{bmatrix} \dualvar^{1} & \dualvar^{2} & \cdots & \dualvar^{N}
\end{bmatrix} \in \R_+^{m\times N}.$$ The following quantities are independent of $\pert$, and thus can be computed offline and cached for future use: 
$$\bar{C}^0 = c\bm{1}_N^\top - A^\top \Lambda, \text{ and } \bar{b}^0 = \Lambda^\top b.$$
Thus, for any online perturbations, we compute the following using the cached offline information:
$$\bar{C}(\pert) = \bar{C}^0 + \pert[c]\mathbf{1}^\top_N - \pert[A]^\top\Lambda, \text{ and } \bar{b}(\pert) = \bar{b}^0 + \Lambda^\top \pert[b].$$
Notice that the cost of the online computations depends on the type of perturbation: $\pert[c]$ requires only broadcasting, $\pert[b]$ requires a matrix-vector product, and finally $\pert[A]$ introduces matrix-matrix products.

\paragraph{Memory and batching.} Although the expressions above involve matrices of large dimensions ($n \times N$, $m\times N$), their evaluation does not require storing all such matrices densely. In particular, we can exploit the natural sparsity of $\Lambda$ and the high degree of redundant information in $L$ and $U$.

We first remark that $\Lambda$ is naturally sparse, given that we select basic optimal nominal dual multipliers. Indeed, for any primal feasible node $v$ with $F_v$ fixed variables, a basic optimal dual satisfies $\|\lambda^v\|_0 \leq \min\{m,n-F_v\}$. Sparse storage of $\Lambda$ can therefore reduce memory requirements and potentially accelerate subsequent matrix operations.

Considering $L$ and $U$, we point out two  sources of redundant information. First, since branching is restricted to integer variables, the bounds of continuous variables remain static across the tree and thus do not require recording.
{Second, at node $v$ of depth $d(v)$, exactly $d(v)$ variable bounds differ from their original values. Specifically $\|L[:,v]-\ell\|_0+\|u-U[:,v]\|_0 = d(v).$}
We therefore store only the deviation matrices, saving significant memory due to sparse representations. 

Nevertheless, cached information such as $\bar{C}^0$ and intermediate arrays such as $\bar{C}(\pert)$ may still require substantial memory, motivating the use of batching. We finish with a comment on how batching balances parallelism and memory. If evaluating all variable contributions to $\mathbf{f}_{LP}(\pert)$ simultaneously requires excessive memory, we can instead process subsets of $K$ variables at a time, computing and accumulating their contributions to the node-level LP duals.
This limits the peak memory required at the cost of some batching overhead.
Conversely, if evaluating $\mathbf{f}_{LP}(\pert)$ requires relatively little memory, we can exploit additional parallelism by processing multiple perturbations $\pert^1, \pert^2, \ldots \pert^S$ simultaneously, improving GPU utilization and throughput.

\subsection{Parallelizing full-tree propagation}\label{subsec:gpu.propagation}
The generalized MILP dual bounds proposed in this paper are computed using the full-propagation dynamic program \eqref{eq:full.propagation}, where bounds are strengthened by recursing up the nominal tree. This introduces parent-child dependencies, preventing the independent
evaluation of all propagated node bounds, as was possible for the elements of $\mathbf{f}_{LP}(\pert)$. Nevertheless, substantial parallelism remains: all nodes at the same depth can be processed simultaneously, since their propagated bounds depend only on their respective subtrees.

We re-express the full-propagation recursion \eqref{eq:full.propagation} in a vectorized form to make this evident. Let $D$ denote the depth of tree $\TT$, and let $V_d$ denote the set of non-leaf nodes at depth $d$, for $d=1,\ldots,D-1$.
Initializing $\mathbf{BB}(\pert)=\mathbf{f}_{LP}(\pert)$, we update the bounds level by level, starting from the deepest level, where square brackets are used to index vectors and minima/maxima are taken elementwise:
\begin{align}\label{eq:level.by.level.prop}
\mathbf{BB}(\pert)[V_d]
=
\max\left\{
\mathbf{f}_{LP}(\pert)[V_d],
\min\left\{
\mathbf{BB}(\pert)[L(V_d)],
\mathbf{BB}(\pert)[R(V_d)]
\right\}
\right\},
\qquad d=D-1,\ldots,1.
\end{align}

This  yields a procedure in which all nodes at a given depth are processed in parallel, requiring $D-1$ propagation rounds.
This favors shallower nominal trees over unbalanced ones.
The propagation step retains a simple computational structure: the vector $\mathbf{f}_{LP}(\pert)$ is computed only once, and subsequent propagation involves only indexed vector operations and elementwise minima and maxima.

\subsection{Batched variable domain reduction}
We now discuss how the variable fixing and bound tightening methods of Section~\ref{sec:inferring.opt.cond} can be efficiently implemented on GPUs.
These methods retain the same computational structure as the generalized MILP dual: independent node-level evaluations followed by bottom-up tree propagation.

Additional parallelism is exploited across candidate variables and their associated domain reduction tests (for example, fixing variables to zero or one, or tightening a continuous variable's upper or lower bound). Moreover, these computations reuse the node-level dual bounds $\mathbf{f}_{LP}(\pert)$ and reduced costs $\bar{C}(\pert)$, already computed when evaluating the generalized dual.
Consequently, the main additional challenge is managing the memory required to perform many candidate domain reductions simultaneously.

We illustrate this with binary variable fixing. Ignoring incompatible node-level fixings momentarily, the matrices holding $f^v_{LP}(\pert;x_i=0)$ and $f^v_{LP}(\pert;x_i=1)$ for all $v \in V(\TT)$ and $i\in\mathcal{B}$ can  respectively be expressed as:
\begin{align}
    &\mathbf{f}_{LP}^0(\pert) =
        \bm{1}_{|\mathcal B|}\mathbf{f}_{LP}(\pert)^\top - [\bar C(\pert)[\mathcal{B},:]]^-\odot U[\mathcal{B},:] \label{eq:gpu.binary.fixing.0}\\
    &\mathbf{f}^1_{LP}(\pert) = \bm{1}_{|\mathcal B|}\mathbf{f}_{LP}(\pert)^\top +[\bar C(\pert)[\mathcal{B},:]]^+\odot(\mathbf{1}_{|\mathcal B|\times N}-L[\mathcal{B},:]).\label{eq:gpu.binary.fixing.1}
\end{align}
These expressions involve only elementwise operations on shared arrays, allowing candidate restrictions to be evaluated in parallel.
Assignments incompatible with node-level bounds are
identified using Boolean masks and treated as having a bound
of $+\infty$.
The propagated bounds, $BB^{\text{root}}(\pert;x_i=s)$ are obtained by the same level-by-level propagation procedure as \eqref{eq:level.by.level.prop}, parallelized across candidate tests. 
If memory requirements become burdensome, tests are processed in smaller batches.

\paragraph{Boolean propagation for variable fixing.} We conclude with an improvement to variable fixing when a primal bound $\bar{z}_{\pert}$ is available at computation time.
Recall that to certify $x_i^* \neq s$, we only need to determine if $BB^{\text{root}}(\pert;x_i=s)>\bar{z}_{\pert}$ rather than compute the numerical value of the propagated bound. To leverage this, define
$$\tau_{i,s}^v = \mathbf{1}\{ f_{LP}^v(\pert;x_i=s)>\bar z_{\pert}\},$$

The full-tree fixing
certificate can then be propagated using Boolean logic:
$$ \bar\tau_{i,s}^v = \tau_{i,s}^v \lor \left(\bar\tau_{i,s}^{L(v)}\land \bar\tau_{i,s}^{R(v)}\right),
$$
with $\bar\tau_{i,s}^v=\tau_{i,s}^v$ at leaf nodes. If $\bar{\tau}_{i,s}^{\text{root}}=1$, then we can safely determine that $x^*_i\neq s$. This Boolean propagation avoids retaining and propagating the full floating-point matrices $\mathbf{f}^0_{LP}(\pert)$ and $\mathbf{f}^1_{LP}(\pert)$, substantially reducing memory usage.

\subsection{Sparse and structured perturbations}

Thus far, the discussion in this section has considered arbitrary perturbations.
However, computational improvements are available when perturbations are sparse or exhibit specific structure.
Sparse perturbations naturally arise in ``what-if'' analyses, where one wishes to evaluate the effect of changing a small subset of the problem data.
Structured perturbations occur when changes are governed by a small number of parameters or are confined to particular rows or right-hand sides.

\paragraph{Sparse perturbations.} The reduced costs $\bar{c}_i^v(\pert)$ can differ from their nominal values only if $\pert[c]$ or $\pert[A]$ affect variable $i$. Define the set of potentially affected variables as:
$$\mathcal{J}(\pert)=\mathrm{supp}(\pert[c])\cup \left\{i\in[n]: \pert[A^i]\neq 0\right\}.$$
Since all variables outside of $\mathcal{J}(\pert)$ have unchanged reduced costs, we can rewrite generalized node-level duals in terms of nominal values:
\begin{align}\label{eq:gpu.incremental.dual}
\dual_{LP}^v(\pert)
=
\dual_{LP}^v(0) +(\dualvar^v)^\top\pert[b]  + \sum_{i\in \mathcal{J}(\pert)}
\left(
\left(\left[\bar c_i^v(\pert)\right]^+-\left[\bar c_i^v(0)\right]^+\right)\ell_i^v +  \left(\left[\bar c_i^v(\pert)\right]^--\left[\bar c_i^v(0)\right]^-\right)u_i^v \right).
\end{align}
The key observation is that when $|\mathcal{J}(\pert)| \ll n$, we can avoid recomputing the contributions of most variables.
By caching $\mathbf{f}_{LP}(0)$ and $\bar{C}^0$, we only need to compute the updated reduced costs $\bar{c}_i^v (\pert)$ for $i \in \mathcal{J}(\pert)$.
These computations remain vectorized, but now operate on substantially smaller matrices, reducing both computation and memory requirements.
Although incremental evaluation remains vectorized and operates on much smaller matrices, it introduces additional indexing and potentially noncontiguous memory accesses, making this approach less effective than direct evaluation \eqref{eq:LP.dual.dense.gpu} when $|\mathcal{J}(\pert)| \approx n$.

\paragraph{Structure.} The structure of a perturbation may be exploited, even if it is dense. We highlight two cases:
\begin{itemize}
    \item \textbf{Pure right-hand side perturbations:} If $\pert = (\mathbf{0}, \pert[b], \mathbf{0})$, then $\mathcal{J}(\pert)=\emptyset$ and node-level dual bounds simplify to
    $\mathbf{f}_{LP}(\pert) = \mathbf{f}_{LP}(0) + \Lambda^\top\pert[b]$.
    \item  \textbf{Low-rank matrix perturbations:} Suppose we wish to perturb only a single row $r$ of $A$, that is $\pert[A] = e_r \pert[a]^\top$, where $e_r$ is the $r$-th unit vector and $\pert[a] \in \mathbb{R}^n.$ Then $\bar{C}(\pert) = \bar C^0 + \pert[c]\mathbf{1}_N^\top - \pert[a]\Lambda[r,:],$ requiring only an outer product to compute the $\pert[A]$ contribution to the reduced cost. More generally, low-rank perturbations can be handled through products involving their low-rank factors.
\end{itemize}

\section{Experimental setup}\label{sec:numerical}
This section describes the test problems and the data we perturb (Section~\ref{subsec:problems}), the nominal instances (Section~\ref{subsec:nominal.instances}), how they are solved and selected (Section~\ref{subsec:selection}), and how the perturbed instances are generated (Section~\ref{subsec:perturbation}).

\subsection{Test problems}\label{subsec:problems}
In this subsection, we provide a brief description of the formulations used in our experiments and the parameters that were perturbed for each problem.

\paragraph{Generalized Independent Set}
Given a graph $G = (V,E)$, where edges are partitioned $E = E_1 \cup E_2$, into non-removable edges $E_1$ and removable edges $E_2$. The GISP with edge costs $c_{ij}$ and node rewards $w_i$ is formulated by:
\begin{align*}
    \min ~ & \sum_{(i,j) \in E_2} (c_{ij}+\pert[c_{ij}])y_{ij} - \sum_{i \in V} (w_i+\pert[w_i]) x_i\\
    \text{s.t. } &x_i + x_j \leq 1, \quad \forall (i,j) \in E_1\\
    & x_i + x_j - y_{ij} \leq 1,\quad  \forall (i,j) \in E_2\\
    & x_i \in \{0,1\}, \quad \forall i \in V\\
    & y_{ij} \in \{0,1\}, \quad  \forall (i,j) \in E_2.
\end{align*}
For this problem, we consider perturbations of the form $\pert=(0,0,\pert[c])$ given by $\pert[w_i]$ and $\pert[c_{i,j}]$, to capture changes in the economic trade-offs governing the selection of nodes and the removal of conflicting edges.

\paragraph{Combinatorial Auction}
Given a set of items $\mathcal{S}$ and a collection of bids $\mathcal{B}$, where each bid $b \in \mathcal{B}$ requests a subset of items $\mathcal{S}_b \subseteq \mathcal{S}$ at a nominal price $p_b$, let $A_{ib} = 1$ if $i \in \mathcal{S}_b$ and $A_{ib} = 0$ otherwise. The CA problem is formulated as:
\begin{align*}
    \max ~ & \sum_{b \in \mathcal{B}} (p_b + \pert[p_b])\, x_b\\
    \text{s.t. } & \sum_{b \in \mathcal{B}} (A_{ib} + \pert[A_{ib}])\, x_b \leq 1, \quad \forall i \in \mathcal{S}\\
    & x_b \in \{0,1\}, \quad \forall b \in \mathcal{B}.
\end{align*}
We perturb the bid prices, $\pert[p_b] \in \R$ for $b \in \mathcal{B}$, and the bid contents, $\pert[A_{ib}] \in \{0, 1 - 2A_{ib}\}$ for $i \in \mathcal{S}$ and $b \in \mathcal{B}$, giving $\pert = (\pert[A], 0, \pert[c])$. \textcolor{black}{These capture variations in bidder valuations and in the composition of bids.}

\paragraph{Capacitated Facility Location}
Given a set of customers $C$ and a set of facilities $F$, let $w_f^s$ denote the opening cost of facility $f$, $w_{c,f}^{t}$ the transportation cost of assigning customer $c$ to facility $f$, $d_c$ the demand of customer $c$, $b_f$ the capacity of facility $f$, and $B$ the minimum total capacity. We define CFLP by:
\begin{align*}
\min_{x,y} ~ & \sum_{c\in C}\sum_{f\in F}(w^{t}_{c,f} + \pert[w^{t}_{c,f}]) x_{c,f} + \sum_{f\in F}(w_f^s + \pert[w_f^s]) y_f \\
\text{s.t. } & \sum_{f\in F}x_{c,f} \geq 1, \quad \forall c\in C  \\
     & \sum_{c\in C}(d_c + \pert[d_c]) x_{c,f} \leq b_f y_f, \quad  \forall f \in F \\
     & \sum_{f\in F}b_f y_f \geq B \\
     & x_{c,f} \leq y_f, \quad \forall c\in C, ~ \forall f\in F\\
     & 0 \leq x_{c,f} \leq 1, \quad \forall c\in C, ~ \forall f\in F\\
     & y_f \in \{0,1\}, \quad  \forall f\in F.
\end{align*}
We perturb customer demands, facility opening costs, and transportation costs, $\pert= (\pert[A], 0, \pert[c])$ given by $\pert[d_c]$, $\pert[w_f^s]$ and $\pert[w_{c,f}^{t}]$, to capture variations in capacity requirements and in the costs underlying location and assignment decisions.

\paragraph{Unit Commitment} Consider $G$ generators, indexed by $g = 1, \ldots, G$, and time steps $t = 1, \ldots, T$. Let $C^F_g$ be the fixed cost, $C^S_g$ the start-up cost and $C^V_g$ the variable cost of generator $g$, $P^{\min}_g$ and $P^{\max}_g$ its minimum and maximum output, $L_g$ and $l_g$ its minimum up and down times, and $D_t$ the demand at time $t$. The variables $x_{gt}$, $w_{gt}$ and $p_{gt}$ indicate whether generator $g$ is on, whether it starts up, and its output at time $t$, and $x_{g0}$ is the given initial status. We define the UC problem by:
\begin{align*}
\min \;\;\;
    & \sum_{t=1}^T \sum_{g=1}^G \left(
        x_{gt} \left( C^F_g +\Delta_{C_g^F}\right)
        + w_{gt} \left(C^S_g+ \Delta_{C_g^S} \right)
        + p_{gt} \left(C^V_g + \Delta_{C_g^V}\right)
    \right)
    \\
\text{s.t.} \;\;\;
    & \sum_{k=t-L_g+1}^t w_{gk} \leq x_{gt}
        & \forall g,\ \forall t=L_g,\ldots,T \\
    & \sum_{k=t-l_g+1}^{t} w_{gk} \leq 1 - x_{g,t-l_g}
        & \forall g,\ \forall t=l_g,\ldots,T \\
    & w_{gt} \geq x_{gt} - x_{g,t-1}
        & \forall g,\ \forall t=1,\ldots,T \\
    & \sum_{g=1}^G p_{gt} \geq D_t + \Delta_{D_t}
        & \forall t \\
    & P^{\min}_g x_{gt} \leq p_{gt} \leq P^{\max}_g x_{gt}
        & \forall g, t \\
    & x_{gt},\, w_{gt} \in \{0, 1\}
        & \forall g, t.
\end{align*}
We perturb electricity demand and the fixed, start-up and variable generation costs, $\Delta=(0, \Delta_b, \Delta_c)$ given by $\Delta_{D_t}$, $\Delta_{C^F_g}$, $\Delta_{C^S_g}$ and $\Delta_{C^V_g}$, to capture variations in system requirements and operating costs that influence commitment and dispatch decisions. 

Table~\ref{tab:perturbations} summarizes which parts of the data are perturbed in each problem. Only the data named above are perturbed; for example, matrix perturbations in CFLP change only the demands $d_c$.
\begin{table}[htbp]
  \centering
  \caption{Data perturbed for each problem}
    \begin{tabular}{cccc}
    \hline
    \textbf{Problem} & \textbf{Matrix} & \textbf{RHS} & \textbf{Objective} \bigstrut\\
    \hline
    GISP &    &    & $\checkmark$ \bigstrut[t]\\
    CA & $\checkmark$   &    & $\checkmark$  \\
    CFLP & $\checkmark$   &    & $\checkmark$  \\
    UC &    & $\checkmark$   & $\checkmark$  \bigstrut[b]\\
    \hline
    \end{tabular}%
  \label{tab:perturbations}%
\end{table}%

\subsection{Nominal instances}\label{subsec:nominal.instances}
We use publicly available instances from Distributional MIPLIB \cite{huang2024distributional}, specifically GISP-medium \cite{ferber2022backdoors,hochbaum1997gisp}, CA-medium \cite{leyton2000towards} and CFLP-medium \cite{cornuejols1991comparison}. For UC, we generate synthetic instances using \texttt{UnitCommitmentGenerator} from \texttt{MIPLearn} 0.4 \cite{SantosXavier2025Jun} with random seed $s=42$. For each number of generators $G\in\{25,\dots,75\}$ and $T=24$ time steps, we generate 5 instances with $P_g^{\max}\sim \mathcal{U}[200,500]$, $P_g^{\min}=\alpha P_g^{\max}$ with $\alpha\sim \mathcal{U}[0.8, 0.9]$, $C_g^{F}\sim \mathcal{U}[0,1000]$, $C_g^{S}\sim \mathcal{U}[0,10000]$, $C_g^{V}\sim \mathcal{U}[0,50]$, $L_g\sim \mathcal{U}\{2,8\}$ and $l_g\sim\mathcal{U}\{2,8\}$. All costs and demands are further randomized by setting \texttt{cost\_jitter} $\sim \mathcal{U}[0.75, 1.25]$ and \texttt{demand\_jitter} $\sim \mathcal{U}[10,10.01]$.

\subsection{Nominal solves and instance selection}\label{subsec:selection} 
To obtain the nominal search trees, we use a custom pure branch-and-bound algorithm \cite{Dey2024May,Shah2025Jul} written in Python 3.11 that uses Gurobi 12.0.3 as the LP solver, with \textit{reliability branching} and reliability parameter $8$. Since no cuts or preprocessing are applied, the resulting trees are pure branch-and-bound trees, as required by Proposition~\ref{valid.dual.ip}. At every node, we store the basic optimal dual solution if the node is primal feasible, the multiplier of Section~\ref{sec:dualselection} if it is primal infeasible, the branching fixings, and the position of the node in the tree. The nominal solves are run on the ISyE HPC cluster with 16 GB of RAM and 1 CPU per instance.

We solve each instance with a time limit of 12 hours and discard instances that are not solved to optimality. For each problem, we sort the remaining instances by the number of nodes in their optimal trees and split them into three groups of equal size: \textit{small} (first third), \textit{medium} (second third) and \textit{large} (last third). From each group, we randomly sample 5 instances with random seed $s=42$, giving the 15 instances per problem reported in Appendix~\ref{appendix:instance.selection}. Table~\ref{tab:data-selection}  below reports the number of instances at each step. Table~\ref{tab:selected.instances} summarizes the selected instances: within each problem class, the size groups have similar numbers of binary variables but differ mainly in the size of the nominal tree, which ranges from about $600$ nodes for small CFLP instances to more than $500{,}000$ nodes for large GISP instances.

\begin{table}[htbp]
  \centering
    \caption{Selected nominal instances: total number of variables, number of binary variables, number of constraints, and number of nodes and depth of the nominal branch-and-bound tree, averaged over the 5 instances of each size group.}

    \begin{tabular}{ccccccc}
    \toprule
    \textbf{Problem} & \textbf{Group} & \textbf{Variables} & \textbf{Bin. Variables} & \textbf{Constraints} & \textbf{Nodes} & \textbf{Depth} \\
    \midrule
       & small & 970 & 970 & 3339 & 345400 & 108 \\
    GISP & medium & 981 & 981 & 3380 & 410047 & 109 \\
       & large & 992 & 992 & 3363 & 527361 & 109 \\
    \midrule
       & small & 1500 & 1500 & 579 & 4834 & 25 \\
    CA & medium & 1500 & 1500 & 579 & 17078 & 30 \\
       & large & 1500 & 1500 & 582 & 48870 & 33 \\
    \midrule
       & small & 20200 & 200 & 20301 & 582 & 88 \\
    CFLP & medium & 20200 & 200 & 20301 & 3073 & 93 \\
       & large & 20200 & 200 & 20301 & 33199 & 98 \\
    \midrule
       & small & 2016 & 1344 & 3150 & 21705 & 123 \\
    UC & medium & 2088 & 1392 & 3275 & 71116 & 117 \\
       & large & 2203 & 1469 & 3453 & 200965 & 131 \\
    \bottomrule
    \end{tabular}%
  \label{tab:selected.instances}%
\end{table}%

\begin{table}[htbp]
  \centering
  \caption{Number of instances at each step}
    \begin{tabular}{cccc}
    \hline
    \textbf{Problem} & \textbf{Initial pool} & \textbf{Solved to optimality} & \textbf{Selected} \bigstrut\\
    \hline
    GISP & 1000 & 113 & 15 \bigstrut[t]\\
    CA & 1000 & 921 & 15 \\
    CFLP & 1000 & 549 & 15 \\
    UC & 250 & 68 & 15 \bigstrut[b]\\
    \hline
    \end{tabular}%
  
  \label{tab:data-selection}%
\end{table}%

\subsection{Data perturbation}\label{subsec:perturbation}
We perturb the selected instances in two ways, which correspond to two uses of sensitivity analysis of mixed-integer programs. The first is \emph{what-if analysis}: a decision maker asks how the optimal value or the optimal decisions would change if a single input were different, for example, if one cost increased or one demand grew. Such questions arise when explaining and stress-testing a solution; for instance, the OptiGuide framework deployed in Microsoft's cloud supply chain answers what-if queries from planners by re-solving the underlying MILP for each query~\cite{li2023large}. We model these questions by \emph{single-coefficient changes}, chosen so that the change works against the nominal optimal solution. 

The second use is \emph{operational reoptimization}: problems such as unit commitment are solved repeatedly as new data arrive, and consecutive instances differ in many coefficients at once, each by a small amount. We model this setting by \emph{{multi-coefficient changes}}, which perturb a fraction of the data simultaneously. Throughout, $\rho \in \{1\%, 5\%, 10\%\}$ denotes the magnitude of a relative change, and only the data listed in Table~\ref{tab:perturbations} are perturbed.

\subsubsection{Single-parameter change}\label{subsec:single.change.perturbations}
{This setting addresses ``what if?'' questions, in which a single parameter of the problem changes. For each selected nominal instance, each parameter family subject to perturbations, and each magnitude $\rho \in \{1\%,5\%, 10\%\}$, we generate four scenarios: two increases and two decreases. An increase or decrease replaces the selected parameter $\theta$ by $(1+\rho)\theta$ or $(1-\rho)\theta$, respectively.}

Let $x^*$ denote the optimal solution of the nominal instance found by branch-and-bound.
Changes that favor $x^*$, namely improving the objective coefficient of a variable that is positive in $x^*$ or worsening the objective coefficient of a variable that is zero in $x^*$, leave $x^*$ optimal, so the new optimal value follows directly from $x^*$.
We therefore change objective coefficients against $x^*$: we make a variable that is positive in $x^*$ less attractive (increase its cost or decrease its reward), or make a variable that is zero in $x^*$ more attractive (decrease its cost or increase its reward). Two scenarios of each kind are generated. For the same reason, when a demand is increased, the increase is chosen large enough to make $x^*$ infeasible, since otherwise $x^*$ remains optimal. The problem-specific rules are:

\begin{itemize}
    \item \textbf{GISP}: we consider perturbing node rewards and removable-edge costs separately. For node rewards, two scenarios increase $w_i$ for a node with $x_i^*=0$, and two decrease $w_i$ for a node with $x_i^*=1$. For edge costs, two scenarios increase $c_{ij}$ for an edge with $y_{ij}^*=1$, and two decrease $c_{ij}$ for an edge with $y_{ij}^*=0$.
    \item \textbf{CA}: for bid prices, two scenarios increase $p_b$ for an unselected bid, $x_b^*=0$, and two decrease $p_b$ for a selected bid,
    $x_b^*=1$. Bid-content changes are binary flips and therefore have no relative magnitude $\rho$. We generate four such scenarios per nominal instance, all involving selected bids. Two scenarios remove an item $i$ from a selected bid $b$ with $x_b^*=1$, changing $A_{ib}$ from $1$ to $0$. Two further scenarios add an item $i$ to a selected bid $b$ with $x_b^* = 1$, changing $A_{ib}$ from $0$ to $1$. For additions, the item is already allocated to another selected bid, so the nominal solution becomes infeasible.
    \item \textbf{CFLP}: cost perturbations affect facility-opening costs $w_f^s$. Two scenarios increase $w_f^s$ for an open facility, $y_f^*=1$, and two decrease $w_f^s$ for a closed facility, $y_f^*=0$. For demand increases, define the nominal capacity slack of facility $f$ by $s_f = b_f y_f^*-\sum_{c\in C}d_c x_{c,f}^*$, and, for each customer $c$, define $\tau_c =
        \min\left\{
            \frac{s_f}{x_{c,f}^*} :
            f\in F,\ x_{c,f}^*>0
        \right\}.$
    We generate two scenarios by selecting a customer satisfying
    $\rho d_c>\tau_c$ and increasing its demand to
    $(1+\rho)d_c$. This makes the nominal solution infeasible at at least one facility. Two further scenarios select a customer without
    this restriction and decrease its demand to $(1-\rho)d_c$. Each demand perturbation changes a single customer parameter $d_c$, which appears in every facility-capacity constraint; it therefore changes multiple constraint-matrix coefficients simultaneously, despite being a single-parameter perturbation.
    \item \textbf{UC}: for cost perturbations, we change the on/off-cost coefficient of a single generator--period pair $(g,t)$. For each magnitude $\rho$, two scenarios increase this coefficient for a pair with $x_{gt}^*=1$, and two decrease it for a pair with $x_{gt}^*=0$. For demand perturbations, let $s_t=\sum_g p_{gt}^*-D_t$ denote the nominal demand slack at period $t$. For each magnitude $\rho$, two scenarios select a period satisfying $\rho D_t>s_t$ and increase its demand to $(1+\rho)D_t$, making $x^*$ infeasible. Two further scenarios select a period without this restriction and decrease its demand to $(1-\rho)D_t$.
\end{itemize}

\subsubsection{Multi-parameter changes}
This setting studies nearby instances by changing multiple parameters simultaneously. A sparsity level specifies the fraction of eligible parameters selected for perturbation. For each selected scalar parameter $\theta$, a magnitude $\rho$ determines the additive change:
$$\pert[\theta] \sim \mathcal{U}[-\rho|\theta|,\rho|\theta|].$$
We generate perturbations separately for each dataset as follows:
\begin{itemize}
    \item \textbf{GISP:} we perturb node rewards and removable-edge costs. We select a fraction $S_c \in \{1\%,5\%,10\%\}$ of these coefficients, perturbing each selected one with magnitude $\rho_c \in \{1\%,5\%,10\%\}.$
    \item \textbf{CA:} we simultaneously perturb bid prices and bid contents. For bid prices, we use sparsity $S_c \in \{1\%,5\%,10\%\}$ and magnitude $\rho_c \in \{1\%,5\%,10\%\}$. For bid contents, we select a fraction $S_A \in \{0.001\%, 0.005\%, 0.01\%\}$ of the item--bid incidence coefficients, and flip each selected coefficient from $0$ to $1$, or from $1$ to $0$. These sparsity levels correspond to $9$, $44$ and $87$ total flips per instance, respectively.
    \item \textbf{CFLP:} we simultaneously perturb cost coefficients and customer demands. For costs, we use sparsity $S_c \in \{1\%,5\%,10\%\}$ and magnitude $\rho_c \in \{1\%,5\%,10\%\}$. For customer demands, we also select sparsity $S_A \in \{1\%,5\%,10\%\}$ and magnitude $\rho_A \in \{1\%,5\%,10\%\}$. Demand sparsity is measured over customer demand parameters, rather than their repeated entries in the constraint matrix.
    \item \textbf{UC:} we simultaneously perturb operating-cost coefficients and electricity demands. For costs, we use sparsity $S_c \in \{1\%,5\%,10\%\}$ and magnitude $\rho_c \in \{1\%,5\%,10\%\}$. For demands, we select $k_b \in \{1,2,3\}$ of the $24$ time periods and perturb their demands with magnitude $\rho_b \in \{1\%,5\%,10\%\}$. The demand perturbations correspond to sparsity levels of approximately $\{4\%,8\%,12\%\}.$
\end{itemize}
For each choice of sparsity and magnitude levels, we generate 3 perturbed instances with different random seeds for the selected indices and magnitudes. This gives $15 \times 9 \times 3 = 405$ perturbed instances for GISP, $15 \times 27 \times 3 = 1215$ for CA, and $15 \times 81 \times 3 = 3645$ for each of CFLP and UC.

\section{Numerical results: What-if analysis (single-parameter perturbations)}\label{sec:whatif}
In this section, we present numerical results of our framework to ``what-if'' analysis.

As discussed in the previous section, in this setting, a nominal MILP has already been solved, and we wish to understand how its optimal objective value and solution change under a single perturbation to the problem data.
We investigate three applications of the proposed GPU-accelerated offline-online framework: efficiently obtaining dual bounds for perturbed problems (Section~\ref{subsec:whatif.duals}), certifying partial optimal solutions (Section~\ref{subsec:whatif.var.fixing}), and accelerating reoptimization (Section~\ref{subsec:whatif.reopt}).

We assume a practitioner submits multiple what-if requests, each corresponding to changes in a single parameter of the problem data.
Our main computational comparison is between evaluating these requests using the stored nominal information and repeatedly solving the corresponding perturbed problems with a commercial solver.

\paragraph{Experiment setup}
We consider the four problem classes detailed in Section~\ref{sec:numerical}, with 15 nominal instances per class. Since each perturbation changes a single parameter, the nominal optimal solution sometimes trivially remains optimal. For example, for a minimization problem, decreasing the objective coefficient of a variable at its upper bound in the nominal optimal solution trivially preserves its optimality. Similarly, restricting the feasible region preserves nominal optimality whenever the nominal optimal solution remains feasible.
As discussed in Section~\ref{subsec:single.change.perturbations}, we therefore generate perturbations based on the nominal optimal solution of each instance and exclude such trivial directions. {The sampling rules are given in Section~\ref{subsec:single.change.perturbations}.}

GPU computations are performed on a single NVIDIA H100 GPU, while the comparison with a commercial solver uses Gurobi 13.0.3 on an Intel Core i9-12900KF processor with 24 threads.

\subsection{Dual bound gaps}\label{subsec:whatif.duals}
We first evaluate the quality of the generalized dual bounds, and compare their evaluation time with the time Gurobi requires to obtain equally strong bounds.
We then investigate sources of dual bound degradation, and the extent to which bound propagation can compress the nominal tree.

\subsubsection{Generalized dual optimality gap}\label{subsubsec:whatif.dual.quality}
\begin{figure}[H]
\centering\includegraphics[width=\linewidth]{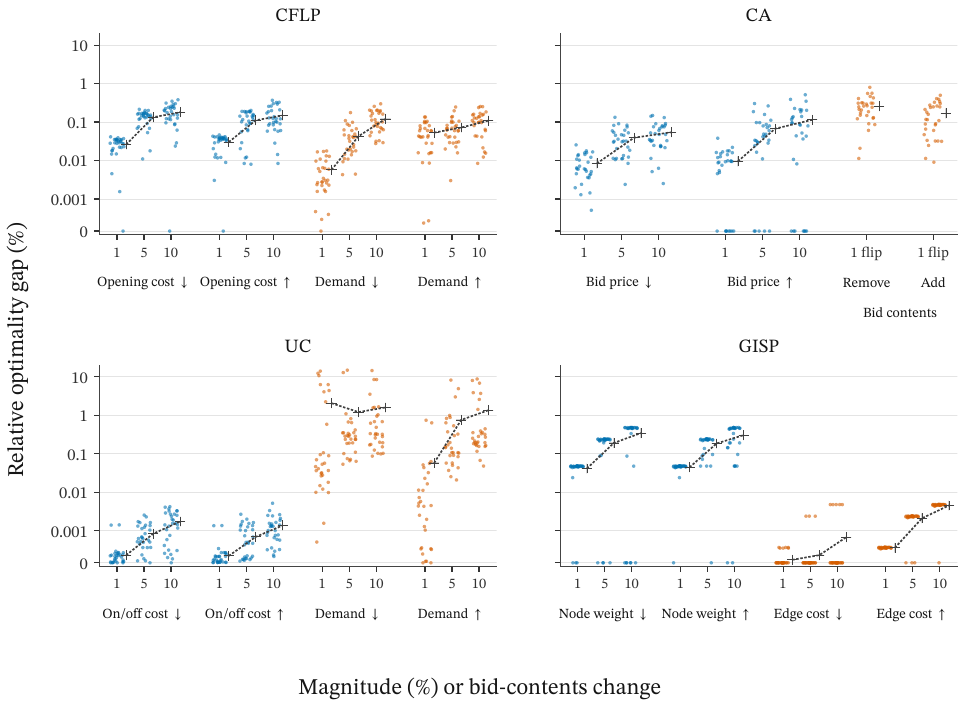}
\caption{Relative optimality gap for sampled single-parameter perturbations. Each point represents a single outcome.}
\label{fig:whatif.dual.individual}
\end{figure}

For each perturbed instance, we report the relative optimality gap of the generalized dual bound $|z(\pert)-f(\pert)|/|z(\pert)|.$
Tables~\ref{tab:whatif.obj.change.mean.dualgap} and \ref{tab:whatif.feas.region.change.mean.dualgap} summarize the mean gaps by problem class, perturbation type, and magnitude.
Figure~\ref{fig:whatif.dual.individual} provides a more detailed view of the individual results, illustrating the variability across instances and perturbation scenarios.

The generalized dual bound is tight in almost all settings: the mean gap is at most $0.35\%$ for every objective perturbation and at most $0.26\%$ for the CFLP demand and CA {bit contents perturbations}, and the mean gap is at most $0.005\%$ for UC on/off costs and GISP edge costs, with many scenarios attaining a zero gap. The mean gap typically increases with the magnitude $\rho$, and Figure~\ref{fig:whatif.dual.individual} shows that this trend holds across individual instances rather than being driven by a few outliers. The main exception is UC demand perturbations, the only setting with mean gaps above $1\%$ and individual gaps exceeding $10\%$, where demand decreases also break the monotone trend, with the largest mean gap at $\rho = 1\%$.

\begin{table}[t]
  \centering
  \caption{{Mean relative optimality gap} (\%) under single-parameter objective perturbations.}
  \label{tab:whatif.obj.change.mean.dualgap}
  \begingroup
  \small
  \setlength{\tabcolsep}{5pt}
  \renewcommand{\arraystretch}{1.12}
  \begin{tabular}{lrrrrrr}
    \toprule
    & \multicolumn{3}{c}{Increase} & \multicolumn{3}{c}{Decrease} \\
    \cmidrule(lr){2-4}\cmidrule(lr){5-7}
    Problem class / parameter & 1\% & 5\% & 10\% & 1\% & 5\% & 10\% \\
    \midrule
    CFLP: {opening cost} & 0.030 & 0.11 & 0.15 & 0.026 & 0.13 & 0.18 \\
    CA: bid price & 0.001 & 0.067 & 0.12 & 0.008 & 0.04 & 0.05 \\
    UC: on/off cost & $2\times 10^{-4}$ & $8\times 10^{-4}$ & 0.001 & $2\times 10^{-4}$ & $9\times 10^{-4}$ & 0.002 \\
    GISP: node weight & 0.045 & 0.19 & 0.30 & 0.042 & 0.19 & 0.35 \\
    GISP: {edge cost} & $5\times 10^{-4}$ & 0.002 & 0.005 & $9\times 10^{-5}$ & $2\times 10^{-4}$ & $8\times 10^{-4}$ \\
    \bottomrule
  \end{tabular}
  \endgroup
  \par\smallskip

\end{table}

\begin{table}[t]
  \centering
  \caption{{Mean relative optimality gap} (\%) under {single-parameter} feasible-region perturbations.}
\label{tab:whatif.feas.region.change.mean.dualgap}
  \begingroup
  \small
  \setlength{\tabcolsep}{5pt}
  \renewcommand{\arraystretch}{1.12}
  \begin{tabular}{lrrrrrr}
    \toprule
    & \multicolumn{3}{c}{Increase} & \multicolumn{3}{c}{Decrease} \\
    \cmidrule(lr){2-4}\cmidrule(lr){5-7}
    Problem class / parameter & 1\% & 5\% & 10\% & 1\% & 5\% & 10\% \\
    \midrule
    CFLP: demand ($A$) & 0.0532 & 0.0720 & 0.109 & 0.00573 & 0.0429 & 0.121 \\
    UC: {demand} ($b$) & 0.0578 & 0.762 & 1.35 & 2.07 & 1.21 & 1.60 \\
    \midrule
    & \multicolumn{3}{c}{Relax (remove)} & \multicolumn{3}{c}{Tighten (add)} \\
    \cmidrule(lr){2-4}\cmidrule(lr){5-7}
    CA: {bid contents} ($A$) & \multicolumn{3}{c}{0.260} & \multicolumn{3}{c}{0.168} \\
    \bottomrule
  \end{tabular}
  \endgroup
  \par\smallskip
\end{table}

\subsubsection{Dual bound speedups}\label{subsubsec:whatif.dual.speedup}
We now compare runtimes to compute generalized dual bounds, comparing the proposed approach with Gurobi.
GPU timings include bound evaluation and transferring the resulting bounds to the CPU, but exclude the offline cost of loading the nominal tree and associated information onto the GPU.
For each perturbation, we measure the time until Gurobi's dual bound first reaches or exceeds $f(\pert)$, the dual bound provided by the GPU computation.
To fairly leverage information from the nominal solve, we provide Gurobi with the nominal optimal solution as a candidate solution, which we observe improves its performance.
Excluding the offline cost is consistent with the intended use of the framework, in which the nominal tree is loaded once and then reused across many what-if requests, so that this cost is amortized. Stopping Gurobi as soon as its dual bound matches $f(\pert)$ measures the time it needs to certify a bound of the same quality, which favors Gurobi, since it is not required to solve the perturbed problem to optimality.

\begin{figure}[H]
\centering\includegraphics[width=0.52\linewidth]{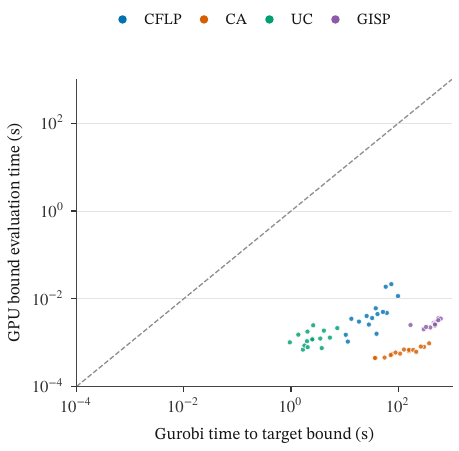}
    \caption{Total time over the perturbations of each nominal instance: time for Gurobi to reach the generalized dual bound $f(\pert)$ (horizontal axis) and time for the GPU to compute it (vertical axis), both in seconds on logarithmic scales. Each point is one nominal instance; the dashed line is the diagonal.}
    \label{fig:whatif.dual.time.individual}
\end{figure}

Table~\ref{tab:whatif.dual.times} reports the total runtimes across all perturbations for each problem class. Note that the Gurobi times are in seconds and the GPU times are in milliseconds.
Figure~\ref{fig:whatif.dual.time.individual} provides a more detailed comparison at the instance level.
The GPU evaluation is between about $2{,}000$ and $240{,}000$ times faster than Gurobi, from UC to CA, and requires well \emph{under a millisecond per perturbation} in every problem class.

\begin{table}[H]
  \centering
  \caption{Total runtime for obtaining the generalized dual bound for all perturbed problems in each problem class.}
  \label{tab:whatif.dual.times}
  \begingroup
  \small
  \setlength{\tabcolsep}{5pt}
  \begin{tabular}{lrrrr}
    \toprule
    Problem class & \# Perturbed problems
      & Gurobi (s)
      & GPU (ms)\\
    \midrule
    CFLP & 360 & 602                    & 93.0 \\
    CA   & 240 & $2.26\times 10^{3}$ & 9.47  \\
    UC   & 360 & 43.3                   & 19.9 \\
    GISP & 360 & $6.78\times 10^{3}$ & 40.7 \\
    
    \bottomrule
  \end{tabular}
  \endgroup
  \par\smallskip
\end{table}

\subsubsection{Dual gap decomposition}\label{subsubsec:whatif.dual.degradation}

The gap between the generalized dual bound $f(\pert)$ and the optimal value $z(\pert)$ can arise from two sources: the nominal branch-and-bound tree may not be sufficient to solve the perturbed instance, and the stored nominal LP dual multipliers may yield weak approximations of the perturbed node LP bounds. We now investigate the relative contributions of these two sources of error.
Let $g(\pert)$ denote the bound obtained by reoptimizing the LP relaxations at the leaves of the nominal tree:
$$g(\pert) = \min_{v\in\LL(\TT)}z^v_{LP}(\pert).$$
Unlike $f(\pert)$, which uses stored nominal dual multipliers, $g(\pert)$ uses optimal dual multipliers.
Since $f(\pert)\leq g(\pert) \leq z(\pert)$, we can decompose the absolute gap as:
\begin{equation}\label{eq:abs.gap.breakdown}z(\pert) - f(\pert) = \underbrace{g(\pert) - f(\pert)}_{\text{nominal dual reuse gap}}+\underbrace{z(\pert)-g(\pert)}_{\text{nominal tree reuse gap}}.
\end{equation}
The first term measures the gap improvement one could attain via access to perturbed dual optimal solutions, and the second term captures the residual gap coming from using a tree that does not solve the perturbed problem.

\begin{figure}[H]
    \centering
    \includegraphics[width=\linewidth]{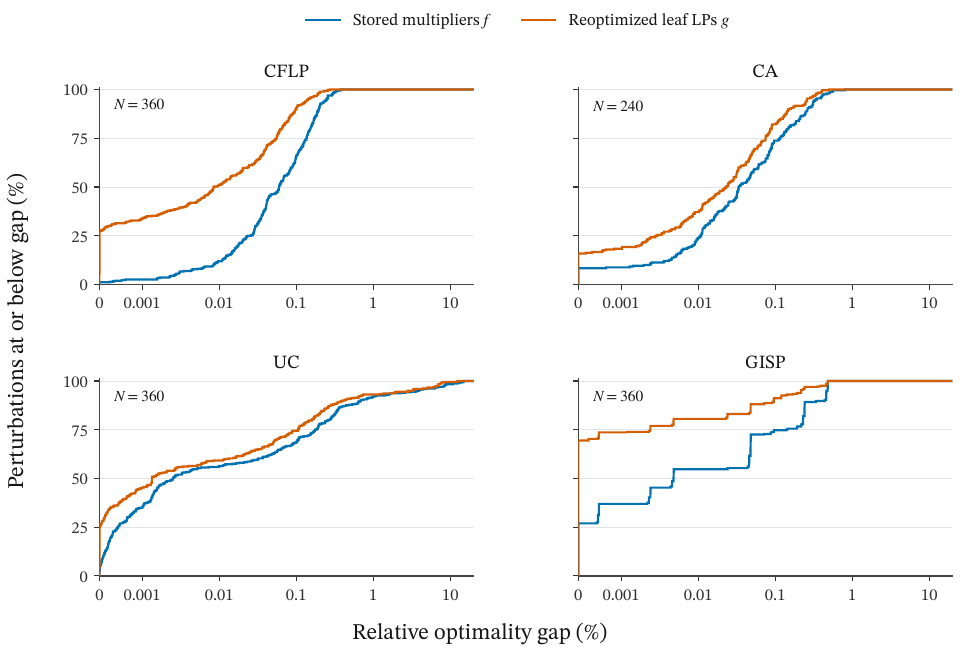}
        \caption{Relative optimality gaps of the generalized dual bound $f(\pert)$ and of the bound $g(\pert)$ obtained by reoptimizing the leaf LPs of the nominal tree. For each value $x$, the curves give the percentage of perturbations whose gap is at most $x\%$; $N$ is the number of perturbations of each problem class.}
        \label{fig:whatif.dual.gap.dist}
\end{figure}

\begin{figure}[H]
\includegraphics[width=\linewidth]{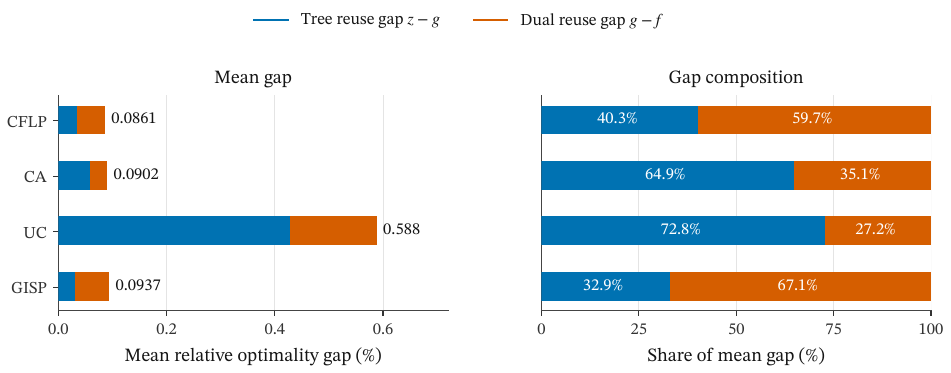}
\caption{Decomposition~\eqref{eq:abs.gap.breakdown} of the gap into the dual reuse gap $g(\pert) - f(\pert)$ and the tree reuse gap $z(\pert) - g(\pert)$. Left: mean relative gap of each term for each problem class; the number gives the mean total gap. Right: each term as a fraction of the total gap $z(\pert) - f(\pert)$.}
        \label{fig:whatif.dual.gap.breakdown}
\end{figure}

Figure~\ref{fig:whatif.dual.gap.dist} reports the gap improvement obtained by reoptimizing the leaf LPs,  while Figure~\ref{fig:whatif.dual.gap.breakdown} examines the relative contributions of the two error terms in the decomposition~\eqref{eq:abs.gap.breakdown}.
Figure~\ref{fig:whatif.dual.gap.breakdown} reports, for each problem class, the dual reuse gap and the tree reuse gap as fractions of the total gap: the fraction of the gap due to reusing the nominal multipliers is $(g(\pert) - f(\pert))/(z(\pert) - f(\pert))$, and the fraction due to reusing the nominal tree is $(z(\pert) - g(\pert))/(z(\pert) - f(\pert))$.

The above analysis separates the problem classes into two groups. For GISP and CFLP, reusing the nominal multipliers accounts for most of the mean gap ($67.1\%$ and $59.7\%$, respectively). Reoptimizing the leaf LPs increases the proportion of scenarios with zero gap from approximately $27\%$ to $70\%$ for GISP and from $2\%$ to $27\%$ for CFLP. For CA and UC, the residual gap dominates ($64.9\%$ and $72.8\%$), so that even exact leaf LPs would leave most of the gap, indicating that the nominal tree itself no longer branches enough for the perturbed problem.

\subsubsection{Bound propagation tree utilization}\label{subsubsec:whatif.dual.tree.compression}

Recall from Section~\ref{subsec:bound.enhancement.techniques} that the generalized dual is obtained by propagating node-level dual bounds up the nominal tree: whenever the local bound of a node is at least the best bound attainable in its subtree, the subtree can be pruned without weakening the final bound. The result is a smaller tree certificate, determined online for each perturbation. We now investigate how much of the nominal branch-and-bound tree is needed by the perturbed instances, measuring the size of the certificate as a fraction of the nominal tree. We measure the size of the propagated certificate as a fraction of the nominal tree.

\begin{figure}[H]
    \centering
    \includegraphics[width=0.52\linewidth]{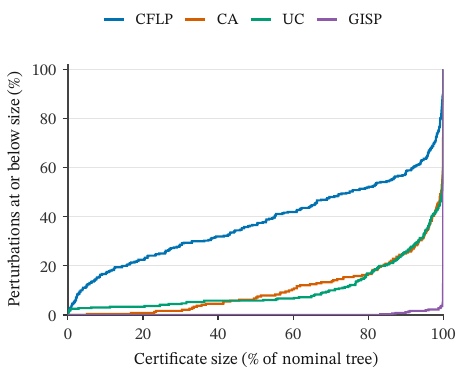}
        
        \caption{Certificate size, that is, the percentage of nodes of the nominal tree retained after propagation. For each value $x$, the curve gives the percentage of single-parameter perturbations whose certificate size is at most $x\%$.}
        \label{fig:whatif.certificate.size.dist}
\end{figure}

\begin{figure}[H]
            \centering
        \includegraphics[width=0.52\linewidth]{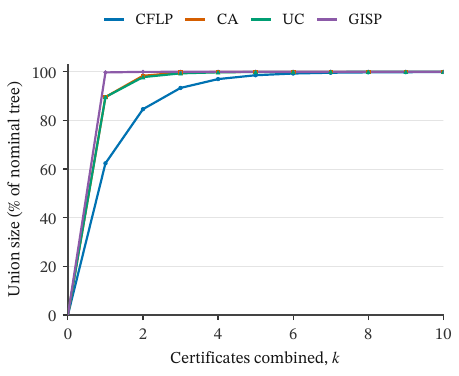}
        \caption{Percentage of nodes of the nominal tree contained in the union of the certificates of $k$ randomly chosen perturbations of the same nominal instance, averaged over instances and random choices.}
        \label{fig:whatif.tree.union}
\end{figure}

\begin{figure}[H]
    \centering
    \includegraphics[width=\linewidth]{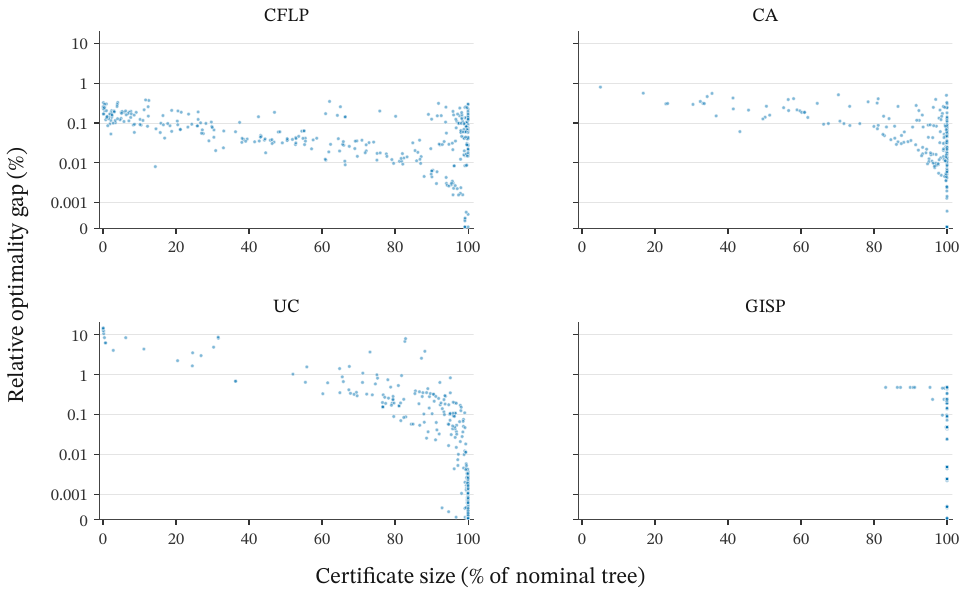}
\caption{Relative optimality gap against certificate size. Each point is one perturbation.}
\label{fig:whatif.dual.gap.vs.certificate.size}
\end{figure}

Figure~\ref{fig:whatif.certificate.size.dist} shows that individual perturbations can require only a small part of the nominal tree: for CFLP, about $20\%$ of the perturbations use at most $15\%$ of it. For CA and UC, however, most perturbations use more than $80\%$ of the tree, and for GISP almost every perturbation uses all of it. Moreover, different perturbations use different parts of the tree: since the certificate depends on $\pert$, which is only known online, retaining the full tree preserves the ability to recover the strongest certificate for any future perturbation. Figure~\ref{fig:whatif.tree.union} shows that little could be saved by trimming it, as the certificates of two to five perturbations together already cover more than $97\%$ of the nominal tree. Finally, Figure~\ref{fig:whatif.dual.gap.vs.certificate.size} shows that small certificates are associated with large gaps, while the smallest gaps are attained with nearly the full tree. A small certificate thus indicates that the stored multipliers deep in the tree yield weak bounds under the perturbation, rather than a saving.


\subsection{Variable fixing}\label{subsec:whatif.var.fixing}

Modern MILP solvers employ sophisticated primal heuristics that often find high-quality feasible solutions quickly~\cite{berthold2025primal}. Moreover, in settings involving repeated optimization of related instances, specialized methods have been shown to obtain good primal solutions at low computational cost~\cite{Shetty1991Jan,Xavier2020Oct,Bertsimas2022Jul,Johnson2020Mar}. We therefore assume that a strong primal bound is available and, to isolate the effectiveness of the proposed variable fixing procedures from the task of obtaining it, use the optimal value of the perturbed problem as the primal cutoff, that is, $\bar z_\Delta = z(\Delta)$. The optimal values are computed separately, and the time required to compute them is excluded from the reported times. Consequently, the results below reflect the fixings attainable with the best possible cutoff; with a weaker bound,  possibly fewer variables would be fixed.

We compare several progressively stronger domain reduction methods. The \textit{single-pass} method evaluates the restricted bounds $f(\pert;x_i=s)$ for all binary variables and $s\in\{0,1\}$ and fixes variables appropriately.
The \textit{binary iteration} method feeds these certified fixings into the node bounds, discards incompatible nodes, and recomputes bounds to identify further fixings.
This procedure is repeated until no new binary variables are fixed, as described in Section~\ref{subsubsec:int.var.fix}.

{The remaining methods introduce additional reductions cumulatively: \textit{+ continuous bounds} tightens continuous variable bounds (Section~\ref{subsubsec:cont.var.tightening}), \textit{+ constraint propagation} further propagates bounds implied by constraints (Section~\ref{subsubsec:constr.bound.prop}), and \textit{+ node implications} also incorporates node-level implied fixings into the constraint propagation step (Section~\ref{subsubsec:constr.bound.prop}). 
For each enhanced method we repeat its enabled reductions and feed their consequences back into
the bounds until no further changes occur.}

\begin{figure}[H]
    \centering
    \includegraphics[width=\linewidth]{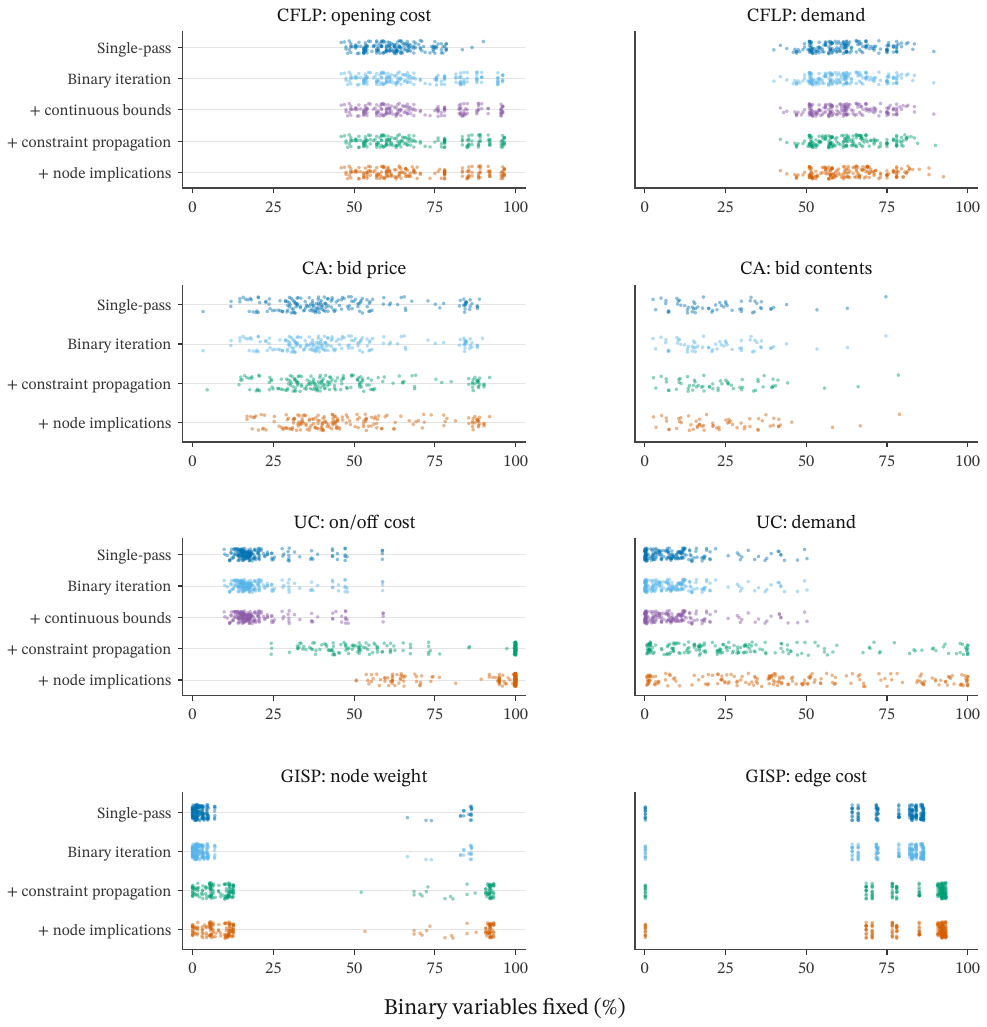}
\caption{Percentage of binary variables fixed for each perturbation, for each method of Table~\ref{tab:whatif.fixing.summary}.}
\label{fig:whatif.dual.fixing.individual}
\end{figure}

Table~\ref{tab:whatif.fixing.summary} reports, for each problem class and perturbation type, the mean percentage of binary variables fixed, the mean GPU time per perturbation, and the mean relative gap after fixing. Figure~\ref{fig:whatif.dual.fixing.individual} shows the percentage fixed for each perturbation. Already the single-pass test of Proposition~\ref{prop:binary.variable.fixing}, which takes between $7$ and $44$ ms per perturbation, fixes a large share of the binaries for CFLP (about $62\%$), GISP edge-cost perturbations ($74\%$), and CA objective perturbations ($44\%$); fewer variables are fixed for CA bit-contents perturbations ($23\%$), UC ($21\%$ and $10\%$), and GISP node-weight perturbations ($11\%$).
\begin{table}[H]
  \centering
  \caption{Variable fixing under single-parameter perturbations, by problem class and perturbation type. Fixed: mean percentage of binary variables fixed. GPU time: mean time per perturbation, excluding preparation of the nominal tree. Gap: mean relative gap after fixing. Methods are cumulative from top to bottom; the row ``+ continuous'' appears only for the problem classes with continuous variables.}
  \label{tab:whatif.fixing.summary}
  \begingroup
\small
\setlength{\tabcolsep}{4pt}
\renewcommand{\arraystretch}{1.10}
\begin{tabular}{@{}lrrrrrr@{}}
\toprule
\multicolumn{7}{@{}l}{\textbf{CFLP}} \\
Method & \multicolumn{3}{c}{{Opening cost}} & \multicolumn{3}{c}{{Demand}} \\
\cmidrule(lr){2-4}\cmidrule(lr){5-7}
 & Fixed (\%) & \shortstack[r]{GPU time\\(ms)} & \shortstack[r]{Relative gap\\(\%)} & Fixed (\%) & \shortstack[r]{GPU time\\(ms)} & \shortstack[r]{Relative gap\\(\%)} \\
\midrule
{Single-pass} & 62.3 & 44 & 0.105 & 61.8 & 15 & 0.0672 \\
{Binary iteration} & 68.7 & 56 & 0.0809 & 61.8 & 29 & 0.0672 \\
{+ continuous bounds} & 68.8 & 93 & 0.0809 & 62.0 & 64 & 0.0671 \\
+ constraint propagation & 69.3 & 434 & 0.0806 & 63.0 & 468 & 0.0548 \\
+ node implications & 69.6 & 954 & 0.0802 & 63.7 & 1187 & 0.0526 \\
\addlinespace[5pt]
\midrule
\multicolumn{7}{@{}l}{\textbf{CA}} \\
Method & \multicolumn{3}{c}{{Bid price}} & \multicolumn{3}{c}{{Bid contents}} \\
\cmidrule(lr){2-4}\cmidrule(lr){5-7}
 & Fixed (\%) & \shortstack[r]{GPU time\\(ms)} & \shortstack[r]{Relative gap\\(\%)} & Fixed (\%) & \shortstack[r]{GPU time\\(ms)} & \shortstack[r]{Relative gap\\(\%)} \\
\midrule
{Single-pass} & 43.9 & 20 & 0.0488 & 23.2 & 7 & 0.214 \\
{Binary iteration} & 44.0 & 25 & 0.0486 & 23.2 & 11 & 0.214 \\
+ constraint propagation & 45.5 & 166 & 0.0486 & 24.4 & 241 & 0.214 \\
+ node implications & 48.5 & 320 & 0.0474 & 25.2 & 440 & 0.214 \\
\addlinespace[5pt]
\midrule
\multicolumn{7}{@{}l}{\textbf{UC}} \\
Method & \multicolumn{3}{c}{{On/off cost}} & \multicolumn{3}{c}{{Demand}} \\
\cmidrule(lr){2-4}\cmidrule(lr){5-7}
 & Fixed (\%) & \shortstack[r]{GPU time\\(ms)} & \shortstack[r]{Relative gap\\(\%)} & Fixed (\%) & \shortstack[r]{GPU time\\(ms)} & \shortstack[r]{Relative gap\\(\%)} \\
\midrule
{Single-pass} & 21.1 & 14 & 0.000882 & 10.1 & 14 & 1.18 \\
{Binary iteration} & 21.2 & 29 & 0.000877 & 10.1 & 29 & 1.18 \\
{+ continuous bounds} & 21.2 & 42 & 0.000877 & 10.1 & 42 & 1.18 \\
+ constraint propagation & 72.7 & 124 & 0.000607 & 32.9 & 881 & 1.11 \\
+ node implications & 89.5 & 55 & 0.000402 & 46.4 & 1011 & 1.11 \\
\addlinespace[5pt]
\midrule
\multicolumn{7}{@{}l}{\textbf{GISP}} \\
Method & \multicolumn{3}{c}{{Node weight}} & \multicolumn{3}{c}{{Edge cost}} \\
\cmidrule(lr){2-4}\cmidrule(lr){5-7}
 & Fixed (\%) & \shortstack[r]{GPU time\\(ms)} & \shortstack[r]{Relative gap\\(\%)} & Fixed (\%) & \shortstack[r]{GPU time\\(ms)} & \shortstack[r]{Relative gap\\(\%)} \\
\midrule
{Single-pass} & 10.6 & 8 & 0.186 & 74.2 & 8 & 0.00138 \\
{Binary iteration} & 10.6 & 15 & 0.186 & 74.2 & 14 & 0.00138 \\
+ constraint propagation & 32.7 & 151 & 0.139 & 80.6 & 61 & 0.00128 \\
+ node implications & 32.7 & 202 & 0.139 & 80.6 & 58 & 0.00128 \\
\bottomrule
\end{tabular}
\endgroup

\end{table}

\subsubsection{Improvements}
\begin{figure}[H]
    \centering
    \includegraphics[width=\linewidth]{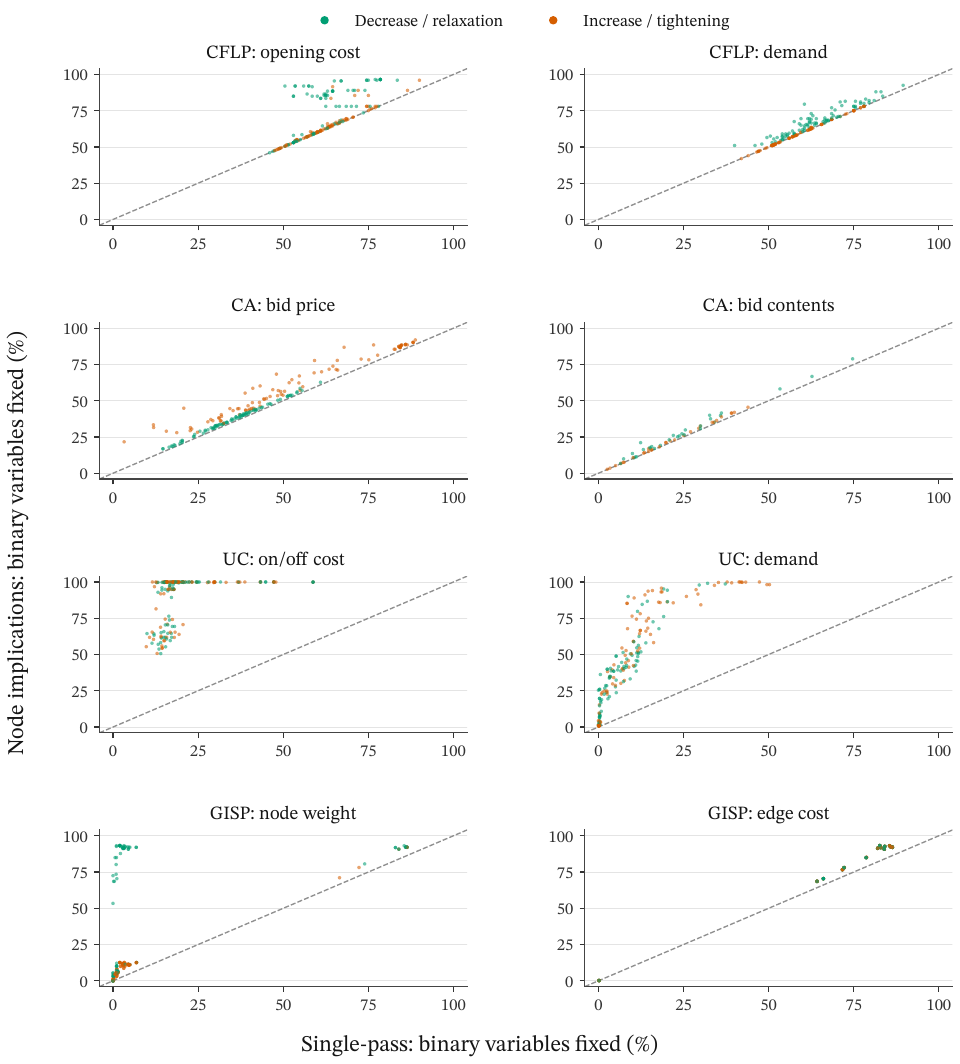}
    \caption{Percentage of binary variables fixed by the single-pass test (horizontal axis) and by the node-implications method (vertical axis), for each perturbation. Points above the diagonal are perturbations for which the improvements fix additional variables.}
    \label{fig:whatif.fixing.single.vs.all}
\end{figure}

\begin{figure}[H]
    \centering
    \includegraphics[width=\linewidth]{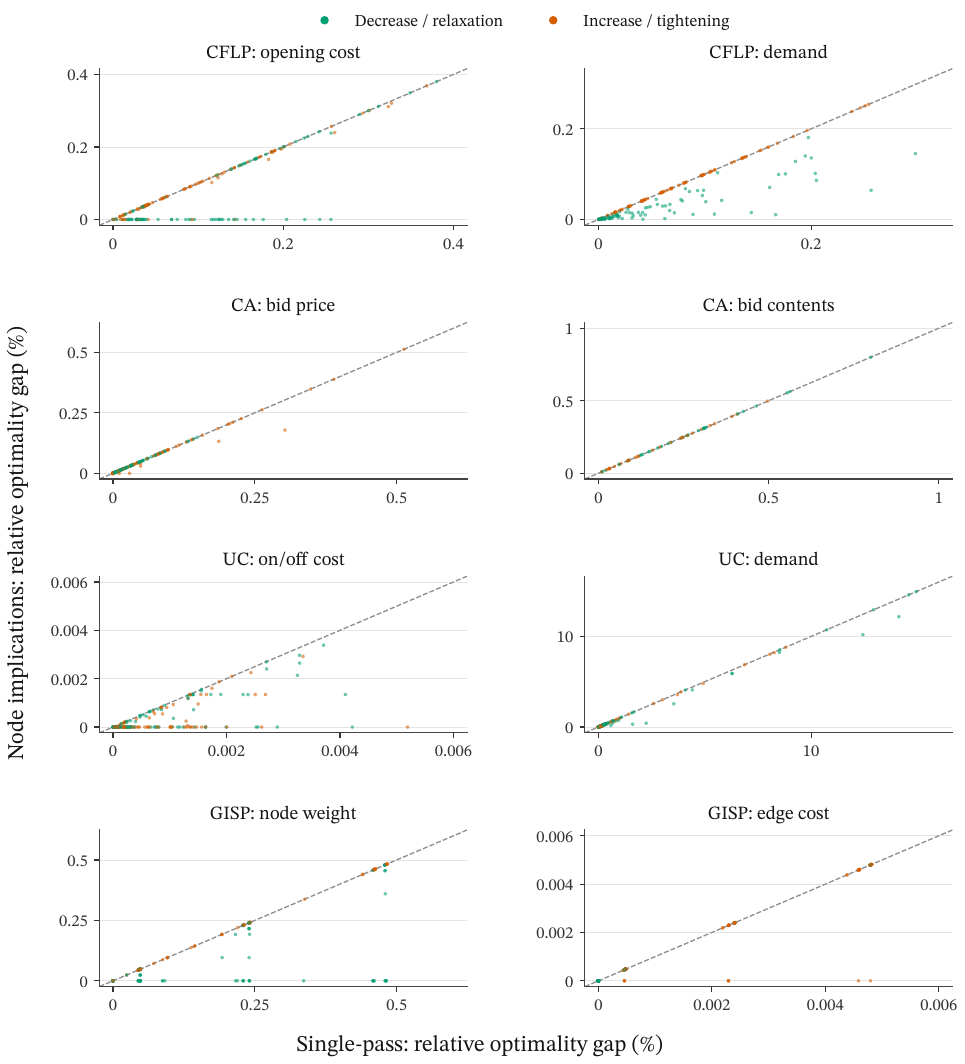}
    \caption{Relative gap of the generalized dual bound after single-pass fixing (horizontal axis) and after the the node-implications method (vertical axis), for each perturbation; colors distinguish decreasing or relaxing perturbations from increasing or tightening ones. Points below the diagonal are perturbations for which the additional fixings strengthen the dual bound, and points on the horizontal axis are those for which the gap is closed.}
    \label{fig:whatif.fixing.dual.reduction}
\end{figure}

The remaining rows of Table~\ref{tab:whatif.fixing.summary} add, one at a time, the iterative fixing of Section~\ref{subsubsec:int.var.fix}, the tightening of continuous bounds of Section~\ref{subsubsec:cont.var.tightening}, constraint propagation as in Section~\ref{subsubsec:constr.bound.prop}, and node implications. Iterative fixing and continuous bounds change little, except for CFLP facility-cost perturbations, where iterative fixing raises the percentage fixed from $62\%$ to $69\%$. Constraint propagation and node implications have the largest effect for UC, where the percentage fixed rises from $21\%$ to $90\%$ for on/off-cost perturbations and from $10\%$ to $46\%$ for demand perturbations, and for GISP node-weight perturbations, where it rises from $11\%$ to $33\%$; these steps increase the GPU time to between $55$ and $1187$ ms per perturbation. Figure~\ref{fig:whatif.fixing.single.vs.all} compares the single-pass and node-implications methods for each perturbation: the gains are concentrated in particular perturbations, such as GISP node-weight decreases, for which the node-implications method fixes up to about $90\%$ of the binaries even when the single-pass test fixes almost none, whereas for most node-weight increases it fixes only about $10\%$. Finally, the fixings also strengthen the dual bound. For CFLP facility-cost, GISP node-weight and UC on/off-cost perturbations, the node-implications method closes the gap entirely for many perturbations, almost all of them decreases (Figure~\ref{fig:whatif.fixing.dual.reduction}). In particular, the mean gap drops from $0.105\%$ to $0.080\%$ for CFLP facility-cost perturbations and from $0.186\%$ to $0.139\%$ for GISP node-weight perturbations.

\subsection{Reoptimization\label{subsec:whatif.reopt}}
We now examine whether the certified fixings of Section~\ref{subsec:whatif.var.fixing} accelerate the solution of the perturbed problems. Since every fixing holds for all optimal solutions, adding the fixings to a perturbed problem does not change its optimal value, so that solving the restricted problem is an exact reoptimization. For each perturbation, we compare the time for Gurobi to solve the perturbed problem with the time to compute the fixings plus the time for Gurobi to solve the perturbed problem with the variables fixed, for single-pass fixing and for the full procedure with node implications.

\begin{figure}[H]
    \centering
    \includegraphics[width=\linewidth]{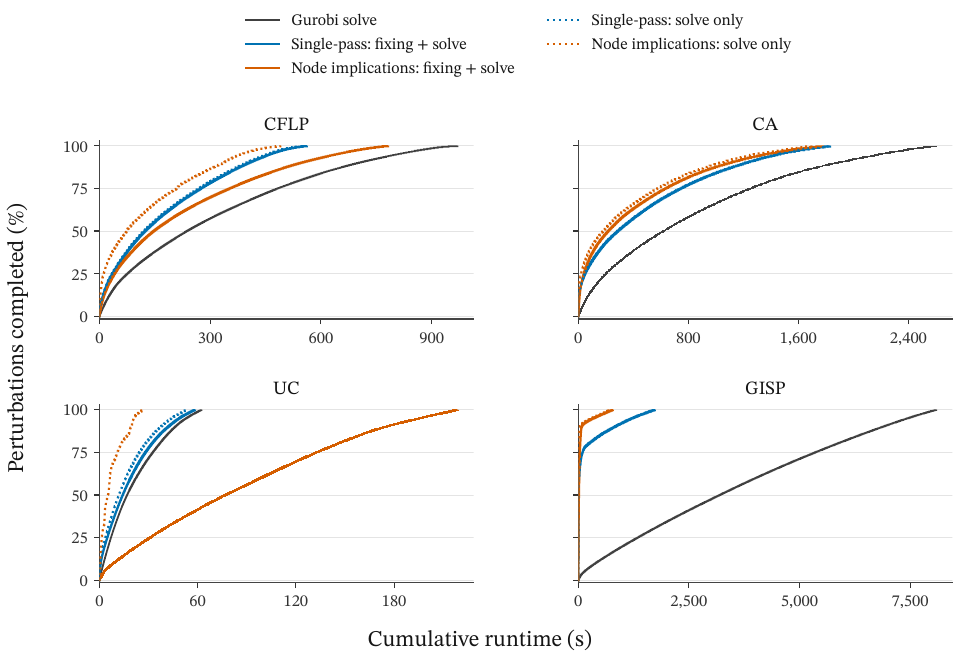}
    \caption{Cumulative runtime over all perturbations of each problem class: solving the perturbed problem with Gurobi, and computing the fixings and then solving the perturbed problem with the fixed variables, for single-pass fixing and for the node-implications method. Dotted lines show the solve time alone, excluding the time to compute the fixings.}
    \label{fig:whatif.runtimes.summary}
\end{figure}

Figure~\ref{fig:whatif.runtimes.summary} reports the cumulative runtime over all perturbations of each problem class. For GISP, fixing reduces the total time by a factor of about $4.6$, from about $8{,}100$~s to about $1{,}750$~s, with single-pass fixing, and by a factor of about $10$, to about $800$~s, with node implications. For CFLP and CA, the reduction is roughly $40\%$ and $30\%$, respectively. UC problems are solved in well under a second, so single-pass fixing gives no gain there. With node implications, the solve becomes more than twice as fast, but the fixing time outweighs this saving. The stronger fixing procedures therefore pay off only when the perturbed problems are hard enough to amortize their cost.

\begin{figure}[H]
    \centering
    \includegraphics[width=0.88\linewidth]{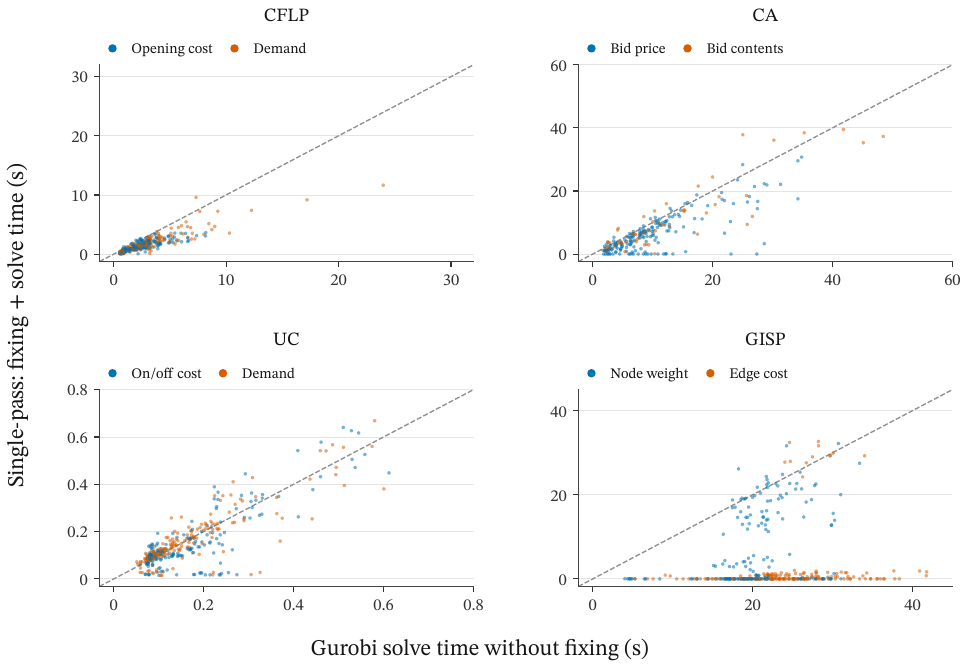}
    \\
\includegraphics[width=0.88\linewidth]{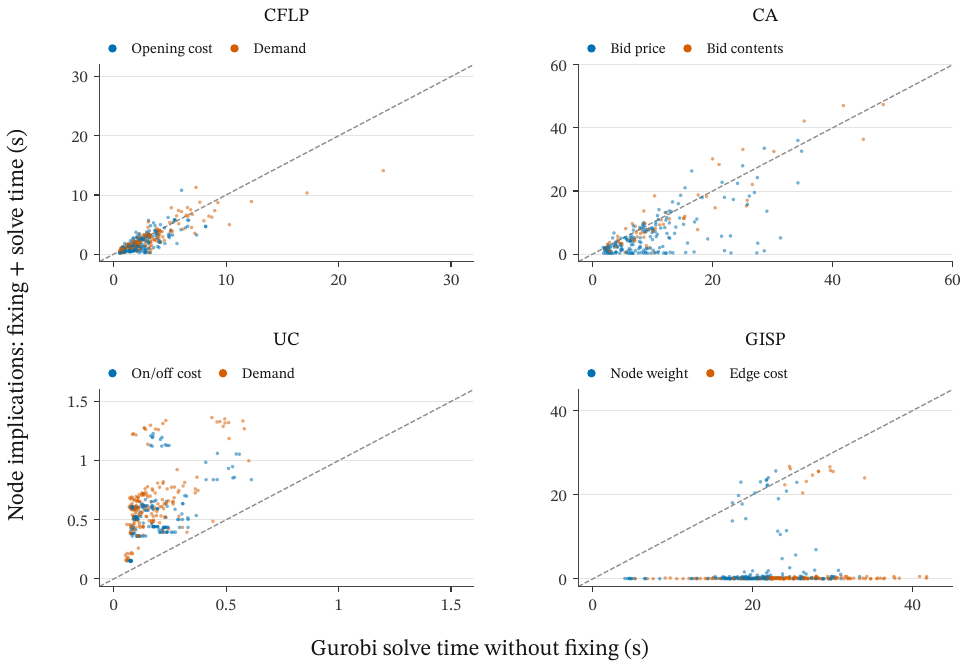}
    \caption{Time for Gurobi to solve each perturbed problem (horizontal axis) against the time to compute the fixings plus the time to solve the perturbed problem with the fixed variables (vertical axis), for single-pass fixing (top) and node implications (bottom). Points below the diagonal are perturbations for which fixing reduces the total time.}
    \label{fig:whatif.runtimes}
\end{figure}

Figure~\ref{fig:whatif.runtimes} compares the two times for each perturbation. A point below the diagonal is a perturbation for which fixing plus solving takes less time than solving the perturbed problem directly. For CFLP, CA and GISP, most points lie below the diagonal. Points near the horizontal axis are perturbations for which the restricted problem is solved almost instantly, because the fixings leave little to decide. This happens for many GISP perturbations, including almost all edge-cost perturbations. For UC with node implications, all points lie above the diagonal, in line with the cumulative results.

\section{Numerical results: Operational reoptimization (multi-parameter perturbations)}\label{sec:dense}

\subsection{Dual bound gaps}
\subsubsection{Generalized dual optimality gaps}

\begin{figure}[H]
    \centering
    \includegraphics[width=\linewidth]{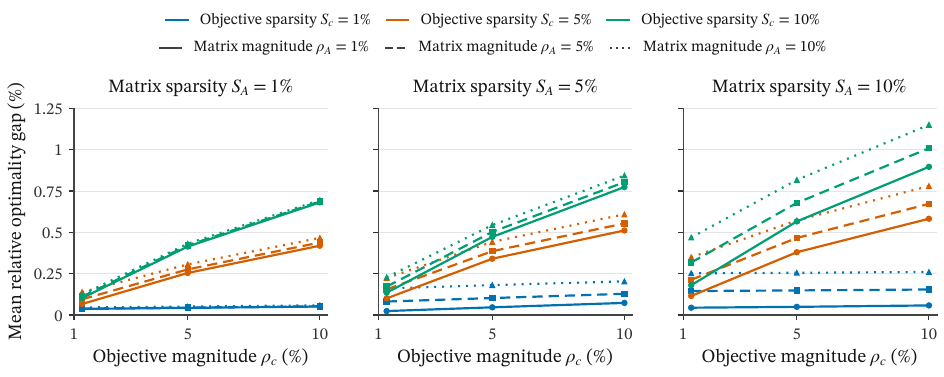}
    \caption{Mean relative optimality gap  for CFLP under simultaneous objective and constraint-matrix perturbation. Plots correspond to fixed matrix sparsity $S_A$, colors to objective sparsity $S_c$, and line styles to matrix perturbation magnitudes $\rho_A$. Each point represents the average across perturbed instances for the corresponding combination of $(S_c,S_A,\rho_c,\rho_A)$.}
\label{fig:dense.cflp.dual.bounds}
\end{figure}

\begin{figure}[H]
    \centering
    \includegraphics[width=\linewidth]{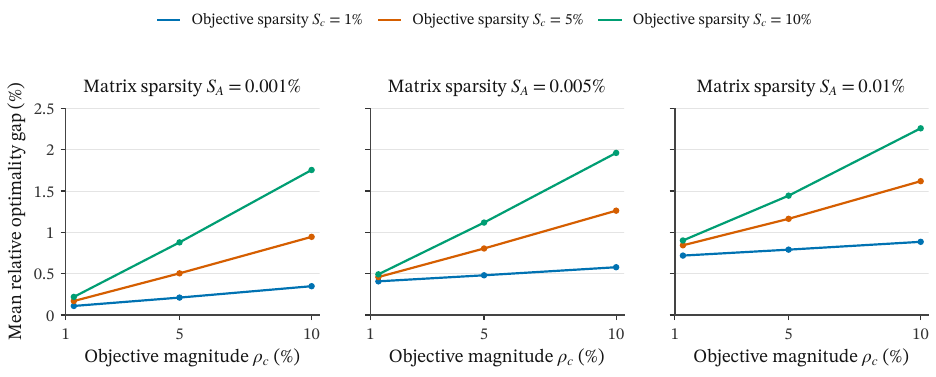}
    \caption{Mean relative optimality gap for CA under simultaneous objective perturbations and matrix binary flips. Plots correspond to fixed matrix sparsity $S_A$ and colors to objective sparsity $S_c$. Each point represents the average across perturbed instances for the corresponding combination of $(S_c,S_A,\rho_c)$.}
\label{fig:dense.ca.dual.bounds}
\end{figure}

\begin{figure}[H]
    \centering
    \includegraphics[width=\linewidth]{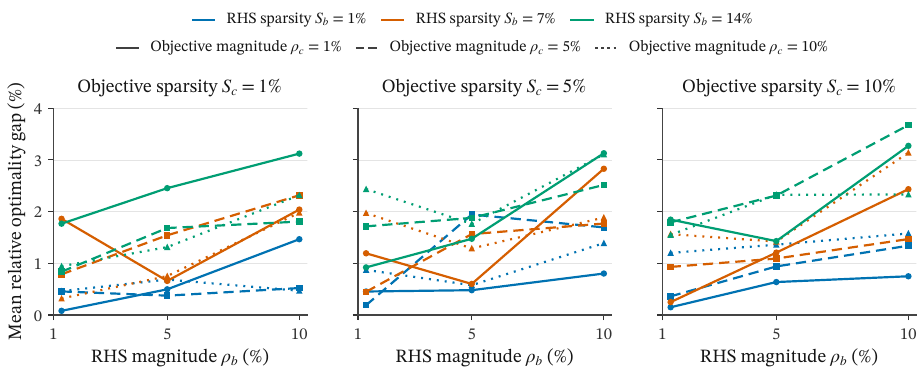}
    \caption{Mean relative optimality gap for UC under simultaneous objective and right-hand side perturbation. Plots correspond to fixed objective sparsity $S_c$, colors to the number of perturbed demand periods $k_b$, and line styles to objective perturbation magnitudes $\rho_c$. Each point represents the average across perturbed instances for the corresponding combination of $(S_c,k_b,\rho_c,\rho_b)$.} 
\label{fig:dense.uc.dual.bounds}
\end{figure}

\begin{figure}[H]
    \centering
    \includegraphics[width=0.5\linewidth]{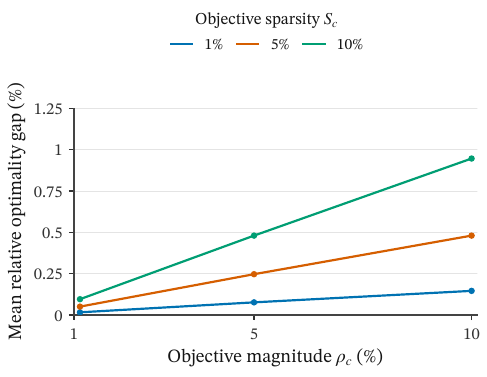}
    \caption{Mean relative optimality gap for GISP under objective perturbation. Each point represents the average across perturbed instances for the corresponding combination of $(S_c,\rho_c)$.}
\label{fig:dense.gisp.dual.bounds}
\end{figure}

We report the relative gap of the generalized dual bound $|z(\pert)-f(\pert)|/|z(\pert)|$, averaged over perturbations of the same sparsity level and magnitude for each dataset in Figures~\ref{fig:dense.cflp.dual.bounds}-\ref{fig:dense.gisp.dual.bounds}. Even when many coefficients change at once, the generalized dual bound remains strong. The mean gap stays below $1.2\%$ for CFLP, $2.3\%$ for CA, $3.7\%$ for UC and $1\%$ for GISP, even in the most extreme setting, in which $10\%$ of the coefficients change by up to $10\%$; for moderate settings, for example with objective sparsity $1\%$, it is typically below $0.3\%$ for CFLP and GISP. In every class, the gap grows roughly linearly with the perturbation magnitude, with a slope that increases with the sparsity, so that the quality of the bound degrades gradually rather than abruptly. The gap is driven mainly by the objective perturbations for CFLP, CA and GISP, whereas for UC it is driven mainly by the right-hand side perturbations, consistent with the single-parameter results of Section~\ref{subsec:whatif.duals}.

\subsubsection{Dual bound speedups}
Table~\ref{tab:whatif.dual.times.dense} reports, for each problem class, the total time our GPU method needs to compute dual bounds for all perturbed problems. It also reports the time Gurobi needs to reach the same bounds. Note that the Gurobi times are in seconds and the GPU times are in milliseconds. The GPU approach is faster in every class, with speedups of about $5\times 10^{4}$ for CA, $1.5\times 10^{4}$ for GISP ($7.52\times 10^{3}$ seconds versus 507 milliseconds), $870\times$ for CFLP ($2.78\times 10^{3}$ seconds versus 3.20 seconds), and $120\times$ for UC (278 seconds versus 2.32 seconds). Per perturbed problem, the GPU method takes between roughly 0.14 and 1.25 milliseconds. Gurobi takes between roughly 0.08 seconds (UC) and 18.6 seconds (GISP). So once the dual information from the base problem is available, dual bounds for all problems with multi-parameter perturbations can be computed almost instantly, whereas Gurobi needs substantially more time to reach the same bounds, especially for CA and GISP.
\begin{table}[H]
  \centering
  \caption{Total runtime for obtaining dual bounds for all multi-parameter perturbations in each problem class.}
  \label{tab:whatif.dual.times.dense}
  \begingroup
  \small
  \setlength{\tabcolsep}{5pt}
  \begin{tabular}{lrrr}
    \toprule
    Problem class & \shortstack[r]{Perturbed problems\\(total)}
      & \shortstack[r]{Gurobi\\(s)}
      & \shortstack[r]{GPU MILP Dual\\(ms)}\\
    \midrule
    CFLP & 3645 & $2.78\times 10^{3}$ & $3.20\times 10^{3}$ \\
    CA   & 1,215 & $9.32\times 10^{3}$ & 173 \\
    UC   & 3645 & 278                  & $2.32\times 10^{3}$ \\
    GISP & 405  & $7.52\times 10^{3}$ & 507 \\
    \bottomrule
  \end{tabular}
  \endgroup
  \par\smallskip
\end{table}

\begin{figure}[H]
    \centering
    \includegraphics[width=0.5\linewidth]{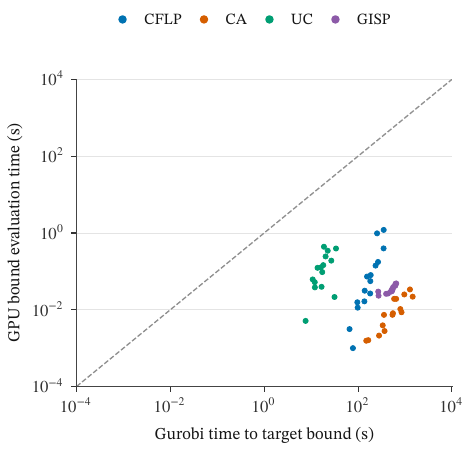}
    \caption{Total time over the multi-parameter perturbations of each nominal instance: time for Gurobi to reach the generalized dual bound $f(\pert)$ (horizontal axis) and time for the GPU to compute it (vertical axis), both in seconds on logarithmic scales. Each point is one nominal instance; the dashed line is the diagonal.}
\label{fig:dense.dual.time.individual}
\end{figure}

\subsubsection{Dual gap decomposition}
\begin{figure}[H]
    \centering
    \includegraphics[width=\linewidth]{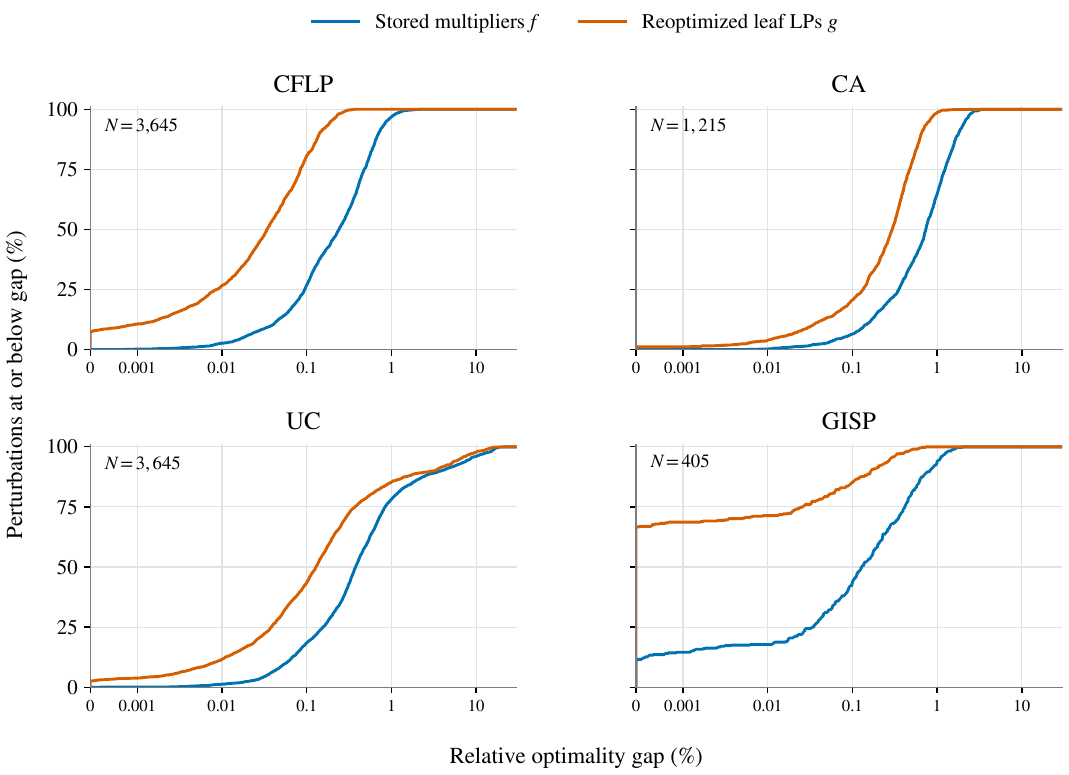}
    \caption{Relative gap of the generalized dual bound $f(\pert)$ and of the bound $g(\pert)$ obtained by reoptimizing the leaf LPs of the nominal tree. For each value of $x$, the curves give the percentage of perturbations whose gap is at most $x\%$, $N$ is the number of perturbations of each problem class.}
\label{fig:dense.reopt.leaves}
\end{figure}

\begin{figure}[H]
    \centering
    \includegraphics[width=\linewidth]{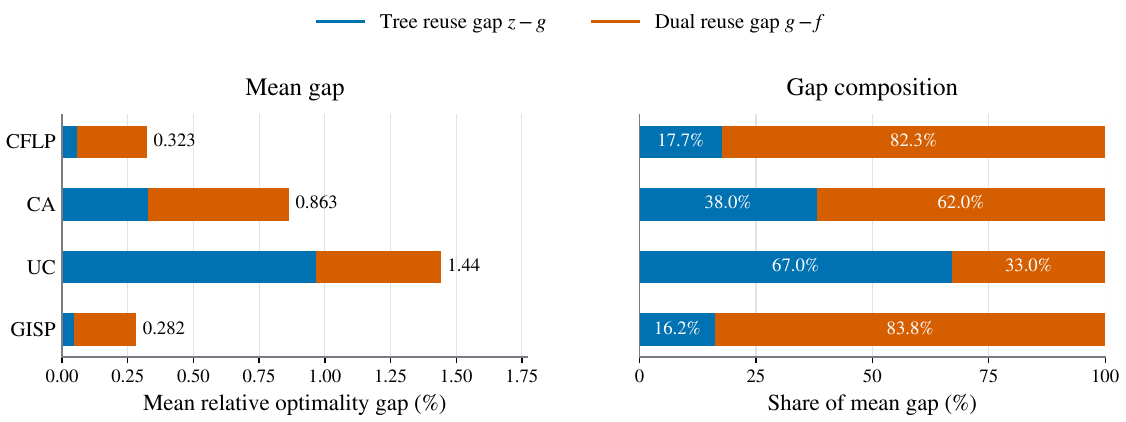}
    \caption{Decomposition of the multi-parameter gap into tree reuse $z(\pert)-g(\pert)$ and dual reuse $g(\pert)-f(\pert)$, as in \eqref{eq:abs.gap.breakdown}. Left: mean relative contribution of each component, with the mean total gap indicated. Right: each component as a share of the mean total relative gap.}
\label{fig:dense.gap.decomp}
\end{figure}

Reoptimizing the leaf LPs raises the share of GISP perturbations with zero gap from about 11\% to about
67\% and substantially reduces CFLP gaps (Figure~\ref{fig:dense.reopt.leaves}). For UC the gain is smaller, so most of the gap again comes from reusing the nominal tree.
{Figure~\ref{fig:dense.gap.decomp} shows that dual reuse accounts for 83.8\% of the mean gap for GISP, 82.3\% for CFLP, and 62.0\% for CA. UC is the exception: tree reuse accounts for 67.0\% of its mean gap. Thus, unlike in the what-if setting, the dominant source of the CA mean gap under
multi-parameter perturbations is dual reuse.}

\subsubsection{Bound propagation tree utilization}
\begin{figure}[H]
    \centering
    \includegraphics[width=0.5\linewidth]{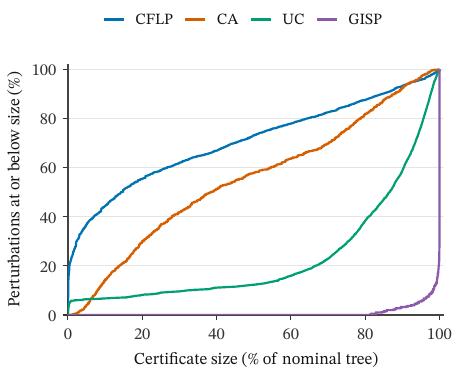}
    \caption{Certificate size, that is, the percentage of nodes of the nominal tree retained after propagation. For each value $x$, the curve gives the percentage of distributed-change perturbations whose certificate size is at most x\%.}
\label{fig:dense.distribution.tree.certificate.sizes}
\end{figure}

\begin{figure}[H]
    \centering
    \includegraphics[width=\linewidth]{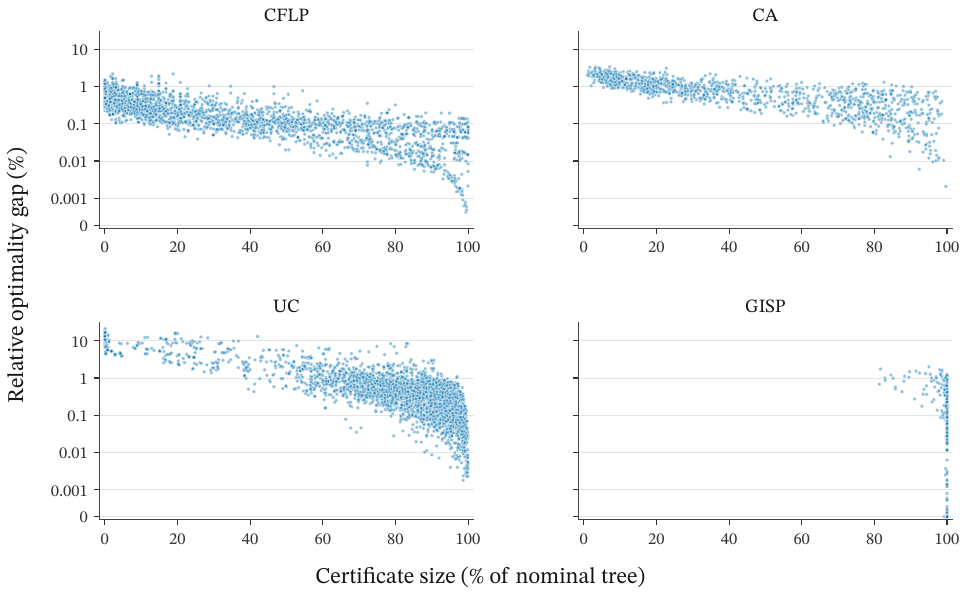}
    \caption{Relative gap $\left(z(\pert)-f(\pert)\right)/|z(\pert)|$ against certificate size. Each point is one perturbation.}
\label{fig:dense.gap.vs.cert.size}
\end{figure}

\begin{figure}[H]
    \centering
    \includegraphics[width=0.5\linewidth]{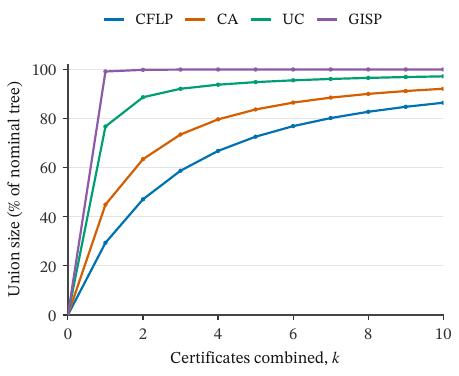}
    \caption{Percentage of nodes of the nominal tree contained in the union of the certificates of $k$ randomly chosen multi-parameter perturbations of the same nominal instance, averaged over instances and random choices.}
\label{fig:dense.distribution.tree.union}
\end{figure}

Compared with the what-if setting, multi-parameter perturbations produce smaller certificates for CFLP and CA, while UC and GISP still use most of the tree (Figure~\ref{fig:dense.distribution.tree.certificate.sizes}). As before, smaller certificates are associated with larger gaps, so they reflect weaker deep-node bounds, not savings (Figure~\ref{fig:dense.gap.vs.cert.size}). The union of certificates grows more slowly than in the what-if setting, but ten certificates still cover 86\% to 100\% of the tree (Figure~\ref{fig:dense.distribution.tree.union}). Hence the full nominal tree must still be stored.

\subsection{Variable fixing}
\begin{figure}[htbp]
\centering
\includegraphics[width=\linewidth]{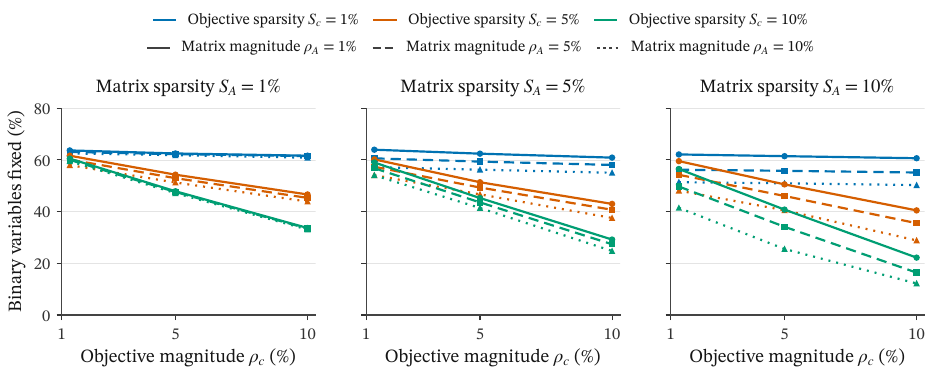}
\caption{Mean percentage of binary variables fixed by the single-pass method for CFLP under {multi-parameter perturbations}.}
\label{fig:dense.varfixing}
\end{figure}

{As in Section~\ref{subsec:whatif.var.fixing}, we use the optimal value of the perturbed problem as the primal cutoff, $\bar{z}_{\pert} = z(\pert)$. We present results for the \textit{single-pass} variable fixing method, showing the effect of sparsity and magnitude of multi-parameter perturbations on variable fixing.}
When all sparsities and magnitudes are 1\%, single-pass fixing fixes about 64\% of the binaries for CFLP, 33\% for CA, and 68\% for GISP, which is comparable to the what-if setting, but only 16\% for UC, compared with 46\% to 90\% there (Figure~\ref{fig:dense.varfixing}). The fraction drops quickly as the sparsity or magnitude grows, falling to about 12\% for CFLP and to nearly zero for CA, UC, and GISP.

\begin{figure}[htbp]
\centering
\includegraphics[width=\linewidth]{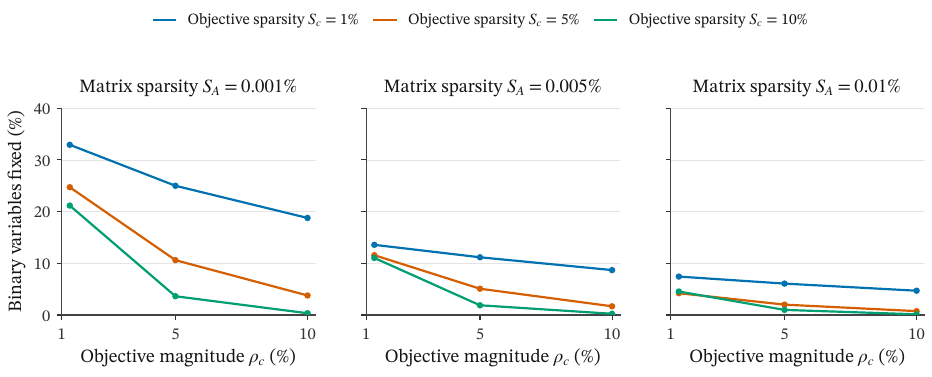}
\caption{Mean percentage of binary variables fixed by the single-pass method for CA under {multi-parameter perturbations}.}
\label{fig:dense.varfixing.ca}
\end{figure}

\begin{figure}[htbp]
\centering
\includegraphics[width=\linewidth]{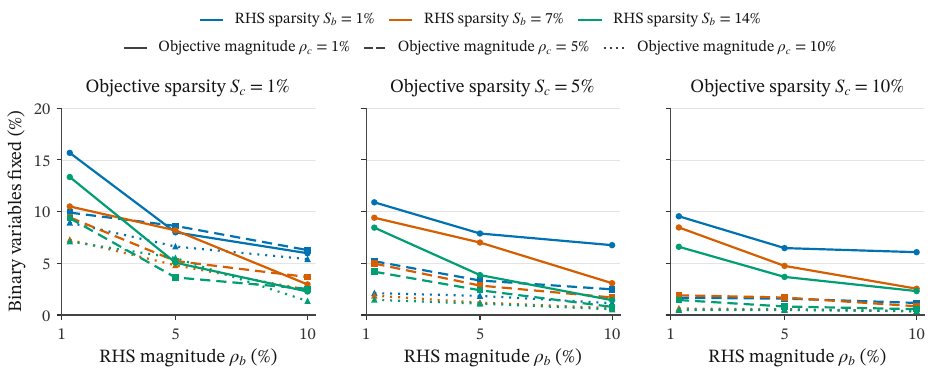}
\caption{Mean percentage of binary variables fixed by the single-pass method for UC under {multi-parameter perturbations}.}
\label{fig:dense.varfixing.uc}
\end{figure}

\begin{figure}[htbp]
\centering
\includegraphics[width=0.52\linewidth]{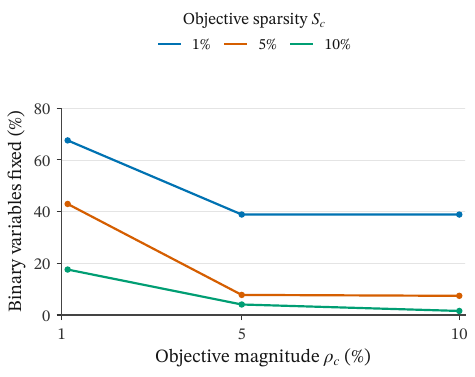}
\caption{Mean percentage of binary variables fixed by the single-pass method for GISP under {multi-parameter perturbations}.}
\label{fig:dense.varfixig.gisp}
\end{figure}

\subsection{Reoptimization}

\begin{figure}[htbp]
    \centering
    \includegraphics[width=\linewidth]{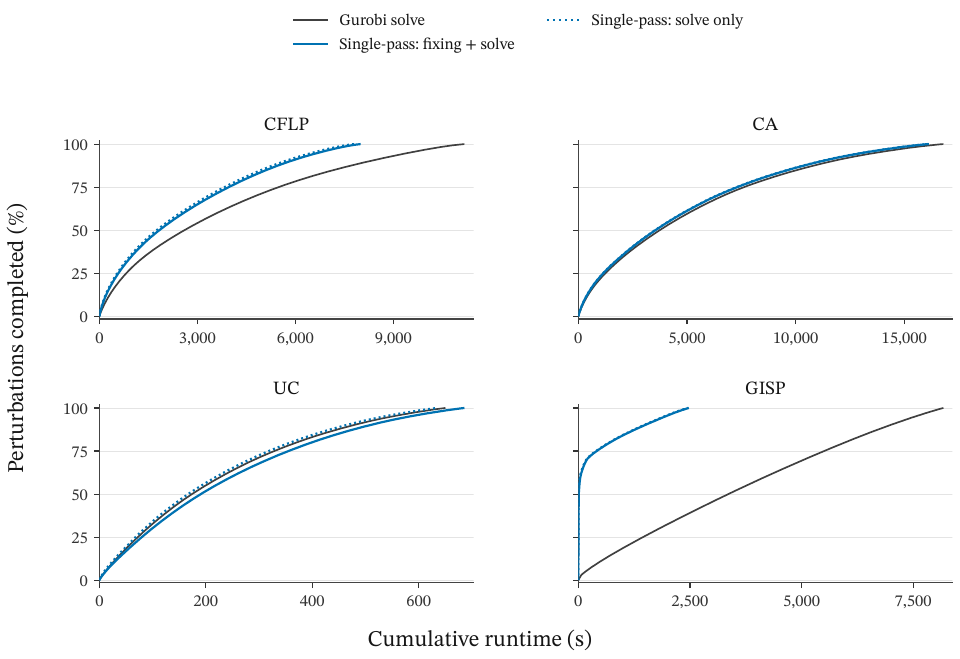}
    \caption{Cumulative runtime over all perturbations of each problem class: solving the perturbed problem with Gurobi, and computing the fixings and then solving the perturbed problem with the fixed variables for single-pass.}
    \label{fig:dense.reopt}
\end{figure}

Fixing variables reduces the total runtime for GISP by a factor of about 3, compared with about 4.6 in the what-if setting (Figure~\ref{fig:dense.reopt}). For the other classes the benefit is smaller than in the what-if setting: the total runtime drops by about 30\% for CFLP, changes little for CA, and rises slightly for UC, which is consistent with the fewer fixings obtained under multi-parameter perturbations.
\section{Conclusion}
\label{sec:conclusion}

In this paper, we developed a framework for obtaining valid dual bounds for mixed-integer linear programs under general perturbations to the problem data by reusing the branch-and-bound tree and LP dual information from a nominal solve. A key feature of the proposed approach is that it is inherently GPU-native: the online computations consist primarily of parallel node-level evaluations and simple tree-based operations, making them particularly well suited for GPU acceleration. Our computational results demonstrate that the method produces extremely tight dual bounds across all four problem classes at a very small computational cost, typically requiring less than a millisecond per perturbation. The results are particularly strong for ``what-if'' analysis involving changes to individual coefficients, where the duality gaps are often a small fraction of a percent. Importantly, the method also performs well under more general perturbations involving simultaneous changes to many coefficients, with the quality of the bounds degrading gradually as the magnitude and sparsity of the perturbations increase. Beyond dual bounds, the framework can certify optimal variable fixings, often fixing a substantial fraction of the binary variables at a relatively small additional cost. However, the benefits of these fixings for subsequent reoptimization with Gurobi are problem-dependent. At one extreme, they lead to nearly ten times reduction in total solution time for GISP; at the other, they provide no net computational improvement for UC. These findings suggest that the most consistent strength of the proposed framework is its ability to rapidly produce high-quality dual bounds, while its benefits for accelerating exact reoptimization depend on the problem structure.

The main limitation of the current approach is its reliance on pure branch-and-bound trees, without cutting planes or presolve. Consequently, the size of the nominal tree that can be generated, stored, and processed restricts the scale of problems to which the method can presently be applied. An important direction for future research is to extend the framework to branch-and-cut trees, building for example on recent work~\cite{becu2024approximating} on incorporating Gomory mixed-integer cuts into branch-and-bound. We hope that incorporating such cuts will substantially reduce the size of the nominal trees, thereby extending the approach to larger and more realistic instances. A second important direction is to develop sensitivity-aware presolve techniques. Standard presolve operations may depend on the nominal problem data and therefore need not remain valid after perturbation. Developing presolve reductions that remain valid over a specified family of perturbations, or that can be efficiently updated when the data change, could substantially improve the scalability of the framework. A third direction is to use multiple trees, and trees larger than the one needed to solve the nominal problem. Our experimental results indicate that such trees can yield stronger dual bounds for the perturbed problems. This must be balanced against the added cost of generating, storing, and processing these trees, so a key question is how to select a small set of trees that gives strong bounds across the family of perturbations. Together, these extensions could bring the proposed GPU-native approach closer to practical large-scale MILP sensitivity analysis and reoptimization.

\section*{Acknowledgements}
Santanu S. Dey gratefully acknowledges support from the Air Force Office of Scientific Research (AFOSR) under grant FA9550261B217.
\bibliography{bib}

@book{berthold2025primal,
  author    = {Berthold, Timo and Lodi, Andrea and Salvagnin, Domenico},
  title     = {Primal Heuristics in Integer Programming},
  publisher = {Cambridge University Press},
  address   = {Cambridge},
  year      = {2025},
  isbn      = {978-1-009-57478-5},
  doi       = {10.1017/9781009574792}
}

@article{Wolsey1981,
  author  = {Wolsey, Laurence A.},
  title   = {Integer programming duality: Price functions and sensitivity analysis},
  journal = {Mathematical Programming},
  volume  = {20},
  number  = {1},
  pages   = {173--195},
  year    = {1981},
  doi     = {10.1007/BF01589344}
}

@article{GuzelsoyRalphs2007,
  author  = {G{\"u}zelsoy, Menal and Ralphs, Ted K.},
  title   = {Duality for mixed-integer linear programs},
  journal = {International Journal of Operations Research},
  volume  = {4},
  number  = {3},
  pages   = {118--137},
  year    = {2007}
}

@incollection{RalphsGuzelsoy2005,
  author    = {Ralphs, Ted K. and G{\"u}zelsoy, Menal},
  title     = {The {SYMPHONY} callable library for mixed integer programming},
  booktitle = {The Next Wave in Computing, Optimization, and Decision Technologies},
  editor    = {Golden, Bruce and Raghavan, S. and Wasil, Edward},
  series    = {Operations Research/Computer Science Interfaces Series},
  volume    = {29},
  pages     = {61--76},
  publisher = {Springer},
  address   = {Boston, MA},
  year      = {2005},
  doi       = {10.1007/0-387-23529-9_5}
}

@inproceedings{GamrathHillerWitzig2015,
  author    = {Gamrath, Gerald and Hiller, Benjamin and Witzig, Jakob},
  title     = {Reoptimization techniques for {MIP} solvers},
  booktitle = {Experimental Algorithms (SEA 2015)},
  series    = {Lecture Notes in Computer Science},
  volume    = {9125},
  pages     = {181--192},
  publisher = {Springer},
  year      = {2015},
  doi       = {10.1007/978-3-319-20086-6_14}
}

@article{becu2024approximating,
  author        = {Becu, Berkay and Dey, Santanu S. and Qiu, Feng and Xavier, {\'A}linson S.},
  title         = {Approximating the {G}omory Mixed-Integer Cut Closure Using Historical Data},
  journal={arXiv preprint arXiv:2411.15090},
  year          = {2024},
  eprint        = {2411.15090},
  archivePrefix = {arXiv},
  primaryClass  = {math.OC},
  url           = {https://arxiv.org/abs/2411.15090}
}

@article{li2023large,
  title   = {Large Language Models for Supply Chain Optimization},
  author  = {Li, Beibin and Mellou, Konstantina and Zhang, Bo and Pathuri, Jeevan and Menache, Ishai},
  journal = {arXiv preprint arXiv:2307.03875},
  year    = {2023}
}

@article{robinson1977characterization,
  author  = {Robinson, Stephen M.},
  title   = {A Characterization of Stability in Linear Programming},
  journal = {Operations Research},
  volume  = {25},
  number  = {3},
  pages   = {435--447},
  year    = {1977},
  doi     = {10.1287/opre.25.3.435}
}

@article{robinson1975stability,
  author  = {Robinson, Stephen M.},
  title   = {Stability Theory for Systems of Inequalities. {Part I}: Linear Systems},
  journal = {SIAM Journal on Numerical Analysis},
  volume  = {12},
  number  = {5},
  pages   = {754--769},
  year    = {1975},
  doi     = {10.1137/0712056}
}

@article{freund1985postoptimal,
  author  = {Freund, Robert M.},
  title   = {Postoptimal Analysis of a Linear Program under Simultaneous Changes in Matrix Coefficients},
  journal = {Mathematical Programming Study},
  volume  = {24},
  pages   = {1--13},
  year    = {1985}
}

@incollection{mills1956marginal,
  author    = {Mills, Harlan D.},
  title     = {Marginal Values of Matrix Games and Linear Programs},
  booktitle = {Linear Inequalities and Related Systems},
  editor    = {Kuhn, Harold W. and Tucker, Albert W.},
  series    = {Annals of Mathematics Studies},
  volume    = {38},
  pages     = {183--193},
  publisher = {Princeton University Press},
  year      = {1956}
}

@article{derval2024sensitivity,
  author       = {Guillaume Derval and Damien Ernst and
                  Quentin Louveaux and Bardhyl Miftari},
  title        = {Sensitivity Analysis for Linear Changes of
                  the Constraint Matrix of a (Mixed-Integer)
                  Linear Program},
    journal={arXiv preprint arXiv:2410.14443},
  year         = {2026},
  eprint       = {2410.14443},
  archivePrefix = {arXiv},
  primaryClass = {math.OC},
  doi          = {10.48550/arXiv.2410.14443},
  url          = {https://arxiv.org/abs/2410.14443}
}

@article{CifuentesDeyXu2026,
  author  = {Cifuentes, Diego and Dey, Santanu S. and Xu, Jingye},
  title   = {Sensitivity Analysis for Mixed Binary Quadratic Programming},
  journal = {Mathematical Programming},
  volume  = {216},
  pages   = {379--424},
  year    = {2026},
  doi     = {10.1007/s10107-025-02265-6}
}

@article{istace2026analyzing,
  title={Analyzing changes in optimal variables in linear programming with uncertain parameters},
  author={Istace, Baptiste and Derval, Guillaume and Miftari, Bardhyl and Louveaux, Quentin},
  journal={arXiv preprint arXiv:2606.18830},
  year={2026}
}

@article{boveroux2026complexity,
  title={On the Complexity of Linear Programs with Parametric Constraint Matrices},
  journal={arXiv preprint arXiv:2610.20422},
  author={Boveroux, Laurie and Louveaux, Quentin and Lodi, Andrea and Carvalho, Margarida},
  year={2026}
}

@article{derval2025efficient,
  title={Efficient LP warmstarting for linear modifications of the constraint matrix},
  author={Derval, Guillaume and Miftari, Bardhyl and Ernst, Damien and Louveaux, Quentin},
  journal={arXiv preprint arXiv:2501.04151},
  year={2025}
}

@article{cook1986sensitivity,
  title={Sensitivity theorems in integer linear programming},
  author={Cook, William and Gerards, Albertus MH and Schrijver, Alexander and Tardos, {\'E}va},
  journal={Mathematical Programming},
  volume={34},
  number={3},
  pages={251--264},
  year={1986},
  publisher={Springer}
}

@article{kilincc2009approximating,
  title={Approximating the stability region for binary mixed-integer programs},
  author={K{\i}l{\i}n{\c{c}}-Karzan, Fatma and Toriello, Alejandro and Ahmed, Shabbir and Nemhauser, George and Savelsbergh, Martin},
  journal={Operations Research Letters},
  volume={37},
  number={4},
  pages={250--254},
  year={2009},
  publisher={Elsevier}
}

@article{schrage1985sensitivity,
author = {Schrage, Linus and Wolsey, Laurence},
title = {Sensitivity Analysis for Branch and Bound Integer Programming},
journal = {Operations Research},
volume = {33},
number = {5},
pages = {1008-1023},
year = {1985},
doi = {10.1287/opre.33.5.1008}
}

@article{Dawande2000Aug,
	author = {Dawande, M. W. and Hooker, J. N.},
	title = {{Inference-Based Sensitivity Analysis for Mixed Integer/Linear Programming}},
	journal = {Oper. Res.},
	year = {2000},
	month = aug,
	publisher = {INFORMS},
	url = {https://pubsonline.informs.org/doi/epdf/10.1287/opre.48.4.623.12420}
}

@misc{huang2024distributional,
      title={Distributional MIPLIB: a Multi-Domain Library for Advancing ML-Guided MILP Methods}, 
      author={Weimin Huang and Taoan Huang and Aaron M Ferber and Bistra Dilkina},
      year={2024},
      eprint={2406.06954},
      archivePrefix={arXiv},
      primaryClass={cs.LG}
}

@inproceedings{leyton2000towards,
  title={Towards a universal test suite for combinatorial auction algorithms},
  author={Leyton-Brown, Kevin and Pearson, Mark and Shoham, Yoav},
  booktitle={Proceedings of the 2nd ACM conference on Electronic commerce},
  pages={66--76},
  year={2000}
}

@inproceedings{ferber2022backdoors,
  title={Learning Pseudo-Backdoors for Mixed Integer Programs},
  author={Ferber, Aaron and Song, Jialin and Dilkina, Bistra and Yue, Yisong},
  booktitle={Integration of Constraint Programming, Artificial Intelligence, and Operations Research},
  pages={91--102},
  year={2022},
  organization={Springer}
}

@article{hochbaum1997gisp,
  title={Forest harvesting and minimum cuts: a new approach to handling spatial constraints},
  author={Hochbaum, Dorit S and Pathria, Anu},
  journal={Forest Science},
  volume={43},
  number={4},
  pages={544--554},
  year={1997},
  publisher={Oxford University Press}
}

@article{cornuejols1991comparison,
  title={A comparison of heuristics and relaxations for the capacitated plant location problem},
  author={Cornu{\'e}jols, G{\'e}rard and Sridharan, Ranjani and Thizy, Jean-Michel},
  journal={European journal of operational research},
  volume={50},
  number={3},
  pages={280--297},
  year={1991},
  publisher={Elsevier}
}

@article{Chen2016Feb,
	author = {Chen, Yonghong and Casto, Aaron and Wang, Fengyu and Wang, Qianfan and Wang, Xing and Wan, Jie},
	title = {{Improving Large Scale Day-Ahead Security Constrained Unit Commitment Performance}},
	journal = {IEEE Trans. Power Syst.},
	volume = {31},
	number = {6},
	pages = {4732--4743},
	year = {2016},
	month = feb,
	publisher = {IEEE},
	doi = {10.1109/TPWRS.2016.2530811}
}

@article{Shetty1991Jan,
	author = {Shetty, Bala},
	title = {{Reoptimization procedures for networks with side constraints}},
	journal = {Comput. Ind. Eng.},
	volume = {20},
	number = {2},
	pages = {243--249},
	year = {1991},
	month = jan,
	issn = {0360-8352},
	publisher = {Pergamon},
	doi = {10.1016/0360-8352(91)90028-5}
}

@article{Xavier2020Oct,
	author = {Xavier, {\ifmmode\acute{A}\else\'{A}\fi}linson S. and Qiu, Feng and Ahmed, Shabbir},
	title = {{Learning to Solve Large-Scale Security-Constrained Unit Commitment Problems}},
	journal = {INFORMS J. Comput.},
	year = {2020},
	month = oct,
	publisher = {INFORMS},
	url = {https://pubsonline.informs.org/doi/10.1287/ijoc.2020.0976}
}

@article{Bertsimas2022Jul,
	author = {Bertsimas, Dimitris and Stellato, Bartolomeo},
	title = {{Online Mixed-Integer Optimization in Milliseconds}},
	journal = {INFORMS J. Comput.},
	volume = {34},
	number = {4},
	pages = {2229--2248},
	year = {2022},
	month = jul,
	issn = {1526-5528},
	publisher = {INFORMS},
	doi = {10.1287/ijoc.2022.1181}
}

@article{Johnson2020Mar,
	author = {Johnson, Emma S. and Ahmed, Shabbir and Dey, Santanu S. and Watson, Jean-Paul},
	title = {{A K-Nearest Neighbor Heuristic for Real-Time DC Optimal Transmission Switching}},
	journal = {arXiv},
	year = {2020},
	month = mar,
	eprint = {2003.10565},
	doi = {10.48550/arXiv.2003.10565}
}

@article{Shah2025Jul,
	author = {Shah, Prachi and Dey, Santanu S.},
	title = {{Improving Full Strong Branching Decisions by Incorporating Additional Information}},
	journal = {arXiv},
	year = {2025},
	month = jul,
	eprint = {2507.09455},
	doi = {10.48550/arXiv.2507.09455}
}

@article{Dey2024May,
	author = {Dey, Santanu S. and Dubey, Yatharth and Molinaro, Marco and Shah, Prachi},
	title = {{A theoretical and computational analysis of full strong-branching}},
	journal = {Math. Program.},
	volume = {205},
	number = {1},
	pages = {303--336},
	year = {2024},
	month = may,
	issn = {1436-4646},
	publisher = {Springer Berlin Heidelberg},
	doi = {10.1007/s10107-023-01977-x}
}

@article{SantosXavier2025Jun,
	author = {Santos Xavier, Alinson and Qiu, Feng and Gu, Xiaoyi and Becu, Berkay and Dey, Santanu S.},
	title = {{MIPLearn: An Extensible Framework for Learning-Enhanced Optimization (Version 0.4)}},
	journal = {Zenodo},
	year = {2024},
	doi = {10.5281/zenodo.15634125}
}

@article{savelsbergh1994preprocessing,
  title={Preprocessing and probing techniques for mixed integer programming problems},
  journal={ORSA Journal on Computing},
  author={Savelsbergh, Martin W. P.},
  volume={6},
  number={4},
  pages={445--454},
  year={1994},
  publisher={INFORMS}
}
\bibliographystyle{abbrv}
\newpage
\appendix

\section{Alternative LP dual}
\label{appendix:alternative.dual}
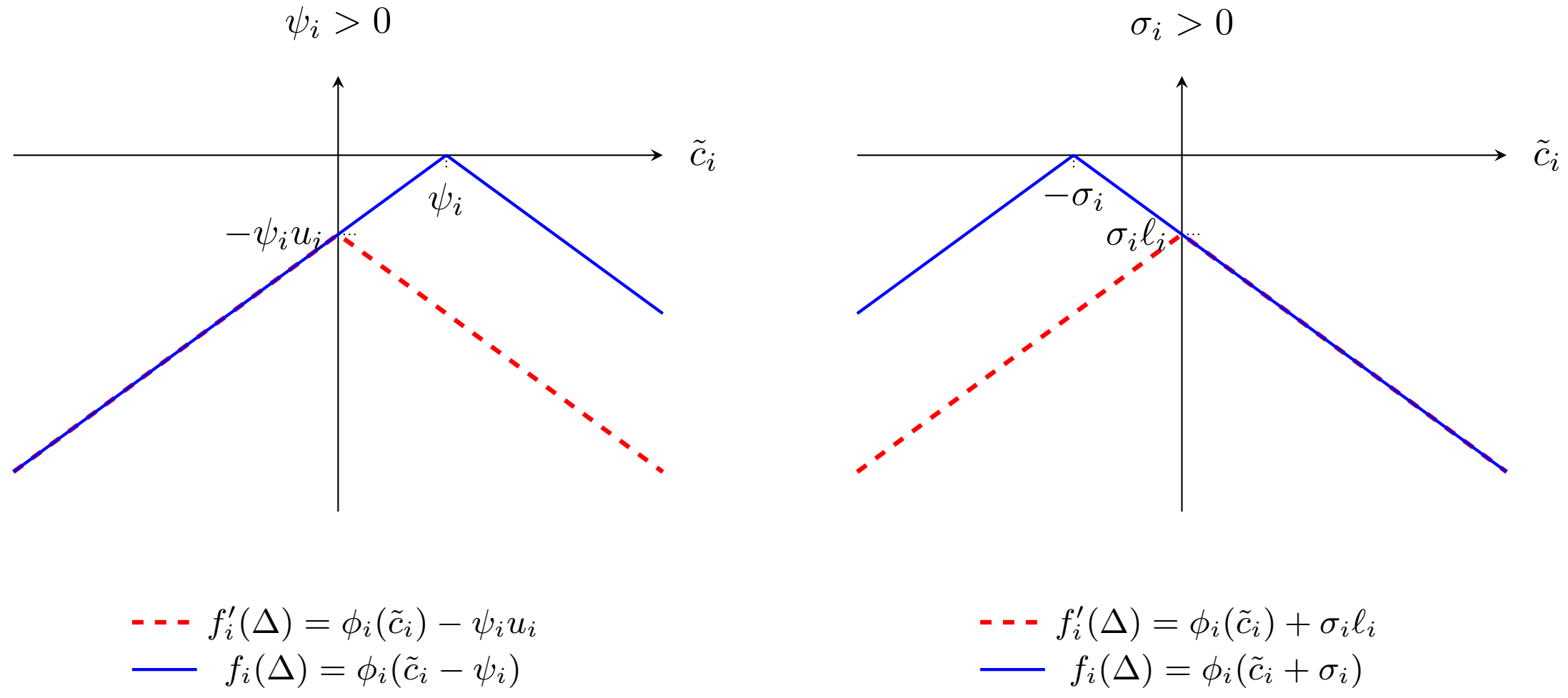
\begin{figure}[H]
\centering

\pgfmathsetmacro{\elli}{-1}   
\pgfmathsetmacro{\ui}{1}      
\pgfmathsetmacro{\psiv}{1}    
\pgfmathsetmacro{\sigv}{1}    

\begin{tikzpicture}
\begin{groupplot}[
    group style={group size=2 by 1, horizontal sep=1.8cm},
    width=0.46\textwidth,
    height=0.34\textwidth,
    axis lines=middle,
    xmin=-3, xmax=3,
    ymin=-4.5, ymax=1.0,
    xtick=\empty,
    ytick=\empty,
    xlabel={$\tilde{c}_i$},
    legend style={draw=none, fill=none, font=\small, at={(0.5,-0.18)}, anchor=north},
    every axis x label/.style={at={(ticklabel* cs:1.02)}, anchor=west},
    every axis y label/.style={at={(ticklabel* cs:1.02)}, anchor=south},
    clip=false,
]

\nextgroupplot[
title={$\psi_i>0$},
]

\addplot[
    red, dashed, very thick, domain=-3:3, samples=401
]
{(x >= 0 ? \elli*x : \ui*x) - \psiv*\ui};
\addlegendentry{$\dual'_i(\pert)=\phi_i(\tilde{c}_i)-\psi_i u_i$}

\addplot[
    blue, thick, domain=-3:3, samples=401
]
{(x-\psiv >= 0 ? \elli*(x-\psiv) : \ui*(x-\psiv))};
\addlegendentry{$\dual_i(\pert)=\phi_i(\tilde{c}_i-\psi_i)$}

\draw[densely dotted] (axis cs:\psiv,0) -- (axis cs:\psiv,-0.18);
\node[below] at (axis cs:\psiv,-0.18) {$\psi_i$};

\draw[densely dotted] (axis cs:0,{-\psiv*\ui}) -- (axis cs:0.18,{-\psiv*\ui});
\node[left] at (axis cs:0,{-\psiv*\ui}) {$-\psi_i u_i$};

\nextgroupplot[
title={$\sigma_i>0$},
]

\addplot[
    red, dashed, very thick, domain=-3:3, samples=401
]
{(x >= 0 ? \elli*x : \ui*x) + \sigv*\elli};
\addlegendentry{$\dual'_i(\pert)=\phi_i(\tilde{c}_i)+\sigma_i \ell_i$}

\addplot[
    blue, thick, domain=-3:3, samples=401
]
{(x+\sigv >= 0 ? \elli*(x+\sigv) : \ui*(x+\sigv))};
\addlegendentry{$\dual_i(\pert)=\phi_i(\tilde{c}_i+\sigma_i)$}

\draw[densely dotted] (axis cs:{-\sigv},0) -- (axis cs:{-\sigv},-0.18);
\node[below] at (axis cs:{-\sigv},-0.18) {$-\sigma_i$};

\draw[densely dotted] (axis cs:0,{\sigv*\elli}) -- (axis cs:0.18,{\sigv*\elli});
\node[left] at (axis cs:0,{\sigv*\elli}) {$\sigma_i \ell_i$};

\end{groupplot}
\end{tikzpicture}

\caption{Coordinatewise comparison between the proposed generalized dual $\dual(\pert)$ and the alternative dual $\dual'(\pert)$. Here $\tilde{c}_i=(\pert[c])_i-(\dualvar_{\nominal}^*)^\top\pert[A]^i$, $\ell_i = -1$ and $u_i = 1$.}
\label{fig:alternative.vs.proposed.dual}
\end{figure}

We compare the proposed bounded Lagrangian dual with the alternative that also dualizes the variable bounds. Let $\sigma\geq0$ and $\psi\geq0$ be the multipliers of $x\geq\ell$ and $-x\geq-u$, respectively. Retaining the box while dualizing these constraints gives
\begin{align*}
\mathcal L'_{\pert}(\lambda,\sigma,\psi)
&:=\min_{\ell\leq x\leq u}\big\{(c+\pert[c])^\top x
-\lambda^\top((A+\pert[A])x-(b+\pert[b]))\\
&\hspace{38mm}-\sigma^\top(x-\ell)-\psi^\top(u-x)\big\}.
\end{align*}
At a nominal optimal dual solution $(\lambda_0^*,\sigma_0^*,\psi_0^*)$, dual feasibility gives
\[
c-A^\top\lambda_0^*=\sigma_0^*-\psi_0^*.
\]
Define the perturbation-only reduced-cost change and the box minimum by
\[
\widetilde c_i(\pert):=[\pert[c]]_i-(\lambda_0^*)^\top[\pert[A]]^i,
\qquad
\phi_i(t):=\min\{\ell_i t,u_i t\}.
\]
Here $\widetilde c_i(\pert)$ is the change in reduced cost, whereas $\bar c_i^v(\pert)$ in Section~\ref{subsec:gpu.node.eval} denotes the full perturbed reduced cost at node $v$. Using $t^-:=\min\{t,0\}$, we have $\phi_i(t)=t^+\ell_i+t^-u_i$. Consequently,
\begin{align*}
f'_{LP}(\pert)
&=\sum_{i=1}^n\big(\phi_i(\widetilde c_i(\pert))
 +(\sigma_0^*)_i\ell_i-(\psi_0^*)_i u_i\big)
 +(\lambda_0^*)^\top(b+\pert[b]),\\
f_{LP}(\pert)
&=\sum_{i=1}^n\phi_i\big(\widetilde c_i(\pert)+(\sigma_0^*)_i-(\psi_0^*)_i\big)
 +(\lambda_0^*)^\top(b+\pert[b]).
\end{align*}
For each coordinate, minimizing a sum over the same interval gives
\begin{align*}
\phi_i\big(\widetilde c_i(\pert)+(\sigma_0^*)_i-(\psi_0^*)_i\big)
&=\min_{\ell_i\leq x_i\leq u_i}
 \big\{\widetilde c_i(\pert)x_i+(\sigma_0^*)_i x_i-(\psi_0^*)_i x_i\big\}\\
&\geq\phi_i(\widetilde c_i(\pert))
 +(\sigma_0^*)_i\ell_i-(\psi_0^*)_i u_i.
\end{align*}
Summing proves $f'_{LP}(\pert)\leq f_{LP}(\pert)$. The proposed dual can therefore be stronger when the perturbation favors moving a variable away from a bound with a positive nominal multiplier. Figure~\ref{fig:alternative.vs.proposed.dual} illustrates the cases of a positive upper-bound multiplier $\psi_i$ and a positive lower-bound multiplier $\sigma_i$, with subscripts denoting nominal optimal multipliers for brevity.

\section{Selected instance details}~\label{appendix:instance.selection}

\begin{table}[htbp]
  \centering
  \caption{Selected nominal instances for combinatorial auctions (CA).}
  \label{tab:CA}
  \small
  \setlength{\tabcolsep}{4pt}
  \begin{tabular}{@{}lrrrrr@{}}
    \toprule
    Instance name & Variables & Binary variables & Constraints & Tree nodes & Depth \\
    \midrule
    train-instance\_703 & 1500 & 1500 & 579 & 3145 & 21 \\
    train-instance\_821 & 1500 & 1500 & 589 & 3231 & 21 \\
    train-instance\_1157 & 1500 & 1500 & 573 & 4259 & 27 \\
    train-instance\_907 & 1500 & 1500 & 579 & 5595 & 27 \\
    train-instance\_443 & 1500 & 1500 & 575 & 7939 & 29 \\
    train-instance\_696 & 1500 & 1500 & 578 & 15009 & 26 \\
    train-instance\_1110 & 1500 & 1500 & 573 & 15077 & 33 \\
    train-instance\_510 & 1500 & 1500 & 586 & 16445 & 33 \\
    val-instance\_893 & 1500 & 1500 & 582 & 17537 & 28 \\
    train-instance\_603 & 1500 & 1500 & 577 & 21323 & 31 \\
    train-instance\_13 & 1500 & 1500 & 589 & 39147 & 31 \\
    train-instance\_262 & 1500 & 1500 & 588 & 39263 & 26 \\
    train-instance\_208 & 1500 & 1500 & 583 & 44965 & 34 \\
    val-instance\_318 & 1500 & 1500 & 580 & 51359 & 33 \\
    train-instance\_695 & 1500 & 1500 & 571 & 69615 & 39 \\
    \bottomrule
  \end{tabular}
\end{table}

\begin{table}[htbp]
  \centering
  \caption{Selected nominal instances for the generalized independent set problem (GISP).}
  \label{tab:GISP}
  \small
  \setlength{\tabcolsep}{4pt}
  \begin{tabular}{@{}lrrrrr@{}}
    \toprule
    Instance name & Variables & Binary variables & Constraints & Tree nodes & Depth \\
    \midrule
    train-newgisp\_erdos\_renyi\_483 & 956 & 956 & 3276 & 296727 & 106 \\
    test-gisp\_erdos\_renyi\_48 & 989 & 989 & 3362 & 331755 & 109 \\
    train-newgisp\_erdos\_renyi\_705 & 941 & 941 & 3327 & 344367 & 108 \\
    train-newgisp\_erdos\_renyi\_367 & 997 & 997 & 3340 & 369523 & 107 \\
    val-gisp\_erdos\_renyi\_22 & 967 & 967 & 3389 & 384627 & 109 \\
    train-newgisp\_erdos\_renyi\_225 & 945 & 945 & 3318 & 391795 & 107 \\
    train-newgisp\_erdos\_renyi\_262 & 1026 & 1026 & 3351 & 397559 & 109 \\
    train-newgisp\_erdos\_renyi\_768 & 967 & 967 & 3386 & 401907 & 109 \\
    train-newgisp\_erdos\_renyi\_767 & 995 & 995 & 3442 & 425239 & 109 \\
    val-gisp\_erdos\_renyi\_25 & 974 & 974 & 3405 & 433733 & 109 \\
    train-newgisp\_erdos\_renyi\_258 & 937 & 937 & 3332 & 449795 & 109 \\
    train-newgisp\_erdos\_renyi\_243 & 1006 & 1006 & 3403 & 481481 & 109 \\
    train-newgisp\_erdos\_renyi\_801 & 984 & 984 & 3356 & 531445 & 109 \\
    train-newgisp\_erdos\_renyi\_405 & 1037 & 1037 & 3393 & 580055 & 109 \\
    val-gisp\_erdos\_renyi\_99 & 995 & 995 & 3331 & 594031 & 109 \\
    \bottomrule
  \end{tabular}
\end{table}

\begin{table}[htbp]
  \centering
  \caption{Selected nominal instances for the capacitated facility location problem (CFLP).}
  \label{tab:CFLP}
  \small
  \setlength{\tabcolsep}{4pt}
  \begin{tabular}{@{}lrrrrr@{}}
    \toprule
    Instance name & Variables & Binary variables & Constraints & Tree nodes & Depth \\
    \midrule
    val-instance\_262 & 20200 & 200 & 20301 & 101 & 22 \\
    train-instance\_513 & 20200 & 200 & 20301 & 231 & 44 \\
    train-instance\_352 & 20200 & 200 & 20301 & 681 & 114 \\
    test-instance\_68 & 20200 & 200 & 20301 & 923 & 129 \\
    train-instance\_546 & 20200 & 200 & 20301 & 975 & 129 \\
    train-instance\_308 & 20200 & 200 & 20301 & 1563 & 27 \\
    val-instance\_305 & 20200 & 200 & 20301 & 1821 & 116 \\
    train-instance\_995 & 20200 & 200 & 20301 & 3205 & 175 \\
    train-instance\_949 & 20200 & 200 & 20301 & 4203 & 41 \\
    val-instance\_526 & 20200 & 200 & 20301 & 4571 & 104 \\
    test-instance\_846 & 20200 & 200 & 20301 & 8109 & 92 \\
    train-instance\_661 & 20200 & 200 & 20301 & 10077 & 62 \\
    train-instance\_981 & 20200 & 200 & 20301 & 22737 & 165 \\
    val-instance\_952 & 20200 & 200 & 20301 & 56173 & 96 \\
    val-instance\_883 & 20200 & 200 & 20301 & 68899 & 76 \\
    \bottomrule
  \end{tabular}
\end{table}

\begin{table}[htbp]
  \centering
  \caption{Selected nominal instances for unit commitment (UC). Here $G$ is the number of generators and $T$ is the number of time periods.}
  \label{tab:UC}
  \small
  \setlength{\tabcolsep}{4pt}
  \begin{tabular}{@{}lrrrrrr@{}}
    \toprule
    Instance name & $(G,T)$ & Variables & Binary variables & Constraints & Tree nodes & Depth \\
    \midrule
    instance\_3 & (25,24) & 1800 & 1200 & 2821 & 3681 & 153 \\
    instance\_11 & (27,24) & 1944 & 1296 & 3017 & 15053 & 93 \\
    instance\_10 & (26,24) & 1872 & 1248 & 2906 & 27925 & 72 \\
    instance\_2 & (25,24) & 1800 & 1200 & 2833 & 30137 & 118 \\
    instance\_64 & (37,24) & 2664 & 1776 & 4171 & 31729 & 178 \\
    instance\_20 & (28,24) & 2016 & 1344 & 3172 & 35463 & 92 \\
    instance\_21 & (29,24) & 2088 & 1392 & 3286 & 62823 & 126 \\
    instance\_23 & (29,24) & 2088 & 1392 & 3273 & 81545 & 62 \\
    instance\_18 & (28,24) & 2016 & 1344 & 3160 & 86265 & 188 \\
    instance\_35 & (31,24) & 2232 & 1488 & 3482 & 89483 & 117 \\
    instance\_32 & (31,24) & 2232 & 1488 & 3503 & 117577 & 112 \\
    instance\_27 & (30,24) & 2160 & 1440 & 3373 & 154973 & 79 \\
    instance\_33 & (31,24) & 2232 & 1488 & 3504 & 211935 & 132 \\
    instance\_24 & (29,24) & 2088 & 1392 & 3273 & 259009 & 147 \\
    instance\_37 & (32,24) & 2304 & 1536 & 3612 & 261329 & 187 \\
    \bottomrule
  \end{tabular}
\end{table}

\end{document}